\documentclass[a4paper,regno,11pt]{amsart}

\usepackage[margin=1in]{geometry}

\usepackage[english]{babel}
\usepackage{graphicx} %figure, includegraphics
\usepackage{tikz}
\usetikzlibrary{arrows}
\usepackage{amsthm} %theoremstyle
\usepackage{amsmath}
\usepackage{amssymb}
\usepackage{amsthm} %proof
\usepackage{mathrsfs} %mathscr
\usepackage{bm,bbm} %mathbm 
\usepackage{pifont}
\usepackage{enumerate} %(i)
\usepackage[shortlabels]{enumitem}
\usepackage{xcolor}
\usepackage{multirow} %\begin{tabular}
\usepackage{subcaption} %subfigure

\usepackage{hyperref} %hyperlinks

\newcommand{\cmark}{\text{\ding{51}}}
\newcommand{\xmark}{\text{\ding{55}}}

\newcommand{\N}{\mathbb{N}}
\newcommand{\Z}{\mathbb{Z}}

\newcommand{\R}{\mathbb{R}}
\newcommand{\C}{\mathbb{C}}
\newcommand{\T}{\mathbb{T}}
\newcommand{\D}{\mathbb{D}}
\newcommand{\B}{\mathbb{B}}
\newcommand{\HP}{\mathbb{H}}
\newcommand{\Sp}{\mathbb{S}}
\newcommand{\Ss}{\mathbf{S}}
\newcommand{\Sc}{\mathbf{D}}
\newcommand{\Efunc}{\mathscr{E}}
\newcommand{\Tfunc}{\mathscr{T}}
\newcommand{\Sfunc}{\mathscr{S}}

\newcommand{\function}[5]{\begin{eqnarray*}
		#1:#2 & \rightarrow & #3 \\
		   #4 & \mapsto     & #5
\end{eqnarray*}}

\newcommand{\velocity}{\mathbf{u}}
\newcommand{\aux}{\mathbf{m}}
\newcommand{\vect}{\mathbf{v}}
\newcommand{\im}{\bm{i}}
\newcommand{\x}{\bm{x}}
\newcommand{\Mixzone}{\Omega_{\mathrm{mix}}}
\newcommand{\Kconstrain}{\mathcal{K}}
\newcommand{\set}{[-1,1]\times\R^2\times\R^2}
\newcommand{\Domain}{\mathscr{D}}
\newcommand{\expected}[1]{\langle #1\rangle}
\newcommand{\dev}[1]{\bar{#1}}
\newcommand{\Atd}{\At_{\rho}}
\newcommand{\Atv}{\At_{\mu}}
\newcommand{\At}{A}
\newcommand{\Sta}{\mathbf{S}}
\newcommand{\rel}{\mathcal{U}}

\newcommand{\dif}{\,\mathrm{d}}

\newcommand{\Div}{\nabla\cdot}
\newcommand{\Curl}{\nabla^\perp\cdot}
\newcommand{\dist}[2]{\mathrm{dist}(#1,#2)}
\newcommand{\TT}{\varphi}
\newcommand{\LL}{L}
\newcommand{\shift}{T}

\def\Xint#1{\mathchoice
	{\XXint\displaystyle\textstyle{#1}}%
	{\XXint\textstyle\scriptstyle{#1}}%
	{\XXint\scriptstyle\scriptscriptstyle{#1}}%
	{\XXint\scriptscriptstyle\scriptscriptstyle{#1}}%
	\!\int}
\def\XXint#1#2#3{{\setbox0=\hbox{$#1{#2#3}{\int}$}
		\vcenter{\hbox{$#2#3$}}\kern-.5\wd0}}

\def\dashint{\Xint-}  %mean integral

\theoremstyle{theorem}
\newtheorem{thm}{Theorem}[section]
\newtheorem{prop}{Proposition}[section]

\newtheorem{lemma}{Lemma}[section]
\theoremstyle{definition}
\newtheorem{defi}{Definition}[section]
\theoremstyle{remark}
\newtheorem{Rem}{Remark}[section]

\usepackage{fancyhdr}
\title[]{H-principle for the 2D incompressible porous media equation with viscosity jump}
\author{Francisco Mengual}
\numberwithin{equation}{section} % label equation with sections
\begin{document}
	
\begin{center}
\begin{minipage}[h]{\textwidth}
	%\centering
	\maketitle
\end{minipage}
\begin{abstract} In this work we extend the results in \cite{LIPM,RIPM} on the 2D IPM system with constant viscosity (Atwood number $\Atv=0$) to the case of viscosity jump ($|\Atv|<1$).
We prove a h-principle whereby (infinitely many) weak solutions in $C_tL_{w^*}^{\infty}$ are recovered via convex integration whenever a subsolution is provided. As a first example, non-trivial weak solutions with compact support in time are obtained. Secondly, we construct mixing solutions to the unstable Muskat problem with initial flat interface.
As a byproduct, we check that the connection, established by Sz\'ekelyhidi for $\Atv=0$ \cite{RIPM}, between the subsolution and the Lagrangian relaxed solution of Otto \cite{O}, holds for $|\Atv|<1$ too.
For different viscosities, we show how a pinch singularity in the relaxation prevents the two fluids from mixing wherever there is neither Rayleigh-Taylor nor vorticity at the interface.
\end{abstract}
\end{center}

%\tableofcontents

\pagenumbering{arabic} 
\setcounter{page}{1}
\section{Introduction and main results}\label{sec:Introduccion}

We deal with the evolution of two incompressible 
fluids with constant densities $\rho^+>\rho^->0$ and viscosities $\mu^+,\mu^->0$ (e.g. water and oil \cite{Mus34}) moving through a 2D porous medium $\Domain$ with constant permeability $\kappa>0$ (or Hele-Shaw cell \cite{ST58}) under the action of gravity $\mathbf{g}=-g\im$, where $\bm{i}=(0,1)$ will also play the roll of the imaginary unit by identifying $\R^2\simeq\C$. Following \cite{O}, we introduce the $\{-1,1\}$-valued variable $\theta(t,\x)$ to indicate whether at time $t\in\R_+$ the pores near $\x=(x_1,x_2)\in\Domain$ are filled with \textbf{phase} $-$ or $+$: 
\begin{equation}\label{phase}\tag{IPM0}
a(t,\x):=\frac{a^++a^-}{2}+\frac{a^+-a^-}{2}\theta(t,\x),
\quad a=\rho,\mu.
\end{equation}
This two-phase flow %dynamic 
can be modelled (\cite{Mus37}) by the \textbf{IPM} (Incompressible Porous Media) system: 
\begin{align}
\partial_t\theta+\nabla\cdot(\theta\velocity) & = 0, \label{IPM:1}\tag{IPM1}\\
\Div\velocity & = 0,\label{IPM:2}\tag{IPM2}\\
\tfrac{\mu}{\kappa}\velocity & = -\nabla p+\rho\mathbf{g}, \label{IPM:3}\tag{IPM3}
\end{align}
in $\R_+\times\Domain$. (IPM0-2) reads as the phase distribution $\theta$ (resp. $\rho$ and $\mu$) is advected by the incompressible flow (coupled with the no-flux boundary condition). (IPM3)
is \textbf{Darcy's law}, which relates the velocity field $\velocity$ of the fluid with the forces acting on it.
By renaming the pressure $p$, Darcy's law can be written in terms of the phase $\theta$ as
\begin{equation}\label{IPM:3.3}
\velocity+\Atv\theta\velocity+\Atd\theta\im= -\nabla p
\tag{$\textrm{IPM3}_\At$},
\end{equation}
where $\Atd$, $\Atv$ are the \textbf{Atwood numbers}
$$
\Atd
:=\kappa g\frac{\rho^+-\rho^-}{\mu^++\mu^-}> 0,
\quad\quad\Atv
:=\frac{\mu^+-\mu^-}{\mu^++\mu^-}\in(-1,1).
$$
Since (IPM0-2) is invariant under the scaling $\theta(\alpha t,\x)$, $\alpha\velocity(\alpha t,\x)$, 
by normalizing ($\alpha=\Atd$) and renaming $p$, we may assume w.l.o.g.~that $\Atd=1$.
Thus, from now on we shall abbreviate $\At\equiv\Atv$.
We have added the tag ``$\At$'' to the reference (IPM3)
 to make explicit the dependence on this parameter. Similarly, we shall abbreviate $(\textrm{IPM}_\At)\equiv(\textrm{IPM0-3}_\At)$.\\ 

\textbf{The main results}. 
The phase jump induces Rayleigh-Taylor (RT) and vorticity
at the interface separating both fluids, which becomes unstable when the RT condition fails (cf.~\S\ref{sec:link}).
In such a case, %it can be observed that 
the two fluids can start to mix on a mesoscopic scale (see e.g.~\cite[pp.~261-267]{WM76} and \cite{Hom87}).
Although unstable configurations in Hydrodynamics are very difficult to model, 
De Lellis-Sz\'ekelyhidi's version of convex integration (\cite{Eulerinclusion,Onadmissibility}) have successfully
describe several examples as the RT instability for ($\textrm{IPM}_0$) \cite{Mixing,Degraded,Piecewise,RIPM}, and the Kelvin-Helmholtz \cite{Sze11} and RT \cite{GKS20} instabilities for the Incompressible Euler equations.\\
\indent In this work we investigate the scope of this view point to the RT instability for ($\textrm{IPM}_\At$) in the case of different viscosities (or mobilities in \cite{O}, cf.~\S\ref{sec:Otto}) which is a recurrent theme in the applied literature.
In short terms, the approach seems to work at least for flat interfaces, but the relaxation presents some unexpected singularities which makes
the project challenging.\\
\indent Before going any further let us present the problem discussed, summarize the main results of this work as well as the technical difficulties, and go back at the end of the introduction with a new link between the mixing regime and the relaxation. Firstly, we present two theorems regarding weak solutions to $(\mathrm{IPM}_\At)$ for any $|\At|<1$ (cf.~Def.~\ref{defi:weaksol}). The first one exhibits lack of uniqueness in the class $C_tL_{w^*}^\infty$.
\begin{thm}\label{thm:Lack} Let $|\At|<1$, $T>0$ and $\Domain=\R^2$ or $\T^2$. There exist infinitely many weak solutions  $(\theta,\velocity)\in  C(\R_+;L_{w*}^{\infty}(\Domain))$ to $(\mathrm{IPM}_\At)$ with $|\theta|=1$ on $(0,T)\times\Domain$ and $\theta=0$ outside.
\end{thm}
Thus, $(\mathrm{IPM}_\At)$ admits non-trivial weak solutions with compact support in time. Opposite to these unphysical solutions, we construct admissible weak solutions to the \textbf{unstable Muskat problem} with initial flat interface. This is $(\mathrm{IPM}_\At)$ starting from the unstable planar phase
\begin{equation}\label{flat}
\theta_0(\x)=\left\lbrace
\begin{array}{rc}
+1, & x_2>0, \\
-1, & x_2<0.
\end{array}\right.
\end{equation}
Similarly to \cite{Mixing,Degraded,Piecewise,RIPM}, we show that these weak solutions start to mix inside a mixing zone $\Mixzone$ which grows linearly in time around $x_2=0$, and that they look macroscopically
almost like the coarse-grained phase, denoted in this paper by $\Theta_{\At}$ (cf.~\eqref{entropy}), introduced by Otto in \cite{O}. For this reason, we shall call them ``$\Theta_{\At}$-mixing solutions'' (cf.~Def.~\ref{def:DMS} and Fig.~\ref{fig:theta0}-\ref{fig:after}).

\begin{thm}\label{thm:DMS} Let $|\At|<1$ and $\Domain=\R^2$ or $(-1,1)^2$. There exist infinitely many $\,\Theta_{\At}$-mixing solutions $(\theta,\velocity)\in  C(\R_+;L_{w*}^{\infty}(\Domain))$ to $(\mathrm{IPM}_\At)$ starting from the unstable planar phase \eqref{flat}.
\end{thm}
While the weak solutions from Theorem~\ref{thm:Lack} 
can not attain the initial datum $\theta_0=0$ in the strong sense, the ones from Theorem \ref{thm:DMS} satisfy $\theta\in C(\R_+;L_{\mathrm{loc}}^p(\Domain))$ for all $1<p<\infty$. Moreover, they are forced to have finite mixing speed (cf.~Prop.~\ref{maxmix}).\\

These theorems are deduced from a more general h-principle (cf.~Thm.~\ref{HP}). In brief, this reads as weak solutions to $(\mathrm{IPM}_\At)$ can be recovered via convex integration whenever a subsolution is provided (cf.~\S\ref{sec:HP}). This subsolution (cf.~Def.~\ref{defi:weaksol}) is a weak solution to a linearised version $(\textbf{T}_\At)$ of $(\textrm{IPM}_\At)$, taking values in a relaxed set $\bar{\rel}_{\At}$ of the corresponding constitutive set ($\Kconstrain$), namely $\rel_{\At}$ is an open set satisfying a perturbation property w.r.t.~$(\textbf{T}_\At,\Kconstrain)$.\\ 
\indent The proof of the h-principle is classical (\cite{Degraded,Onadmissibility,RIPM}) but difficulties arise as the parameter $\At$, which originally looks innocent, turns the relation between the components of the subsolution less explicit, which ends up hampering considerably the proof of the hypothesis $\textrm{(H1)-(H3)}_p$ required therein (cf.~\S\ref{sec:HPproof}). 
For instance, the $L^p$-boundedness property $\textrm{(H3)}_p$ becomes non-trivial for $0<|\At|<1$  (cf.~Lemmas~\ref{L2bound} and \ref{Linf}). 
A more delicate issue is the relaxation $\bar{\rel}_\At$. We take $\bar{\rel}_\At=\Kconstrain^{lc,\Lambda_\At}$ $\equiv$ $\Lambda_\At$-lamination hull of $\Kconstrain$, which we compute explicitly (cf.~\eqref{rel} and \S\ref{sec:Relaxation}). However, since it is not obvious that such $\bar{\rel}_{\At}$ is closed under weak*-convergence (not even that $\bar{\rel}_{\At}$ is equal to the functional $\Lambda_{\At}$-convex hull of $\Kconstrain$) we refine the Baire category argument to adapt the proof of the h-principle we follow \cite{Degraded,Onadmissibility} to our situation (cf.~Rem.~\ref{Rem:stability}).\\ 
\indent While the relaxation $\bar{\rel}_0$ only narrows at $\Kconstrain$, for different viscosities $\bar{\rel}_{\At}$ develops a pinch singularity far away from $\Kconstrain$. Up to our knowledge, this kind of singularity outside the constitutive set $\Kconstrain$ does not appear in other examples in Hydrodynamics. 
This necessarily complicates the existence of long $\Lambda_{\At}$-segments as the perturbation property (H2) requires. To our surprise, they do exist even if $\bar{\rel}_{\At}$ is very narrow far away from $\Kconstrain$. Remarkably, the use of Complex Analysis becomes very helpful, reducing considerably some tedious computations and providing a nice geometric interpretation in terms of the automorphisms of the unit disc (cf.~Rem.~\ref{Rem:geometry}).\\
\indent In order to find bounded velocities, Sz\'ekelyhidi computed cleverly the relaxation of some $\Kconstrain_M\Subset\Kconstrain$ for $\At=0$. 
In the case of viscosity jump the parameter $\At$ introduces an asymmetry that makes less clear what restriction of $\Kconstrain$ may return a simple relaxation (cf.~Rem.~\ref{Rem:asymmetry}). The way of arguing is somewhat original as first we guess (inspired by an identity in \cite{RIPM}) a shape for $\bar{\rel}_{\At,M}$, and then find $\Kconstrain_{\At,M}\Subset\Kconstrain$ satisfying  $(\Kconstrain_{\At,M})^{lc,\Lambda_{\At}}=\bar{\rel}_{\At,M}$. \\
\indent The proof of the perturbation property (H2) for $\rel_{\At,M}$ presents some added difficulties compared to $\At=0$ (cf.~Lemma~\ref{lemma:ZMseg}). 
The main obstacle is that one of the inequalities bounding $\rel_{\At,M}$, which is just a restriction on $\velocity$ for $\At=0$, depends on $\aux$ (relaxation of the non-linear term $\theta\velocity$) for $0<|\At|<1$. 
Geometrically, the projection $\rel_{\At,M}(\theta,\velocity)\equiv\{\aux\in\R^2\,:\,(\theta,\velocity,\aux)\in\rel_{\At,M}\}$, which is given by the intersection of three balls for $\At=0$, is also restricted by a half-plane for $0<|\At|<1$ (cf.~Fig.~\ref{fig:geo}).
This causes that $\rel_{\At,M}(\theta,\velocity)$ collapses as $|\velocity|$ grows, in contrast to the case $\At=0$ (cf.~Fig.~\ref{geo:0}-\ref{geo:MM}). 
Furthermore, the pinch singularity 
becomes further complicated since
the new inequalities defining $\rel_{\At,M}$ can interfere
with it (cf.~Rem.~\ref{Rem:M*}). All this makes the choice of the $\Lambda_{\At}$-segments cumbersome in some of the cases (see e.g.~\eqref{alpha}\eqref{balphasol}).

\subsection{A link between the mixing regime and the relaxation}\label{sec:link}
The aim of this section is to analyse the physical implications of the pinch singularity that arises at $\rel_{\At}$. In a nutshell, it prevents the two fluids from mixing wherever there is neither Rayleigh-Taylor nor vorticity (equiv. $\nabla p$ and $\velocity$ are continuous) at the interface. Let us explain this in more detail.\\
\indent The Muskat problem describes $(\textrm{IPM}_\At)$ under the assumption that there is a time-dependent moveable interface $\mathbf{z}(t)$ separating $\Domain$ in two disjoint open sets $\Omega_\pm(t)$ $\equiv$ region occupied by the  fluid with phase $\pm$ at time $t$.
Let us denote $f^\uparrow$ ($f^\downarrow$) by the limit of $f(\mathbf{z}+\varepsilon\partial_s\mathbf{z}^\perp)$ as $\varepsilon\uparrow 0$ ($\varepsilon\downarrow 0$), and also $[f]:=f^\uparrow-f^\downarrow$ by the jump of $f=\theta$, $\velocity$, $p$ along $\mathbf{z}$.\\  
\indent The Biot-Savart system ($\textrm{IPM2-3}_{\At}$) determines $p$ and $\velocity$ in terms of $\mathbf{z}$ and $[\theta]$. On the one hand, the incompressibility condition ($\textrm{IPM2}$) implies that $\velocity=\nabla^\perp\psi$ for some stream function $\psi$, and so the vorticity $\omega:=\Curl\velocity=\Delta\psi$. On the other hand,
by applying $\Div$ and $\Curl$ on Darcy's law ($\textrm{IPM3}_{\At}$), we deduce that  both $\Delta p$ and $\Delta \psi$ are Dirac measures supported on $\mathbf{z}$
$$\Delta p=\sigma\delta_{\mathbf{z}},
\quad\quad
\Delta\psi=\varpi\delta_{\mathbf{z}},$$
for some scalar functions $\sigma$ $\equiv$ Rayleigh-Taylor and $\varpi$ $\equiv$ vorticity strength. 
Thus, both $p$ and $\psi$ (and so $\velocity$) are recovered from $\sigma$ and $\varpi$ respectively by means of Potential Theory,
namely they are harmonic outside $\mathbf{z}$ and have well-defined traces. Moreover, $p$ and $\psi$ are continuous ($[p]=[\psi]=0$) 
but have discontinuous gradients along $\mathbf{z}$ ($*$ $\equiv$ complex conjugate)
$$[\nabla p] = -\im\frac{\sigma}{\partial_s\mathbf{z}^*},
\quad\quad
[\nabla\psi] = -\im\frac{\varpi}{\partial_s\mathbf{z}^*},$$
and so $[\velocity]=\im[\nabla\psi]$. Observe $\sigma=-[\nabla p]\cdot\partial_s\mathbf{z}^\perp$
and $\varpi=[\velocity]\cdot\partial_s\mathbf{z}$. 
Thus, (the jump along $\mathbf{z}$ of) Darcy's law $(\textrm{IPM3}_{\At})$
reads as
\begin{equation}\label{RTKH}
\varpi+\sigma\im
=-[\theta](\At\breve{\velocity}+\im)^*\partial_s\mathbf{z},
%\quad
%[\velocity]
%+[\theta](\At\breve{\velocity}+\im)=-[\nabla p],
\end{equation}
where $\breve{\velocity}:=\tfrac{1}{2}(\velocity^\uparrow+\velocity^\downarrow)$ is the \textbf{mean velocity} along $\mathbf{z}$. Observe that both $\sigma$ and $\varpi$ vanish if and only if $\At\breve{\velocity}+\im = 0$. As we shall see, these are precisely the states where $\rel_{\At}$ pinches.\\
\indent Finally, (IPM1) turns out to be a free boundary problem, namely $\mathbf{z}$ is driven by the Birkhoff-Rott integrodifferential equations
\begin{equation}\label{Cauchyz}
\partial_t\mathbf{z}=\breve{\velocity}(\mathbf{z})+r\partial_s\mathbf{z},
\quad\quad\mathbf{z}|_{t=0}=\mathbf{z}_0,
\end{equation}
where $r$ represents the re-parametrization freedom, $\breve{\velocity}(\mathbf{z})=\mathscr{B}(\mathbf{z},\varpi(\mathbf{z}))$ with
$$\mathscr{B}(\mathbf{z},\varpi)(t,\alpha)^*
=\frac{1}{2\pi\im}\mathrm{PV}\!\int\frac{\varpi(t,\beta)}{\mathbf{z}(t,\alpha)-\mathbf{z}(t,\beta)}\dif\beta,$$
and, by \eqref{RTKH}, $\varpi(\mathbf{z})$ is given by the (implicit) equation $\varpi(\mathbf{z})=-[\theta](\At\mathscr{B}(\mathbf{z},\varpi(\mathbf{z}))+\im)\cdot\partial_s\mathbf{z}$. Similarly, $\sigma(\mathbf{z})
=[\theta](\At\mathscr{B}(\mathbf{z},\varpi(\mathbf{z}))+\im)\cdot\partial_s\mathbf{z}^\perp$.\\ 
\indent In brief, this Cauchy problem \eqref{Cauchyz} for $\mathbf{z}$ is well-posed provided the Rayleigh-Taylor (also called Saffman-Taylor \cite{ST58}) condition for the Muskat problem, $\sigma>0$, holds (\cite{Amb04,Amb14,CCG11,GGPS18,Mat18,Mat19,SCH04}). The geometric meaning of $\sigma(\mathbf{z})>0$ is not evident since the dependence on $\mathbf{z}$ is highly implicit.
The situation is simpler for equal viscosities ($\At=0$) or flat interfaces ($\velocity=0$) because $[\theta]\partial_s\mathbf{z}_1>0$ just requires the heavier fluid to remain below the lighter. The Muskat problem for $\At=0$ has been widely studied in the literature (see the survey \cite{Gan17} and the references therein).\\
\indent 
When the RT condition fails the free boundary can turn into a growing strip, $\Mixzone$ $\equiv$ mixing zone, where the phases start to mix on a mesoscopic scale.
In the last years
this kind of mixing solutions have been constructed by means of convex integration in the RT unstable regime (\cite{Mixing,Degraded,Piecewise,RIPM}). They are driven by a two-scale dynamic: one dealing with the evolution of the pseudo-interface, which may describe the macroscopic fingering phenomenon, and other dealing with the laminar-turbulent transition region $\Mixzone$ around the pseudo-interface.\\
\indent In \cite{Mixing,Piecewise} the authors discovered that
mixing solutions also exist in the RT stable regime  provided the velocity is discontinuous, i.e.~when $\varpi\neq 0$. Inspired by \cite{Sze11},
we speculate it may describe a turbulence zone of spiral vortices, usually observed in the Kelvin-Helmholtz instability. 
We remark in passing that, since there are initial data $\mathbf{z}_0$ for which both \eqref{Cauchyz} is solvable and mixing solutions exist, a main unsolved question is to identify a selection criterion among them which leads to a unique physical solution.\\
\indent In short, it seems that
the mixing phenomenon may be triggered at least by two mechanisms:
$\sigma < 0$ or $\varpi\neq 0.$
By \eqref{RTKH}, one of these 
is awake at some 
point of the interface 
$\mathbf{z}(s)$ if $$-[\theta](\At\breve{\velocity}(\mathbf{z}(s))+\im)^*\partial_s\mathbf{z}(s)\in\mathcal{M},$$ 
where $\mathcal{M}:=\R^2\setminus\mathcal{L}$ $\equiv$ \textbf{mixing regime} and $\mathcal{L}:=\{\varpi+\sigma\im\,:\,\sigma\geq 0=\varpi\}$. Conversely, the open half-line $\mathcal{L}^\circ=\{\varpi+\sigma\im\,:\,\sigma> 0=\varpi\}$ classifies the points where the interface is RT stable and there is not vorticity. Remarkably, we have found that the relaxation $\rel_{\At}$ (for different viscosities) excludes $\partial\mathcal{L}=\{0\}$: a pinch a singularity arises at $\At\breve{\velocity}+\im=0$ (cf.~\eqref{rel}) representing the points $\mathbf{z}(s)$ where $\sigma=0=\varpi$. In other words, this relaxation approach prevents the two fluids from mixing wherever both $\nabla p$ and $\velocity$ are continuous.\\

\textbf{Organization of the paper}. 
We start Section \ref{sec:HP} recalling briefly the background of the problem.
After this, we present the h-principle from which Theorems \ref{thm:Lack}-\ref{thm:DMS} are deduced. The proof of this h-principle appears in Section \ref{sec:HPproof}. In Section~\ref{sec:Relaxation} we compute $\bar{\rel}_{\At}$, $\bar{\rel}_{\At,M}$ and show some of their properties. 
With the aim of figuring out how these $\Theta_{\At}$-mixing solutions may look like, we introduce a toy random walk in Appendix~\ref{sec:RW}  (Fig.~\ref{fig:theta0}-\ref{fig:after}). Finally, we recall in Appendix~\ref{sec:Otto} some properties of $\Theta_{\At}$ as well as the transition to the stable planar phase in the confined domain $\Domain=(-1,1)^2$.

\section{H-principle for $(\textrm{IPM}_{\At})$}\label{sec:HP}

We start this section with a brief explanation of the strategy we shall follow, the convex integration method, to help better understand the main results of this work. This method was introduced in Hydrodynamics by De Lellis and Sz\'ekelyhidi in \cite{Eulerinclusion} for the incompressible Euler equations (IE) (see e.g.~\cite{Gromov} for the background in Differential Geometry and \cite{MS03} in PDEs and Calculus
of Variations).\\

Following \cite{LIPM,RIPM}, we introduce a new variable $\aux$ to encode the non-linear term $\theta\velocity$. Thus, if we denote $z=(\theta,\velocity,\aux)
\in[-1,1]\times\R^2\times\R^2$, 
this two-phase flow can be interpreted as a \textbf{differential inclusion} $(\textbf{T}_{\At},\Kconstrain)$ in the spirit of Tartar (\cite{Tartar,Tartar2}) as 
\begin{align}
\Div\textbf{T}_\At(z)=0 \label{DI1}\tag{$\textbf{T}_\At$},\\
z\quad\Kconstrain\textrm{-valued} \label{DI2}\tag{$\Kconstrain$},
\end{align}
in $\R_+\times\Domain$, that is, a linear differential system $(\textbf{T}_\At)$ coupled with a non-linear pointwise constraint \eqref{DI2},
where $\textbf{T}_\At:\R^5\rightarrow\R^{3\times 3}$ is the (injective) linear map
\begin{equation}
\textbf{T}_\At(z):=\left(\begin{array}{ccc}
\theta & \aux_1 & \aux_1 \\
0 & \velocity_1 & \velocity_2 \\
0 & \velocity_2+\At\aux_2+\theta & -\velocity_1-\At\aux_1
\end{array}\right),
\end{equation}
and $\Kconstrain$ is the \textbf{constitutive set}
\begin{equation}\label{K}
\Kconstrain:=\{z\in\set\,:\,|\theta|=1,\,
\aux=\theta\velocity\}.
\end{equation}
Notice that $(\textbf{T}_{\At},\Kconstrain)$ is more demanding than $(\textrm{IPM}_\At)$ because this does not require $|\theta|=1$.\\

Roughly speaking, if an (hypothetical) solution $z$ to $(\textbf{T}_{\At},\Kconstrain)$ is averaged somehow, call the result $\breve{z}$, then $\breve{z}$ solves $(\textbf{T}_{\At},\breve{\Kconstrain}_\At)$ for some set $\breve{\Kconstrain}_\At$.  
It is natural to assume that the fluctuation $z'=z-\breve{z}$ is a highly oscillatory solution (in $\Mixzone$) to $(\textbf{T}_\At)$, thus $z'$  may look (locally) like a \textbf{plane wave} $\bar{z}h(k\xi\cdot(t,\x))$ for some $\bar{z}\in\R^5$, $\xi\in\R\times\Sp^1$, $h\in C^1(\T)$ with $\int h=0$ and $k\gg 1$. The set of directions $\bar{z}$ for which there is a plane wave solving $(\textbf{T}_\At)$ is the \textbf{wave cone} of $(\textbf{T}_\At)$
\begin{equation}\label{WC}
\Lambda_\At
:=\{\bar{z}\in\R^5\,:\,\exists\xi\in\R\times\Sp^1\textrm{ so that }\textbf{T}_\At(\bar{z})\xi=0\}.
\end{equation}
All this suggests that the optimal choice of $\breve{\Kconstrain}_\At$ is $\Kconstrain^{\Lambda_\At}$ $\equiv$ $\Lambda_\At$-convex hull of $\Kconstrain$ (\cite[Def.~4.3]{Rigidity}).
However, when the explicit computation of $\Kconstrain^{\Lambda_\At}$ is unattainable due to the high complexity and dimensionality, it is more practical to consider a simpler but still large enough subset $\breve{\Kconstrain}_\At$ of $\Kconstrain^{\Lambda_\At}$ (see \cite{LIPM,ActiveScalar} and also \cite[\S4]{hprincipleFD}). 
When these correcting terms $z'$ can be constructed and the set $\breve{\Kconstrain}_\At$ satisfies some geometric and functional properties (cf.~\S\ref{sec:HP}) the convex integration method yields a \textbf{homotopy-principle} \cite[\S5]{RIPM} whereby the problem of finding solutions is reduced to find a subsolution, a solution $\breve{z}$ to $(\textbf{T}_\At,\breve{\Kconstrain}_\At)$.
Schematically,
\begin{equation}\label{HPscheme}
\begin{array}{ccccc}
& (\textbf{T}_\At,\Kconstrain) & \overset{\textrm{relaxation}}{\longrightarrow} & (\textbf{T}_\At,\breve{\Kconstrain}_\At) & \\
& \rotatebox[origin=c]{-90}{$\displaystyle\longrightarrow$} & \textrm{h-principle} & \rotatebox[origin=c]{-90}{$\displaystyle\longrightarrow$} & \\
\textrm{solution} & z & \underset{\substack{\textrm{convex} \\ \textrm{integration}}}{\longleftarrow} & \breve{z} & \textrm{subsolution}
\end{array}
\end{equation}

These ideas have been implemented successfully for $\mu^+=\mu^-$ (\cite{Mixing,Degraded,LIPM,Piecewise,RIPM}) but not for $\mu^+\neq\mu^-$. 
Let us recall the previous results for $\At=0$ we want to generalize for $|\At|<1$.\\

\textbf{Brief overview of the case $\At=0$}. In \cite{LIPM}, C\'ordoba, Faraco and Gancedo discovered that the convex integration method developed in \cite{Eulerinclusion} for (IE) could be adapted to prove lack of uniqueness in $L^\infty(\R_+\times\T^2)$
for $(\textrm{IPM}_0)$. In addition, they noticed that, in contrast to \cite{Eulerinclusion},  $\Kconstrain^{\Lambda_0}$ does not agree with $\Kconstrain^{\mathrm{co}}$. To overcome this extra difficulty the authors resorted to the theory of laminates. 
Remarkably, this result was generalized for a class of active scalar equations by Shvydkoy in \cite{ActiveScalar} (see \cite{IV15} for improvements of the
regularity).\\ 
\indent Later in \cite{RIPM} Sz\'ekelyhidi computed explicitly $\Kconstrain^{\Lambda_0}= \bar{\rel}_0$, with $\rel_{0}$ the open set of states $z=(\theta,\velocity,\aux)\in\set$ satisfying
\begin{equation}\label{rel0}
|2(\aux-\theta\velocity)+(1-\theta^2)\im|<(1-\theta^2),
\end{equation}
thus providing a h-principle \eqref{HPscheme} for $\breve{\Kconstrain}_0=\Kconstrain^{\Lambda_0}$ (see \cite{Kno12} for a generalization in a class of active scalar equations).
Another advantage of this computation is that
it allows to identify compatible boundary and initial conditions in order to obtain admissible solutions, opposite to those paradoxical examples with compact support in time. As a promising application in evolution of microstructures, Sz\'ekelyhidi constructed  weak solutions in $L^\infty(\R_+\times(-1,1)^2)$ to the unstable Muskat problem with initial flat interface $\mathbf{z}_0(s)=(s,0)$.
Remarkably, he observed that the subsolution $\breve{\theta}_\alpha$ (for any $0<\alpha<1$, being $c=2\alpha$ the rate of expansion of the mixing zone) 
that naturally arises in this scenario is closely related to the relaxation introduced in \cite{O} (see also \cite{GO12,O2}). In this paper Otto dealt with the general case $|\At|<1$. Since this is the motivation of this work, we have thought appropriate to sketch briefly this approach in Appendix \ref{sec:RW}.\\ 
\indent In short, after introducing a Lagrangian relaxation of $(\textrm{IPM}_\At)$, Otto obtained a unique (relaxed) solution (cf.~\S\ref{sec:RW}-\ref{sec:Otto})
\begin{equation}\label{entropy}
\Theta_{\At}(t,\x)
=\left\lbrace
\begin{array}{cl}
+1, & \hspace{1.35cm}x_2>c_{\At}^+ t,\\[0.1cm]
\frac{x_2+\At t}{ t+\At  x_2+\sqrt{(1-\At^2) t(t+\At x_2)}}, & -c_{\At}^-t< x_2<c_{\At}^+ t,\\[0.1cm]
-1, &  -c_{\At}^-t>x_2,
\end{array}\right.
\quad\textrm{where}\quad
c_{\At}^{\pm}=\frac{2}{1\mp\At},
\end{equation} 
which aims to capture the macroscopic properties of (exact) solutions to $(\textrm{IPM}_{\At})$, thus giving a prediction of the actual shape and evolution of the mixing profile. This $\Theta_{\At}$ is indeed the (unique) entropy solution (\cite[(3.72)]{O})
of the conservation law (or Burgers type equation)
\begin{equation}\label{conservationlaw}
\partial_t\Theta=\partial_{x_2}\left(\frac{1-\Theta^2}{1-\Theta\At}\right),\quad\quad
\Theta|_{t=0}=\theta_0.
\end{equation}
\indent The link between the approaches of Sz\'ekelyhidi and Otto for $\At=0$ is given by 
\begin{equation}\label{SO}
\breve{\theta}_{\alpha}(t)=\Theta(\alpha t),
\quad t\in\R_+,
\end{equation}
(for any $0<\alpha<1$) where $\Theta\equiv\Theta_0$. The interpretation given in \cite{RIPM} of \eqref{SO} is that, although weak solutions are clearly not unique due to the symmetry breakdown, the uniqueness result of Otto can be understood as selecting the subsolution with maximal mixing zone (cf.~Prop.~\ref{maxmix}).\\ 
\indent At this point we remark that a natural question that arises here is if \eqref{SO} defines a subsolution in the general case $|\At|<1$. As we shall see in Theorem~\ref{thm:DMS2}, this is the case.\\
\indent Continuing the overview of the case $\At=0$, Castro, C\'ordoba and Faraco \cite{Mixing} applied this h-principle to construct weak solutions to the unstable Muskat problem for non-flat interfaces $\mathbf{z}_0(s)=(s,f_0(s))$ with $f_0\in H^5(\R)$, by taking the subsolution as $\breve{\theta}_\alpha(t,\x)
=\Theta(\alpha t,\x-f(t,x_1)\im)$ with $f$ a suitable evolution of $f_0$. Moreover, they showed that these solutions indeed mix inside the mixing zone, thus justifying the name ``mixing solution''. In \cite{Piecewise} F\"{o}rster and Sz\'ekelyhidi obtained a similar result for $f_0\in C_*^{3,\gamma}(\R)$ with a simpler proof by taking piecewise constant subsolutions approaching the linear profile of $\Theta$ adapted to $f_0$.\\
\indent Recently, the h-principle presented in \cite{Onadmissibility} was adapted in \cite{Degraded} to measure, in terms of weak*-continuous quantities, the proximity of the weak solutions coming from the convex integration scheme to the subsolution $\breve{z}$, thus selecting those which retain more information from $\breve{z}$, thereby emphasizing the fact that the subsolution aims to be the macroscopic solution (cf.~Rem.~\ref{diagram}). For this reason, the authors called them ``degraded mixing solutions'' (here $\Theta_{0}$-mixing solutions).\\

\textbf{Our extension to the case $|\At|<1$}. 
With the aim of generalizing these results, we follow \cite{Degraded,RIPM} to prove a h-principle for the system $(\mathrm{IPM}_\At)$, which additionally provides weak solutions in the stronger class $C_tL_{w*}^{\infty}$. In order to prove it we need to check three hypothesis. The first one (H1) is the existence of localized plane waves of $(\textbf{T}_\At)$, which is checked similarly to \cite{LIPM,RIPM}.\\ 
The second and more delicate part of this work is to compute a large enough set $\breve{\Kconstrain}_\At$ satisfying the perturbation property (H2). This is the $\Lambda_\At$-lamination hull of $\Kconstrain$,  $\Kconstrain^{lc,\Lambda_\At}=\bar{\rel}_{\At}$ with $\rel_{\At}$ the open set of states $z=(\theta,\velocity,\aux)\in\set$ satisfying
\begin{equation}\label{rel}
|2(1-\theta\At)(\aux-\theta\velocity)+(1-\theta^2)(\At\velocity+\im)|<(1-\theta^2)|\At\velocity+\im|.
\end{equation}
Observe that \eqref{rel} generalizes \eqref{rel0}.
Notice that each slice $\rel_{\At}(\theta,\velocity)$ is an (open) disc of radius proportional to $(1-\theta^2)|\At\velocity+\im|$. Thus, while for $\At=0$ the relaxation $\rel_0$ only narrows as $|\theta|\uparrow 1$ (i.e.~$z$ tends to $\Kconstrain$), for $0<|\At|<1$ a pinch singularity arises at $\At\velocity+\im=0$ far away from $\Kconstrain$. As we saw in Section \ref{sec:link}, these are the states for which both $\sigma$ and $\varpi$ vanish.\\
The last one $\textrm{(H3)}_\infty$ requires finding bounded subsets $\rel_{\At,M}$ of $\rel_{\At}$ satisfying (H2), which is further laborious than the unbounded case.\\ 
\indent Before embarking on this task (\S\ref{sec:HPproof}-\ref{sec:Relaxation}) we present the statement of our h-principle and we prove Theorems~\ref{thm:Lack}-\ref{thm:DMS} as corollaries.

\begin{defi}\label{defi:weaksol} Let $L_{\Sta}^\infty(\Domain)$ be the (weak*) closed linear subspace of $L_{w*}^{\infty}(\Domain)$ consisting of functions $z=(\theta,\velocity,\aux)$ satisfying the Biot-Savart system
\begin{align}
\int_{\Domain}\velocity\cdot\nabla\phi\dif \x&=0,
\quad\forall\phi\in C_c^1(\bar{\Domain}),
\label{T:2}\tag{\textbf{T}2}\\
\int_{\Domain}(\velocity+\At\aux+\theta\im)\cdot\nabla^\perp\psi\dif\x&=0,
\quad\forall\psi\in C_c^1(\Domain).
\label{T:3}\tag{$\textbf{T}3_\At$}
\end{align}
%with $\int_{\Domain}\velocity\dif\x=0$ for bounded $\Domain$'s.
Notice that (\textbf{T}2) includes the no-flux boundary condition.\\
\indent Let $\theta_0\in L^\infty(\Domain;[-1,1])$ and $T>0$. We say that $\breve{z}=(\breve{\theta},\breve{\velocity},\breve{\aux})\in C([0,T];L_{\Sta}^\infty(\Domain;\bar{\rel}_\At))$ is a \textbf{subsolution} to $(\textrm{IPM}_\At)$ starting from $\theta_0$ if, at each $t\in[0,T]$,
\begin{equation}\tag{\textbf{T}1}
\int_0^t\int_{\Domain}(\breve{\theta}\partial_t\phi+\breve{\aux}\cdot\nabla\phi)\dif\x\dif\tau
=\int_{\Domain}\theta(t)\phi(t)\dif\x
-\int_{\Domain}\theta_0\phi_0\dif\x,
\quad\forall\phi\in C_c^1(\R_+\times\bar{\Domain}).
\label{T:1}
\end{equation}
In particular, a pair $(\theta,\velocity)\in C([0,T];L_{w*}^{\infty}(\Domain;[-1,1]\times\R^2))$ 
is a \textbf{weak solution} to $(\textrm{IPM}_\At)$ if $z=(\theta,\velocity,\theta\velocity)$ is a subsolution to $(\textrm{IPM}_\At)$.\\ 
\indent Let $\breve{z}$ be a subsolution to $(\textrm{IPM}_\At)$ and $\emptyset\neq\Mixzone\subset[0,T]\times\Domain$ open.  
Let us denote $\Mixzone(t)\equiv\{\x\in\Domain\,:\,(t,\x)\in\Mixzone\}$.
We say that $\breve{z}$ is \textbf{strict} w.r.t.~$\Mixzone$ if it is perturbable inside
\begin{equation}\label{strict:1}
\breve{z}\in C(\Omega_{\mathrm{mix}};\rel_\At),
\end{equation}
and exact outside %($\Mixzone(t)\equiv\{\x\in\Domain\,:\,(t,\x)\in\Mixzone\}$)
\begin{equation}\label{strict:2}
\breve{\aux}=\breve{\theta}\breve{\velocity}\quad\textrm{a.e.~in }\Domain\setminus\Mixzone(t),\,\forall t\in[0,T].
\end{equation}
In particular,  we say that $\breve{z}$ is \textbf{admissible} w.r.t.~$\Mixzone$ if it satisfies \eqref{strict:1}, \eqref{strict:2} and
\begin{equation}
|\breve{\theta}|=1\quad\textrm{a.e.~in }\Domain\setminus\Mixzone(t),\,\forall t\in[0,T].
\end{equation}
\end{defi}

\begin{defi}\label{Efuncdef}
In the setting of Theorem \ref{HP} below we need to fix some arbitrary
$\gamma\in[0,1)$, space and time error functions $\Sfunc\in C([0,1];[0,1])$ and $\Tfunc\in C([0,T];[0,1])$ 
with $\Sfunc(0)=\Tfunc(0)=0$ and $\Sfunc(r),\Tfunc(t)>0$ for $r,t>0$. 
With them we define the error function w.r.t.~$\Mixzone$
$$
\Efunc(t,R):=\Sfunc\left(1\wedge\sup_{\x\in R}\dist{\x}{\partial\Mixzone(t)}\right)\Tfunc(t)\frac{1\wedge|R|^\gamma}{|R|},
$$
with %$\expected{R}$ $\equiv$ center of mass of $R$ and 
$|R|$ $\equiv$ area of the bounded rectangle $R\subset\Mixzone(t)$.
% and a $\wedge$ b $\equiv$ $\min(a,b)$.
\end{defi}

\begin{Rem}
The first two terms $\Sfunc$ and $\Tfunc$ defining $\Efunc$ were introduced in \cite{Degraded} to show that the error in Theorem~\ref{HP}\ref{degraded} below depends on the distance to the (space-time) boundary of the mixing zone, and the parameter $\gamma$ to refine this estimate for small rectangles.
However, for simplicity one may consider $\Efunc(t,R)=\Tfunc(t)/|R|$ since it contains relevant information and it is easier to understand in a first reading (cf.~\cite[Rem.~1.1]{Degraded}).
\end{Rem}

\begin{thm}[H-principle for $(\textrm{IPM}_\At)$]\label{HP} 
Let $|\At|<1$, $T>0$, $\emptyset\neq\Omega_{\mathrm{mix}}\subset(0,T]\times\Domain$ open and $\Efunc$ as in Def.~\ref{Efuncdef}. 	
Suppose there is a strict subsolution 
$\breve{z} 
%\in C([0,T];L_{\Sta}^\infty(\Domain;\bar{\rel}_\At))
$ to $(\mathrm{IPM}_\At)$ w.r.t.~$\Omega_{\mathrm{mix}}$.
Then, there exist infinitely many weak solutions $(\theta,\velocity)
%\in C([0,T];L_{w*}^{\infty}(\Domain))
$ to $(\mathrm{IPM}_\At)$ satisfying that, at each $t\in[0,T]$:
\begin{enumerate}[(a)]
\item They agree with $\breve{z}$ outside $\Omega_{\mathrm{mix}}$
$$(\theta,\velocity)(t)=(\breve{\theta},\breve{\velocity})(t)
\quad \textrm{in }\Domain\setminus\Omega_{\mathrm{mix}}(t).$$
\item\label{Mixingprop} For every (bounded) open $\emptyset\neq\Omega\subset\Mixzone(t)$,
$$\int_{\Omega}(1-\theta(t,\x)^2)\dif \x=0<\int_{\Omega}(1-\theta(t,\x))\dif \x\int_{\Omega}(1+\theta(t,\x))\dif\x.$$
\item\label{degraded} For every bounded rectangle $\emptyset\neq R\subset\Omega_{\mathrm{mix}}(t)$,
$$
\left|\dashint_{R}[\mathbf{F}(z)-\mathbf{F}(\breve{z})](t,\x)\dif\x\right|\leq\Efunc(t,R),
$$
for $\mathbf{F}=\mathrm{id}$ or $\mathbf{P}(\breve{z})
:=\breve{\velocity}\cdot(\breve{\velocity}+\At\breve{\aux}+\breve{\theta}\im)$, where $\dashint_{R}\equiv\tfrac{1}{|R|}\int_R$ and $z=(\theta,\velocity,\theta\velocity)$.
\end{enumerate}
In addition, if $\breve{z}$ is admissible w.r.t.~$\Mixzone$, then $\theta\in C([0,T];L_{\mathrm{loc}}^p(\Domain;\{-1,1\}))$ for all $1<p<\infty$. 
\end{thm}

The choice of $\breve{z}$ in Theorem~\ref{thm:Lack} is related to \cite{LIPM,RIPM}, but in order to guarantee the weak*-continuity of the non-linearity $\theta\velocity$ we have chosen a time dependent $\breve{\aux}$.
\begin{proof}[Proof of Theorem~\ref{thm:Lack}] By Theorem \ref{HP}, we consider $\Mixzone=(0,T)\times\Domain$ and $\breve{z}=(0,0,\breve{\aux})$ with $\breve{\aux}\in C([0,T];\R^2)$ satisfying  $\breve{\aux}(0)=\breve{\aux}(T)=0$
and $|2\breve{\aux}(t)+\im|<1$ for all $t\in(0,T)$.
\end{proof}

Similarly, Theorem~\ref{thm:DMS} can be proved as a corollary of the above h-principle. Before writing the proof, let us reformulate it with the new terminology.

\begin{thm}\label{thm:DMS2}  Let $|\At|<1$, $\Domain=\R^2$ and $0<\alpha<1$. Then $\breve{z}_{\At,\alpha}$ with
\begin{equation}\label{Osubsol}
\breve{\theta}_{\At,\alpha}(t)
=\Theta_{\At}(\alpha t),\quad t\in\R_+,
\end{equation}
$\breve{\velocity}_{\At,\alpha}=0$ and $\breve{\aux}_{\At,\alpha}$ given by \eqref{Msubsol},
is an admissible subsolution to $(\mathrm{IPM}_\At)$ w.r.t.
\begin{equation}\label{mixzone}
\Mixzone=\{(t,\x)\in\R_+\times\Domain\,:\,
-\alpha c_{\At}^-t< x_2<\alpha c_{\At}^+ t\}.
\end{equation}
For $\Domain=(-1,1)^2$ the same holds except that \eqref{mixzone} is only valid until $\Mixzone(t)$ meets either the lower or upper boundary of $(-1,1)^2$. After this, $\Mixzone(t)$ starts to reduce until it ends up collapsing and the stable planar phase is reached (cf.~\S\ref{sec:transition}).
\end{thm}
\begin{defi}\label{def:DMS}
We say that the weak solutions 
$(\theta,\velocity)
%\in C(\R_+;L_{w*}^{\infty}(\Domain))
$ coming from the h-principle  applied to this $\breve{z}_{\At,\alpha}$ are $\Theta_{\At}$-\textbf{mixing solutions} %of degree $\Efunc$ and mixing speed $\alpha$ 
to $(\mathrm{IPM}_\At)$ starting from the unstable planar phase \eqref{flat}.
%In the setting of Def.~\ref{def:DMS}, 
For $\Domain=\R^2$ let us denote
\begin{equation}\label{Omegapm}
\Omega_{\pm}=\{(t,\x)\in\R_+\times\Domain\,:\,
\pm x_2>\alpha c_{\At}^{\pm} t\}.
\end{equation}
As in Thm.~\ref{thm:DMS2}, for $\Domain=(-1,1)^2$ \eqref{Omegapm} changes once $\Mixzone(t)$ hits either  $x_2=-1$ or $1$ (cf.~\S\ref{sec:transition}).\\
Thus, at each $t\in\R_+$, these $\Theta_{\At}$-mixing solutions satisfy:
\begin{enumerate}[(a)]
	\item\label{DMS:1} Non-mixing outside $\Omega_{\mathrm{mix}}$:
	$$(\theta,\velocity)(t)=(\pm 1,0)
	\quad\textrm{in }\Omega_{\pm}(t).$$
	\item\label{DMS:2} Mixing inside $\Omega_{\mathrm{mix}}$: For every (bounded) open $\emptyset\neq\Omega\subset\Mixzone(t)$, 
	$$	\int_{\Omega}(1-\theta(t,\x)^2)\dif\x=0<\int_{\Omega}(1-\theta(t,\x))\dif \x\int_{\Omega}(1+\theta(t,\x))\dif\x.$$
	\item\label{DMS:3} $\Theta_{\At}$-macroscopic behaviour: For every bounded rectangle $\emptyset\neq R=S\times tL\subset\Omega_{\mathrm{mix}}(t)$,
	$$
	\left|\dashint_{R}\theta(t,\x)\dif \x-\expected{L}_{\At,\alpha}\right|\leq\Efunc(t,R)
	\quad\textrm{where}\quad
	\expected{L}_{\At,\alpha}:=\dashint_{L}\Theta_{\At}(\alpha,x_2)\dif x_2.
	$$
	\item\label{DMS:4} For $f(\theta,\velocity)=\velocity$, $\theta\velocity$ and $\mathbf{P}(\theta,\velocity,\theta\velocity)$, and every bounded rectangle $\emptyset\neq R\subset\Omega_{\mathrm{mix}}(t)$,
	$$
	\left|\dashint_{R}f(\theta,\velocity)(t,\x)\dif\x\right|\leq\Efunc(t,R).
	$$
\end{enumerate}
\end{defi}
\begin{Rem} The properties \ref{DMS:1}\ref{DMS:2} justify the adjective ``mixing'' and \ref{DMS:3} the tag ``$\Theta_{\At}$'' (cf.~Rem. \ref{uncernainity} and Prop.~\ref{Oprop} for a explicit computation of $\expected{L}_{\At,\alpha}$). The property \ref{DMS:4} shows that $\breve{\velocity}_{\At,\alpha}=0$ can be interpreted as the macroscopic velocity too, and also that the ``power balance'' $\mathbf{P}$ (cf.~\cite[(14)]{Degraded}), which is a quadratic quantity, is almost preserved. 
\end{Rem}

\begin{Rem}\label{diagram} In \cite[Rem.~5]{RIPM}, the interpretation that $\breve{\theta}_{0,\alpha}$ represents the coarse-grained phase follows from the fact that there is a sequence of exact solutions %$(\theta_k)$ satisfying  
$\theta_k\overset{*}{\rightharpoonup}\breve{\theta}_{0,\alpha}$. Here, the property \ref{DMS:3} closes the diagram \eqref{HPscheme} in the sense that it provides an explicit relaxation for each exact solution separately. Schematically, if we denote $X_{\At,\alpha}$ by the space of these $\Theta_{\At}$-mixing solutions with mixing speed $\alpha$, then we have
\begin{center}
	\begin{tikzpicture}
	\node (1){$X_{\At,\alpha}$};     
	\node (2) [node distance=2.5cm, right of=1] {$\breve{\theta}_{\At,\alpha}$}; 
	\path[->] (1) edge [bend left, auto=left] node {average} (2);
	\path[->] (2) edge [bend left, auto=left] node {h-principle} (1);
	\end{tikzpicture}
\end{center}
where the upper arrow means that $\breve{\theta}_{\At,\alpha}$ can be recovered from each $\theta\in X_{\At,\alpha}$ by averaging it over horizontal lines as follows
$$\breve{\theta}_{\At,\alpha}(t,\x)=
\lim_{M\rightarrow\infty}\dashint_{R_M(\x)}\theta(t,\x')\dif \x',\quad
(t,\x)\in\R^2\times\R_+,$$
with $R_M(\x)=\x+(-M,M)\times(-M^{-\delta},M^{-\delta})$ for some arbitrary $\delta\in(0,1)$. 
\end{Rem}
\begin{proof}[Proof of Theorem~\ref{thm:DMS2}] Consider
$\breve{\theta}=\breve{\theta}(t,x_2)$, $\breve{\velocity}=0$ and $\breve{\aux}$ to be determined. The condition \eqref{strict:1} reads as $\breve{z}$ maps continuously $\Mixzone$ into
$$|2(1-\breve{\theta}\At)\breve{\aux}+(1-\breve{\theta}^2)\im|
<(1-\theta^2).$$
This suggests to take, for some $0<\alpha<1$,
\begin{equation}\label{Msubsol}
\breve{\aux}= -\alpha\frac{1-\breve{\theta}^2}{1-\breve{\theta}\At}\im.
\end{equation}
On the one hand, $(\textrm{\textbf{T}2-3}_\At)$ is automatically satisfied. On the other hand, \eqref{T:1} reads as
\begin{equation}\label{conservationlaw:alpha}
\partial_t\breve{\theta}=\alpha\partial_{x_2}\left(\frac{1-\breve{\theta}^2}{1-\breve{\theta}\At}\right).
\end{equation}
The (unique) entropy solution of the above scalar conservation law is \eqref{Osubsol}. Finally, it is clear that $\breve{z}$ is admissible w.r.t.~$\Mixzone$.  
\end{proof}

We conclude this section by extending Prop.~4.3 in \cite{RIPM} to the general case $|\At|<1$. 
Roughly speaking this reads as, among subsolutions $\breve{z}$ to $(\textrm{IPM}_\At)$ starting from \eqref{flat} with planar symmetry, the borderline case $\alpha=1$ in Thm.~\ref{thm:DMS2} maximizes the mixing zone. As suggested in \cite{RIPM}, this may serve as a selection criterion. We remark in passing that, inspired by \cite{Sze11}, the intermediate case $\alpha=\tfrac{1}{2}$, which maximizes the energy dissipation rate for the Kelvin-Helmholtz instability, may contain relevant physical information and then should be explored in future works. \\
\indent Let us assume that $\partial_{x_1}\breve{z}=0$ and that both fluids are at rest ($\breve{\velocity}=0$) outside $\Mixzone$. Then, $(\textrm{\textbf{T}2-3}_\At)$ implies that
\begin{equation}\label{max:1}
\breve{\velocity}=-\At\breve{\aux}_1.
\end{equation}
Notice that $\At=0$ yields $\breve{\velocity}=0$. Indeed, in \cite{RIPM} $\breve{\velocity}=0$ follows from the slighter assumption $\partial_{x_1}\breve{\theta}=0$. Although Proposition~\ref{maxmix} below holds in the class $\breve{\velocity}=0$ too, we find more natural the condition \eqref{max:1} here.\\ 
\indent As in \cite{RIPM}, on the confined domain $(-1,1)^2$ the no-flux boundary condition implies $\breve{\velocity}=0$. Therefore, Prop.~4.3 in \cite{RIPM} can be extended analogously for $\Domain=(-1,1)^2$.
However, if we remove the vertical walls, say $\Domain=\T\times(-1,1)$, then \eqref{max:1} requires some extra computations. Let us see it. Notice that $\At\breve{\velocity}+\im\neq0$ because $\breve{\velocity}_2=0$.
Then, since $\breve{z}$ is $\bar{\rel}_\At$-valued, the following inequality holds (a.e.)
\begin{equation}\label{max:2}
\left|\frac{2(1-\breve{\theta}\At)(\breve{\aux}-\breve{\theta}\breve{\velocity})}{\At\breve{\velocity}+\im}+(1-\breve{\theta}^2)\right|\leq (1-\breve{\theta}^2).
\end{equation}
By taking the real part of \eqref{max:2} and applying \eqref{max:1}, we get
$$
-\frac{1-\breve{\theta}^2}{1-\breve{\theta}\At}
\leq
\Re\left(\frac{\breve{\aux}-\breve{\theta}\breve{\velocity}}{\At\breve{\velocity}+\im}\right)
=-\frac{\Re((\breve{\aux}+\breve{\theta}\At\breve{\aux}_1)(\At^2\breve{\aux}_1+\im))}{1+(\At^2\breve{\aux}_1)^2}
=\frac{\breve{\aux}_2-(1+\breve{\theta}\At)(\At\breve{\aux}_1)^2}{1+(\At^2\breve{\aux}_1)^2},
$$
and so
\begin{equation}\label{max:4}
\begin{split}
\breve{\aux}_2&\geq -\frac{1-\breve{\theta}^2}{1-\breve{\theta}\At}(1+(\At^2\breve{\aux}_1)^2)+(1+\breve{\theta}\At)(\At\breve{\aux}_1)^2\\
&=-\frac{1-\breve{\theta}^2}{1-\breve{\theta}\At}+(\At\breve{\aux}_1)^2\frac{1-\At^2}{1-\breve{\theta}\At}\geq -\frac{1-\breve{\theta}^2}{1-\breve{\theta}\At}.
\end{split}
\end{equation}
The rest follows similarly to \cite{RIPM}.
Let us denote $\varOmega_{\pm}=\{(t,\x)\in\R_+\times\Domain\,:\,\pm x_2>c_\At^{\pm} t\}$.
By approximation, $\phi^{\pm}(t,\x)=(\pm x_2-c_\At^{\pm}t)\vee 0$  is a valid test function. 
 Then, since
$$c_{\At}^{\pm}|\varOmega_{\pm}|=1=\pm\int_{\Domain}\theta_0\phi_0^{\pm}\dif\x,$$
by evaluating (\textbf{T}1) with $\phi^{\pm}$ we obtain
$$\int_{\varOmega_{\pm}}(c_\At^{\pm}(1\mp\breve{\theta})+\breve{\aux}_2)\dif\x\dif t=0.$$
Finally, since \eqref{max:4} implies
$$c_\At^{\pm}(1\mp\breve{\theta})+\breve{\aux}_2
\geq(1\mp\breve{\theta})\left(\frac{2}{1\mp\At}-\frac{1\pm\breve{\theta}}{1-\breve{\theta}\At}\right)
=\frac{(1\mp\breve{\theta})^2(1\pm\At)}{(1-\breve{\theta}\At)(1\mp\At)}
\geq 0,$$
necessarily $\breve{\theta}=\pm 1$ in $\varOmega_{\pm}$. In summary, at least for bounded and rectangular $\Domain$'s (cf.~\cite{CGO14}), either with or without vertical boundaries, the following holds.

\begin{prop}\label{maxmix} Let $\breve{z}$ be a subsolution to $(\mathrm{IPM}_\At)$ starting from \eqref{flat} w.r.t.~some $\Mixzone$ and satisfying \eqref{max:1}. Then $\breve{\theta}=\pm 1$ in $\varOmega_{\pm}$, i.e.~$\Mixzone\subset\{(t,\x)\in\R_+\times\Domain\,:\,
-c_{\At}^-t< x_2<c_{\At}^+ t\}$.
\end{prop}

\section{Proof of the h-principle}\label{sec:HPproof}

In this section we prove Theorem~\ref{HP}. To this end, we need to check the following three hypothesis (cf.~\cite{Degraded,Onadmissibility,RIPM}). We do so for $p=\infty$ and also for $p=2$ on $\Domain=\T^2$. Although $L^\infty(\T^2)\subset L^2(\T^2)$, the direct proof (Prop.~\ref{L2bound}) for $p=2$ shows that $\bar{\rel}_\At$ is somehow sharp.\\

\textbf{(H1) Localized plane waves}. Let $0\neq h\in C^1(\T;[-1,1])$ with $\int h=0$. There is a cone $\Lambda\subset\R^5$ so that, for all $\bar{z}\in\Lambda$ and $\psi\in C_c^\infty(\R^3)$ there is $\xi\in\R\times\Sp^1$ for which there are smooth solutions to $(\textbf{T}_\At)$ of the form
$$z_k(t,\x)=\bar{z}h(k\xi\cdot(t,\x))\psi(t,\x)+\mathcal{O}(k^{-1}),$$ 
with $k\in\N$ and $\mathcal{O}$ depending on $|\bar{z}|,|\xi|$ and
$\{|D^\beta\psi(t,\x)|\,:\,1\leq|\beta|\leq 2\}$.\\

\textbf{(H2) Perturbation property}. There is an open set $U\subset\set$ and a function $\Phi\in C(]0,1];]0,1])$ such that, for all $z\in U$ there is $\bar{z}\in\Lambda$ with $\bar{\theta}=1$
for which
$$z+\lambda\bar{z}\in U,\quad|\lambda|\leq \Phi(1-\theta^2).$$

$\textbf{(H3)}_p$ \textbf{Weak*-compactness}. The space $L_{\Sta}^p(\Domain;\bar{U})$ is $L^p$-bounded.\\

%From now on we identify $\R^2\simeq\C$ as usual $\mathbf{v}=(v_1,v_2)=v_1+v_2\im$. Thus, without any ambiguity we shall denote $\mathbf{v}^*=v_1-v_2\im$, $\mathbf{v}^\perp=\im\mathbf{v}
%=-v_2+v_1\im$ and $\mathbf{v}\cdot\mathbf{w}=\Re(\mathbf{v}\mathbf{w}^*)=
%v_1w_1+v_2w_2$.\\
\indent Let us start checking (H1). 
Since
$$\det\mathbf{T}_\At(\bar{z})
=-\bar{\theta}\bar{\velocity}\cdot(\bar{\velocity}+\At\bar{\aux}+\bar{\theta}\im)
=\tfrac{1}{4}\bar{\theta}(|\At\bar{\aux}+\bar{\theta}\im|^2-|2\bar{\velocity}+\At\bar{\aux}+\bar{\theta}\im|^2),$$
from the definition of the wave cone \eqref{WC} it follows that
\begin{equation}\label{Lambda}
\Lambda_\At=\varLambda_0\cup\varLambda_1,
\end{equation}
with $\varLambda_j\equiv\varLambda_{\At,j}$ given by
\begin{align*}
\varLambda_0&:=\{\bar{z}\in\R^5\,:\,\bar{\theta}= 0,\,\bar{\velocity}=-\At\bar{\aux}\},\\
\varLambda_1&:=\{\bar{z}\in\R^5\,:\,\bar{\theta}\neq 0,\,\bar{\velocity}=\bar{\omega}(\At\bar{\aux}+\bar{\theta}\im)\textrm{ for some }\bar{\omega}\in\Ss\},
\end{align*}
where
\begin{equation}\label{S}
\Ss:=\{\bar{\omega}\in\R^2\,:\,|2\bar{\omega}+1|=1\},
\end{equation}
that is, $\Ss$ is the sphere centered at $-\tfrac{1}{2}$ with radius $\tfrac{1}{2}$. We shall also consider the interior of its convex hull $\Sc=(\Ss^{\mathrm{co}})^\circ=\{\omega\in\R^2\,:\,|2\omega+1|< 1\}$. Both can be expressed in terms of the unit sphere $\Sp$ and the unit disc $\D$ as 
$\Sp=\shift\Ss$ and $\D=\shift\Sc$
where $T:\bar{\Sc}\rightarrow\bar{\D}$ is the translation
\begin{equation}\label{shift}
T\omega:=2\omega+1.
\end{equation}

\begin{lemma}\label{LPW} $\mathrm{(H1)}$ holds for $\Lambda=\Lambda_\At$.
\end{lemma}
\begin{proof}
\textit{Step 1. Construction of a potential}: Let us suppose that $z=(\theta,\velocity,\aux)$ is a smooth localized solution to $(\textbf{T}_\At)$. Then, by \eqref{T:2}, $\velocity=\nabla^\perp f$ for some smooth $f$. If we write $\aux$ in its Hodge's decomposition, $\aux=\nabla^\perp\varphi+\nabla g$ for some smooth $\varphi,g$, then \eqref{T:1} and $(\textrm{\textbf{T}3}_\At)$ read as
\begin{align*}
\partial_t\theta+\Delta g&=0,\\
\Delta(f+\At \varphi)+\partial_{x_1}\theta&=0.
\end{align*}
Notice that $\theta=\Delta\phi$ for some smooth $\phi$. Hence, $g=-\partial_t\phi$ and $f=-(\partial_{x_1}\phi+A\varphi)$. In summary,
$$
\theta=\Delta\phi,\quad\quad
\velocity=-\nabla^\perp(\partial_{x_1}\phi+\At\varphi),\quad\quad
\aux=\nabla^\perp\varphi-\partial_t\nabla\phi.
$$
This suggests to consider the following potential
$$P(\phi,\varphi)
:=(\Delta\phi,-\nabla^\perp(\partial_{x_1}\phi+\At\varphi),\nabla^\perp\varphi-\partial_t\nabla\phi).$$
Since $\velocity+\At\aux+\theta\im=\nabla(\partial_{x_2}-\At\partial_t)\phi$, 
it satisfies $\Div\mathbf{T}_{\At}(P(\phi,\varphi))=0$ for all $\phi,\varphi\in C^3(\R^3)$. \\
\indent\textit{Step 2. Construction of $z_k$}: Let us take $H\in C^3(\T)$ such that $H''=h$.\\
Given $\bar{z}=(\bar{\theta},\bar{\velocity},\bar{\aux})\in\Lambda$ and $k\in\N$,
we consider
$$\phi_k(t,\x)=\tfrac{a}{k^2}H(k\xi\cdot(t,\x))
,\quad\varphi_k(t,\x)=\tfrac{b}{k}H'(k\xi\cdot(t,\x)),$$
with $\xi=(\xi_0,\zeta)\in\R\times\Sp^1$ and $a,b\in\R$ to be determined. This choice yields
$$P(\phi_k,\varphi_k)(t,\x)
=
(a,-\im(a\zeta_1+b\At)\zeta,(b\im-a\xi_0)\zeta)
h(k\xi\cdot\x).$$
Then, to prove (H1) we need to find $\xi$, $a$, $b$ satisfying
\begin{equation}\label{P:1}
(a,-\im(a\zeta_1+b\At)\zeta,(b\im-a\xi_0)\zeta)=(\bar{\theta},\bar{\velocity},\bar{\aux}).
\end{equation}
The first column in \eqref{P:1} reads as $a=\bar{\theta}$.
Firstly assume that $\bar{z}=\varLambda_0$, i.e.~$a=0$ and $\bar{\velocity}=-\At\bar{\aux}$. Hence, the second and third column in \eqref{P:1} are equivalent to $\bar{\aux}=b\zeta^\perp$.
Thus, we take $b=|\bar{\aux}|$ and $\zeta\in\Sp^1$ such that $\bar{\aux}=b\zeta^\perp$.
Secondly assume that $\bar{z}\in\varLambda_1$, i.e.~$a\neq 0$ and there is $\bar{\omega}\in\Ss$ so that $\bar{\velocity}=\bar{\omega}(\At\bar{\aux}+a\im)$. 
Hence, for the third column in \eqref{P:1}, $\bar{\aux}=(b\im-a\xi_0)\zeta$, necessarily $\xi_0=-a^{-1}\bar{\aux}\cdot\zeta$ and
$b=\bar{\aux}\cdot\zeta^\perp$. 
Now, the second column in \eqref{P:1} reads as $\bar{\velocity}=-\im(a\zeta_1+b\At)\zeta
=-\zeta^\perp(\At\bar{\aux}+a\im)\cdot\zeta^\perp$.
Since $\bar{\velocity}=\bar{\omega}(\At\bar{\aux}+a\im)$, $\zeta$ is given by the equation
$$
\bar{\omega}(\At\bar{\aux}+a\im)
=-\zeta^\perp(\At\bar{\aux}+a\im)\cdot\zeta^\perp.
$$
If $\bar{\omega}(\At\bar{\aux}+a\im)=0$, we take $\zeta\parallel(\At\bar{\aux}+a\im)$. Otherwise, we take ($|\bar{\omega}|^2=-\bar{\omega}_1$)
$$\zeta^\perp=
\pm\frac{\bar{\omega}}{|\bar{\omega}|}\frac{\At\bar{\aux}+a\im}{|\At\bar{\aux}+a\im|}.$$
\indent Finally, we consider $z_k=P(\phi_k\psi,\varphi_k\psi)$ because
$$z_k-\bar{z}h\psi=P(\phi_k\psi,\varphi_k\psi)-P(\phi_k,\varphi_k)\psi
=\mathcal{O}(k^{-1}),$$
as we wanted.
\end{proof}

\begin{lemma}\label{H2rel} $\mathrm{(H2)}$ holds for $U=\rel_\At$.
\end{lemma}

We will prove this lemma in Section \ref{sec:unbounded}. Now, we check $\textrm{(H3)}_2$ on $\Domain=\T^2$. To this end, it is convenient to normalize $L_{\Sta}^2(\T^2;\bar{\rel}_\At)$ by imposing $\int\velocity=0$ therein.

\begin{prop}\label{L2bound} The space $L_{\Sta}^2(\T^2;\bar{\rel}_\At)$ is $L^2$-bounded.
\end{prop}
\begin{proof}
Let $z\in L_{\Sta}^2(\T^2;\bar{\rel}_\At)$.
On the one hand, since $z$ is $\bar{\rel}_\At$-valued, we will see in Lemma~\ref{lemma:U}\ref{Ldirection'} that $\aux$ can be expressed (a.e.) as
\begin{equation}\label{L2:1}
\aux=\theta\velocity+\frac{(1-\theta^2)(\At\velocity+\im)\omega}{1+\omega\theta\At}
=\frac{(\theta+\omega\At)\velocity+(1-\theta^2)\im\omega }{1+\omega\theta\At},
\end{equation}
for some $\bar{\Sc}$-valued $\omega$.
Hence, by applying
\begin{equation}\label{L2:2}
\left|\frac{\theta+\omega\At}{1+\omega\theta\At}\right|^2
%=\frac{\theta^2+\At^2|\omega|^2+2\theta\At\omega_1}{1+\theta^2\At^2|\omega|^2+2\theta\At\omega_1}
=1-(1-\theta^2)\frac{1-|\omega|^2\At^2}{|1+\omega\theta\At|^2}
\leq 1,
\end{equation}
the triangle inequality yields
\begin{equation}\label{L2:4}
|\aux|
\leq |\velocity|+\frac{1-\theta^2}{1-|\theta\At|}
\leq |\velocity|+(1+|\theta|).
\end{equation}
On the other hand, since $(\textrm{\textbf{T}2-3}_\At)$ is written in the Fourier side as
\begin{equation*}
\hat{\velocity}(k)\cdot k=0,\quad
(\hat{\velocity}+\At\hat{\aux}+\hat{\theta}\im)(k)\cdot k^\perp=0,
\quad k\in\Z^2,
\end{equation*}
and we have normalized $\hat{\velocity}(0)=0$, the velocity $\velocity$ is given by
$$\hat{\velocity}(k)
=-\frac{k^\perp}{|k|^2}(\At\hat{\aux}+\hat{\theta}\im)(k)\cdot k^\perp,\quad k\in\Z^2.$$
Therefore, Plancherel's identity and the triangle inequality yield
\begin{equation}\label{L2:3}
\|\velocity\|_2\leq\|\At\aux+\theta\im\|_2
\leq |\At|\|\aux\|_2+\|\theta\|_2.
\end{equation}
This concludes the proof since $|\theta|\leq 1$ and because \eqref{L2:4}\eqref{L2:3} imply
$$\|\velocity\|_2
\leq\frac{\|\theta\|_2+|\At|\|1+|\theta|\|_2}{1-|\At|},
\quad\quad
\|\aux\|_2\leq
\frac{\|\theta\|_2+\|1+|\theta|\|_2}{1-|\At|}.$$
\end{proof}

Thus, $\textrm{(H1)-(H3)}_2$ hold on $\Domain=\T^2$. In order to prove it for $p=\infty$ we need to find bounded $U$'s satisfying (H2). To this end, we will prove the following lemma in Section \ref{sec:bounded}.

\begin{lemma}\label{H2bounded} For any $R>0$ there is a bounded open subset $U$ of $\rel_\At$ satisfying $\mathrm{(H2)}$ and 
$$\{z\in\rel_\At\,:\, |\velocity|< R\}\subset U.$$
Obviously, $\mathrm{(H3)}_\infty$ holds for $U$.
\end{lemma}

\begin{Rem}\label{Rem:stability}
At this point we have all the ingredients to apply the h-principle in \cite{Degraded}, except we do not know if $L_{\Sta}^p(\Domain;\bar{\rel}_\At)$ is (weak*) closed. 
Although we have not been able to show it, we have noticed that the proof of this h-principle can be adapted to $\bar{\rel}_\At$. In brief, the original proof uses this property to show that a certain set ``$J^{-1}(0)$'' consists of functions $z$ solving $(\mathbf{T}_\At,\Kconstrain)$. Here, we overcome this obstacle by checking that the residual subset ``$X_J$''of $J^{-1}(0)$ satisfies this requirement.
\end{Rem}

\begin{proof}[Proof of Theorem~\ref{HP}]
Let $\breve{z}\in C([0,T];L_{\Sta}^p(\Domain;\bar{\rel}_\At))$ be a strict subsolution to $(\mathrm{IPM}_\At)$ w.r.t.~$\Mixzone$. For $p=2$ we take $U=\rel_\At$ and for $p=\infty$ we take $U$ from Lemma~\ref{H2bounded} in such a way that $|\breve{\velocity}|<R$.
\indent Now, let us recall how ``$X_0$'' is defined in \cite{Degraded}. A subsolution $z\in C([0,T];L_{\Sta}^p(\Domain;\bar{\rel}_\At))$ belongs to $X_0$ if it agrees with $\breve{z}$ outside $\Mixzone$
$$z=\breve{z}\quad\textrm{a.e.~in }\Domain\setminus\Mixzone(t),\,\forall t\in[0,T],$$
and it is perturbable inside
$$z\in C(\Omega_{\mathrm{mix}};U).$$
In addition, we ask $z$ to satisfy the following property. There is $C(z)\in(0,1)$ so that, at each $t\in[0,T]$, for $\mathbf{F} = \mathrm{id}$ and $\mathbf{P}$,
$$\left|\dashint_{R}[\mathbf{F}(z)-\mathbf{F}(\breve{z})](t,\x)\dif\x\right|\leq C(z)\Efunc(t,R),$$
for every bounded rectangle $\emptyset\neq R\subset\Mixzone(t)$. By $\textrm{(H3)}_p$, the closure $X$ of $X_0$ in $C([0,T];L_{\Sta}^p(\Domain))$ is a completely metrizable space.\\ 
\indent Given $\Omega\Subset\Domain$ open and $I=[t_1,t_2]\subset[0,T]$, the relaxation-error functional is defined in \cite{Degraded} as
\function{J}{X}{\R_+}{z}{\sup_{t\in I}\int_{\Omega}(1-\theta(t,\x)^2)\dif\x ,}
which is well defined because, by convexity, $|\theta|\leq 1$ for states in $X$. Indeed, $J$ is upper-semicontinuous, and so the set $X_J$ of continuity points of $J$ is residual (countable intersection of open dense sets). Then, following \cite{Degraded}, the hypothesis $\textrm{(H1)-(H3)}_p$ imply that $X_J\subset J^{-1}(0)$. 
In contrast to \cite{Degraded}, here we can not use that $L_{\Sta}^p(\Domain;\bar{\rel}_\At)$ is (weak*) closed to ensure that the functions in $J^{-1}(0)$ are $\Kconstrain$-valued in $I\times\Omega$. However, we shall prove that $X_J$ satisfies this requirement. \\
\indent Given $z\in X_J$ let $(z_k)\subset X_0$ converging to $z$. Fix $t\in I$. We claim that $\theta_k(t)\rightarrow\theta(t)$ in $L^q(\Omega)$ for every $1<q<\infty$. Indeed, since $J(z)=0$ and
$$|\|\theta(t)\|_q^q-\|\theta_k(t)\|_q^q|
=\int_\Omega(1-|\theta_k(t,\x)|^q)\dif\x
\leq C_q\int_\Omega(1-\theta_k(t,\x)^2)\dif\x
\leq C_q J(z_k)\rightarrow C_qJ(z)=0,$$
the claim follows by convexity. 
Now take $1<q<p$ and denote $f=(\aux-\theta\velocity)$, $f_k=(\aux_k-\theta\velocity_k)$ and $\tilde{f}_k=(\aux_k-\theta_k\velocity_k)$.
On the one hand, by convexity and applying $f_k(t)\overset{*}{\rightharpoonup}f(t)$, we get
\begin{equation}\label{HP:1}
\int_\Omega|\aux-\theta\velocity|^q(t,\x)\dif\x
=\|f(t)\|_q^q
\leq\liminf_{k}\|f_k(t)\|_q^q
=\liminf_{k}\int_\Omega|\aux_k-\theta\velocity_k|^q(t,\x)\dif\x.
\end{equation}
On the other hand, by applying the inverse triangle inequality, we obtain
\begin{equation}\label{HP:2}
|\|f_k(t)\|_q-\|\tilde{f}_k(t)\|_q|^q\leq\|f_k(t)-\tilde{f}_k(t)\|_q^q=\int_\Omega|\theta-\theta_k|^q|\velocity_k|^q(t,\x)\dif\x\rightarrow 0,
\end{equation}
where the last convergence follows from H\"{o}lder's inequality  and $\textrm{(H3)}_p$. Finally, by applying \eqref{HP:1}\eqref{HP:2} and that $z_k$ is $\bar{\rel}_\At$-valued, we deduce
\begin{align*}
\int_\Omega|\aux-\theta\velocity|^q(t,\x)\dif\x
&\leq\liminf_k\int_\Omega|\aux_k-\theta\velocity_k|^q(t,\x)\dif\x\\
&=\liminf_k\int_\Omega|\aux_k-\theta_k\velocity_k|^q(t,\x)\dif\x\\
&\leq\liminf_k\int_\Omega(1-(\theta_k)^2)^q\frac{|\At\velocity_k+\im|^q}{(1-\theta_k\At)^q}(t,\x)\dif\x=0,
\end{align*}
and so $\aux=\theta\velocity$. Therefore, $z(t)$ is $\Kconstrain$-valued on $\Omega$. The rest follows as in \cite{Degraded}.
\end{proof}

\section{The relaxation}\label{sec:Relaxation}

%In this section we prove Lemmas~\ref{H2rel}-\ref{H2bounded}. 
First of all let us recall several notions in Lamination Theory. Given a set $K$ and a cone $\Lambda$ in $\R^N$, the $\Lambda$-lamination of order $1$ of $K$ is
\begin{equation}\label{def:K1L}
K^{1,\Lambda}:=\{\tfrac{1+s}{2}z_1+\tfrac{1-s}{2}z_2\,:\,s\in[-1,1],\,z_1,z_2\in K\textrm{ s.t. }z_1-z_2\in\Lambda\},
\end{equation}
and, inductively, the $\Lambda$-lamination of order $n\geq 2$ of $K$ is
$$K^{n,\Lambda}:=(K^{n-1,\Lambda})^{1,\Lambda}.$$
This generates an ascending chain of sets $K\subset K^{1,\Lambda}\subset K^{2,\Lambda}\subset\cdots$ whose limit $K^{lc,\Lambda}:=\bigcup K^{n,\Lambda}$ is the $\Lambda$-\textbf{lamination hull} of $K$. 
This is contained in the $\Lambda$-\textbf{convex hull} of $K$ which is defined as follows: A state $z\in\R^N$ does not belong to $K^\Lambda$ if there is a $\Lambda$-convex function $f$ (meaning that $\lambda\mapsto f(z_0+\lambda\bar{z})$ is convex for all $z_0\in\R^N$ and $\bar{z}\in\Lambda$) so that $f\leq 0$ on $K$ and $f(z)>0$.\\
\indent From now on we consider $\Kconstrain$ and $\Lambda_{\At}$ given in \eqref{K} and \eqref{Lambda} respectively.
In order to alleviate the notation we shall omit the tag ``$\At$'' wherever we do not need to distinguish between the cases $\At=0$ and $\At\neq 0$. Thus, we shall abbreviate $\mathbf{T}\equiv\mathbf{T}_{\At}$, $\Lambda\equiv\Lambda_\At$ and $\rel\equiv\rel_\At$.\\
\indent This section is split in three parts. Firstly we compute $\Kconstrain^{1,\Lambda}$ since it contains the key to understand the relaxation. Secondly we prove Lemmas \ref{H2rel} (\S\ref{sec:unbounded}) and \ref{H2bounded} (\S\ref{sec:bounded}). Finally we check that $\Kconstrain^{lc,\Lambda}=\bar{\rel}$ and $(\Kconstrain_M)^{lc,\Lambda}=\bar{\rel}_M$ (\S\ref{sec:Lambdahull}).

\begin{lemma}\label{lemma:K1} Let $z=(\theta,\velocity,\aux)\in\set$. The following are equivalent:
\begin{enumerate}[(a)]
	\item\label{K1:L} $z\in\Kconstrain^{1,\Lambda}$.
	\item\label{K1:L1} $z\in\Kconstrain^{1,\varLambda_{1}}$.
	\item\label{K1:1} There are $(\bar{\aux},\bar{\omega})\in\R^2\times\Ss$ so that
	$$
	\velocity=\bar{\aux}+\theta\bar{\velocity},\quad\quad
	\aux=\bar{\velocity}+\theta\bar{\aux},
	$$
	where $\bar{\velocity}=\bar{\omega}(\At\bar{\aux}+\im)$, or equivalently,
	$$
	\velocity
	=\LL_{\theta\bar{\omega}}(\bar{\aux})
	:=\bar{\aux}+\theta\bar{\omega}(\At\bar{\aux}+\im),\quad\quad
	\aux=\theta\velocity+(1-\theta^2)(\At\bar{\aux}+\im)\bar{\omega}.
	$$
	\item\label{Ldirection} There is $\bar{\omega}\in\Ss$ so that
	$$(1+\bar{\omega}\theta\At)(\aux-\theta\velocity)=(1-\theta^2)(\At\velocity+\im)\bar{\omega}.$$
	\item\label{Lposition} $z\in\partial\rel$, that is,
	$$|2(1-\theta\At)(\aux-\theta\velocity)+(1-\theta^2)(\At\velocity+\im)|
	=(1-\theta^2)|\At\velocity+\im|.$$
	\item\label{fA} $f(z)=0$, where
	$$f(z):=|2(1-\theta\At)(\aux-\theta\velocity)+(1-\theta^2)(\At\velocity+\im)|-(1-\theta^2)|\At\velocity+\im|.$$
	\item\label{gA} $g(z)=0$, where
	\begin{align*}
	g(z)&:=
	((1-\theta\At)(\aux-\theta\velocity)+(1-\theta^2)(\At\velocity+\im))\cdot(\aux-\theta\velocity)\\
	&\hspace{0.1cm}=(\aux+\At\velocity+\im-\theta(\velocity+\At\aux+\theta\im))\cdot(\aux-\theta\velocity)\\
	&\hspace{0.1cm}=(\aux+\At\velocity+\im)\cdot(\aux-\theta\velocity)-\theta(\velocity+\At\aux+\theta\im)\cdot\aux
	-\theta\det\mathbf{T}_{\At}(z).
	\end{align*}
\end{enumerate}
\end{lemma}
\begin{proof} 
By definition \eqref{def:K1L} a state $z=(\theta,\velocity,\aux)$ belongs to $\Kconstrain^{1,\Lambda}$ if and only if there are $s\in[-1,1]$, $z_1,z_2\in\Kconstrain$ so that $z_1-z_2\in\Lambda$ and
\begin{equation}\label{lemma:K1:1}
z=\tfrac{1+s}{2}z_1+\tfrac{1-s}{2}z_2
=\expected{z}+s\bar{z},
\end{equation}
where $\expected{z}\equiv\tfrac{z_1+z_2}{2}$ and $\bar{z}\equiv\tfrac{z_1-z_2}{2}$. Since $z_j\in\Kconstrain$, we have $|\theta_j|=1$ and $\aux_j=\theta_j\velocity_j$ for $j=1,2$.\\
\indent$\ref{K1:L}\Leftrightarrow\ref{K1:L1}$: Let us assume that $\bar{\theta}=0$ ($\bar{z}\in\varLambda_0$). On the one hand, $\theta_1=\theta_2=\theta$. Hence,
$\aux_j=\theta\velocity_j$ for $j=1,2$, and so $\bar{\aux}=\theta\bar{\velocity}$. On the other hand, $\bar{\velocity}=-\At\bar{\aux}$. Thus, necessarily $\bar{z}=0$ ($z_1=z_2$). Therefore, $\Kconstrain^{1,\varLambda_0}=\Kconstrain$.\\ 
\indent$\ref{K1:L1}\Leftrightarrow\ref{K1:1}$: Now let us assume that $\bar{\theta}\neq 0$ ($\bar{z}\in\varLambda_1$). On the one hand, w.l.o.g.~(relabelling if necessary) we may assume that $\theta_1=-\theta_2=1$.
Hence $\aux_1=\velocity_1$ and $\aux_2=-\velocity_2$, and so $\expected{\aux}=\bar{\velocity}$ and $\bar{\aux}=\expected{\velocity}$. Thus, \eqref{lemma:K1:1} reads as
\begin{equation}\label{lemma:K1:2}
(\theta,\velocity,\aux)
=(0,\bar{\aux},\bar{\velocity})
+s(1,\bar{\velocity},\bar{\aux}).
\end{equation}
On the other hand, there is $\bar{\omega}\in\Ss$ so that $\bar{\velocity}=\bar{\omega}(\At\bar{\aux}+\im).$ Thus, \eqref{lemma:K1:2} reads as
\begin{align*}
\theta&=s,\\
\velocity&=\bar{\aux}+\theta\bar{\velocity}
=L_{\theta\bar{\omega}}(\bar{\aux}),\\
\aux&=\bar{\velocity}+\theta\bar{\aux}
=\theta\velocity+(1-\theta^2)(\At\bar{\aux}+\im)\bar{\omega}.
\end{align*}
\indent$\ref{K1:1}\Leftrightarrow\ref{Ldirection}$:
By definition, the map $\LL_{\theta\bar{\omega}}$ satisfies the identity
\begin{equation}\label{L:1}
\At\LL_{\theta\bar{\omega}}(\bar{\aux})+\im=(1+\bar{\omega}\theta\At)(\At\bar{\aux}+\im).
\end{equation}
This concludes the proof because $\bar{\aux}=\LL_{\theta\bar{\omega}}^{-1}(\velocity)$ and $(1+\bar{\omega}\theta\At)\neq 0$.\\
\indent$\ref{Ldirection}\Leftrightarrow\ref{Lposition}$: 
Although this equivalence can be checked directly by elementary computations, let us give a shorter geometric proof.
For any $b\in\D$ let us consider the automorphism of the shifted disc $\bar{\Sc}$ 
\begin{equation}\label{MobT}
\TT_b(\omega):=\frac{(1-b)\omega}{1+\omega b}.
\end{equation}
This can be expressed in terms of the classical automorphism of the unit disc $\D$ (recall \eqref{shift})
$$\tilde{\varphi}_a(z):=\frac{z-a}{1-a^*z},$$
as $\TT_{b}(\omega)=\shift^{-1}\tilde{\varphi}_{a(b)}(\shift\omega)$
where $a(b)=\frac{b}{2-b}\in\D$. From Complex Analysis it is well-known that $\TT_{b}\in\mathrm{Aut}(\Ss)$ and also $\TT_{b}\in\mathrm{Aut}(\Sc)$.
Thus, \ref{Ldirection} reads as
$$(1-\theta\At)(\aux-\theta\At)
=(1-\theta^2)(\At\velocity+\im)\TT_{\theta\At}(\bar{\omega}).$$
This concludes the proof since $\TT_{\theta\At}\in\mathrm{Aut}(\Ss)$.\\
\indent$\ref{Lposition}\Leftrightarrow\ref{fA}$: Trivial.\\
\indent$\ref{fA}\Leftrightarrow\ref{gA}$: This follows from
\begin{equation}\label{fg}
4(1-\theta\At)g(z)
=|2(1-\theta\At)(\aux-\theta\velocity)+(1-\theta^2)(\At\velocity+\im)|^2-(1-\theta^2)^2|\At\velocity+\im|^2,
\end{equation}
and the fact that $(1-\theta\At)>0$.
\end{proof}

\begin{lemma}\label{lemma:U} Let $z=(\theta,\velocity,\aux)\in(-1,1)\times\R_{\At}^2\times\R^2$ where $\R_{\At}^2:=\{\velocity\in\R^2\,:\,\At\velocity+\im\neq 0\}$.\\ The following are equivalent:
\begin{enumerate}[(a)]
	\addtocounter{enumi}{3}
	\item\label{Ldirection'} There is $\omega\in\Sc$ so that
	$$(1+\omega\theta\At)(\aux-\theta\velocity)=(1-\theta^2)(\At\velocity+\im)\omega.$$
	\item\label{Lposition'} $z\in\rel$, that is,
	$$|2(1-\theta\At)(\aux-\theta\velocity)+(1-\theta^2)(\At\velocity+\im)|
	<(1-\theta^2)|\At\velocity+\im|.$$	
	\item\label{fA'} $f(z)<0$.
	\item\label{gA'} $g(z)<0$.
\end{enumerate}
\end{lemma}
\begin{proof} $\ref{Ldirection'}\Leftrightarrow\ref{Lposition'}$: Analogously to the proof of the equivalence $\ref{Ldirection}\Leftrightarrow\ref{Lposition}$ in Lemma~\ref{lemma:K1}, this follows from the fact that $\TT_{\theta\At}\in\mathrm{Aut}(\Sc)$. $\ref{Lposition'}\Leftrightarrow\ref{fA'}$: Trivial. 
$\ref{fA'}\Leftrightarrow\ref{gA'}$: This follows from \eqref{fg} and $(1-\theta\At)>0$.
\end{proof}

\begin{Rem}\label{Rem:geometry} The equivalences  $\ref{Ldirection}\Leftrightarrow\ref{Lposition}$ are trivial for $\At=0$ because $\TT_{0}=\mathrm{id}$ (cf.~\eqref{MobT}). For a general $|\At|<1$, $\rel_{\At}$ can be understood as $(-1,1)\times\R_{\At}^2\times\Sc$ via the change of variables
\begin{eqnarray*}
\rel_{\At} & \simeq & (-1,1)\times\R_{\At}^2\times\Sc\\
(\theta,\velocity,\aux) & \leftrightarrow & (\theta,\velocity,\omega)
\end{eqnarray*}	
given by
$$\aux=\theta\velocity+(1-\theta^2)(\At\velocity+\im)\bigcirc
\quad\textrm{where}\quad
\bigcirc
=\underbrace{\frac{\omega}{1+\omega\theta\At}}_{\textit{(d)}}
=\underbrace{\frac{\TT_{\theta\At}(\omega)}{1-\theta\At}}_{\textit{(e)}}.$$
Thus, given $z\in\rel_\At$ near to some $z_0\in\Kconstrain^{1,\Lambda_\At}=\partial\rel_\At$, while $\omega\in\Sc$ 
is near to the direction $\bar{\omega}(z_0)\in\Ss=\partial\Sc$ (coupled with $\bar{\aux}=\LL_{\theta\bar{\omega}}^{-1}(\velocity)$) used to construct $z_0$ in Lemma~\ref{lemma:K1}\ref{Ldirection}, the transformation $\TT_{\theta\At}(\omega)$ represents the position of $\aux$ in the ball defined by Lemma~\ref{lemma:U}\ref{Lposition'}.
\end{Rem}

\subsection{Proof of Lemma~\ref{H2rel}}\label{sec:unbounded}
This follows from the below stronger version of Lemma~\ref{H2rel}.

\begin{lemma}\label{lemma:Zseg} There is $d_{\At}>0$ such that, for all $z\in\rel$ there is $\bar{z}\in\Lambda$ with $\bar{\theta}=1$ for which
$$z+\lambda\bar{z}\in\rel,\quad |\lambda|\leq d_A(1-\theta^2).$$
\end{lemma}
\begin{proof} Given $z=(\theta,\velocity,\aux)\in\rel$, let $\bar{z}=(1,\bar{\velocity},\bar{\aux})\in\varLambda_1$, that is, $\bar{\velocity}=\bar{\omega}(\At\bar{\aux}+\im)$ for some $(\bar{\aux},\bar{\omega})\in\R^2\times\Ss$ to be determined. 
Since $\rel$ is open, there is $\epsilon(z,\bar{z},\rel)>0$ so that 
$z_{\lambda}\equiv z+\lambda\bar{z}\in\rel$ for all $|\lambda|\leq\epsilon$, that is, 
\begin{equation}\label{rel:6}
|\theta_\lambda|<1,\quad\quad
\At\velocity_{\lambda}+\im\neq 0,
\end{equation}
and there is $\omega_{\lambda}\in\Sc$ satisfying (Lemma~\ref{lemma:U}\ref{Ldirection'})
\begin{equation}\label{rel:1}
(1+\omega_\lambda\theta_\lambda\At)(\aux_\lambda-\theta_\lambda\velocity_\lambda)
=(1-(\theta_\lambda)^2)(\At\velocity_\lambda+\im)\omega_\lambda,
\end{equation}
for all $|\lambda|\leq\epsilon$.
To prove Lemma~\ref{lemma:Zseg} we must find some $\bar{z}$ making $\epsilon$ big enough, namely $\epsilon(1-\theta^2,\At)$. Roughly speaking, if $z$ is far from $\partial\rel$, $\epsilon$ is controlled easily. Conversely, if $z$ is close to $\partial\rel$, a priori $\epsilon$ is comparable to $\dist{z}{\partial\rel}$, unless we take $\bar{z}$ somehow ``parallel'' to $\partial\rel$. In light of Remark~\ref{Rem:geometry}, it seems suitable to consider $\bar{\aux}=L_{\theta\bar{\omega}}^{-1}(\velocity)$ with $\bar{\omega}\approx\omega_0$ to be determined. Let us see that this choice works. 
We split the proof in two steps. Firstly (\textit{step 1}) we prove the statement by assuming a claim. 
Secondly (\textit{step 2}) this claim is proved by elementary computations.\\
\indent\textit{Step 1. Claim}: Let us take $\bar{\aux}=L_{\theta\bar{\omega}}^{-1}(\velocity)$ with $\bar{\omega}\in\Ss$ to be determined. Then, \eqref{rel:6} holds for all $|\lambda|\leq\frac{1}{2}(1-\theta^2)$ and \eqref{rel:1} is equivalent to
\begin{equation}\label{rel:4}
\lambda\alpha|\shift\bar{\omega}-\shift\omega|
<(1-\theta^2)(1-|\shift{\omega}|^2),
\end{equation}
where $\omega\equiv\omega_0$ and $|\alpha|\leq \alpha_{\At}$ for some constant $\alpha_{\At}>0$. 
We shall prove this claim in the \textit{step 2}.\\
Assume that this claim is true.
Hence, if we make the change of variables $\lambda=d(1-\theta^2)$ for $d\in\R$, \eqref{rel:4} reads as
\begin{equation}\label{rel:5}
d\alpha|\shift\bar{\omega}-\shift\omega|<(1-|\shift\omega|^2).
\end{equation} 
If $|\shift\omega|\leq\frac{1}{2}$ ($z$ is far from $\partial\rel$) we take $\bar{\omega}=0\in\Ss$ and then \eqref{rel:5} holds for every $|d|\leq\frac{1}{2\alpha_\At}$.\\
If $\frac{1}{2}<|\shift\omega|<1$ ($z$ is close to $\partial\rel$) we take $\shift\bar{\omega}=\frac{\shift\omega}{|\shift\omega|}$ and then \eqref{rel:5} reads as
$$d\alpha<(1+|\shift{\omega}|),$$
which holds for every $|d|\leq \frac{1}{\alpha_\At}$. Therefore, we can take $d_\At=\frac{1}{2\alpha_{\At}}$.\\
\indent\textit{Step 2. Proof of the claim}:
Since $\bar{\velocity}=\bar{\omega}(\At\bar{\aux}+\im)$ and $\bar{\aux}=\LL_{\theta\bar{\omega}}^{-1}(\velocity)$,
Lemma~\ref{lemma:K1}\ref{K1:1} and \eqref{L:1} yield
$$
\bar{\velocity}=
\bar{\omega}\frac{\At\velocity+\im}{1+\bar{\omega}\theta\At},\quad\quad
\bar{\aux}
=\velocity-\theta\bar{\velocity}.
$$
Let us expand the factors of \eqref{rel:1} in terms of $\lambda$. They are
\begin{subequations}
\label{rel:7}
\begin{align}
\aux_\lambda-\theta_\lambda\velocity_\lambda
&=(\aux-\theta\velocity)+\lambda(\bar{\aux}-(\theta\bar{\velocity}+\velocity))-\lambda^2\bar{\velocity}
=(\aux-\theta\velocity)+(\theta^2-(\theta_\lambda)^2)\bar{\velocity},\label{rel:7:1}\\
\At\velocity_\lambda+\im
&=(\At\velocity+\im)+\lambda\bar{\velocity}
=(\At\velocity+\im)\frac{1+\bar{\omega}\theta_\lambda\At}{1+\bar{\omega}\theta\At}.\label{rel:7:2}
\end{align}
\end{subequations}
Since $z\in\rel$, we have $|\theta|<1$ and $\At\velocity+\im\neq 0$. Then, by \eqref{rel:7:2}: $|\theta_{\lambda}|<1\Rightarrow\At\velocity_{\lambda}+\im\neq 0$. Therefore, \eqref{rel:6} is equivalent to $|\theta+\lambda|<1$, and this holds for all $|\lambda|\leq\frac{1}{2}(1-\theta^2)$.\\
By \eqref{rel:7}, if we multiply \eqref{rel:1} by $(1+\bar{\omega}\theta\At)(1+\omega\theta\At)/(\At\velocity+\im)$, we get
\begin{equation}\label{rel:2}
\begin{split}
(1+\omega_\lambda\theta_\lambda\At)
((1+\bar{\omega}\theta\At)(1-\theta^2)\omega+(1+\omega\theta\At)(\theta^2-(\theta_\lambda)^2)\bar{\omega})\\
=(1-(\theta_\lambda)^2)(1+\omega\theta\At)(1+\bar{\omega}\theta_\lambda\At)\omega_\lambda.
\end{split}
\end{equation}
Hence, by applying the following identities
\begin{align*}
(1+\bar{\omega}\theta\At)(1-\theta^2)\omega+(1+\omega\theta\At)(\theta^2-(\theta_\lambda)^2)\bar{\omega}
&=(1-\theta^2)(\omega-\bar{\omega})
+(1+\omega\theta\At)(1-(\theta_\lambda)^2)\bar{\omega},\\
(1+\bar{\omega}\theta_\lambda\At)\omega_\lambda
&=(\omega_\lambda-\bar{\omega})+(1+\omega_\lambda\theta_\lambda\At)\bar{\omega},
\end{align*}
\eqref{rel:2} reads as 
\begin{equation}\label{rel:3}
(1-\theta^2)(1+\omega_\lambda\theta_\lambda\At)(\bar{\omega}-\omega)=
(1-(\theta_\lambda)^2)(1+\omega\theta\At)(\bar{\omega}-\omega_\lambda).
\end{equation}
Since (recall \eqref{shift}) $w=\frac{1}{2}(\shift w-1)$ for all $w\in\R^2$, \eqref{rel:3} reads as
$$(1-\theta^2)((2-\theta_{\lambda}\At)+\theta_{\lambda}\At\shift\omega_{\lambda})(\shift\bar{\omega}-\shift\omega)
=(1-(\theta_\lambda)^2)((2-\theta\At)+\theta\At\shift\omega)(\shift\bar{\omega}-\shift\omega_\lambda),$$
or equivalently, $\zeta\shift\omega_{\lambda}=\eta$ where we have abbreviated
\begin{align*}
\zeta&\equiv(1-(\theta_\lambda)^2)((2-\theta\At)+\theta\At\shift\omega)+(1-\theta^2)\theta_\lambda\At(\shift\bar{\omega}-\shift\omega),\\
\eta&\equiv(1-(\theta_\lambda)^2)((2-\theta\At)+\theta\At\shift\omega)\shift\bar{\omega}
-(1-\theta^2)(2-\theta_\lambda\At)(\shift\bar{\omega}-\shift\omega).
\end{align*}
In this way: $|\eta|<|\zeta|\Rightarrow\omega_{\lambda}\in\Sc$.
Let us write the inequality $|\eta|^2<|\zeta|^2$. Since $|\shift\bar{\omega}|=1$, the term $(1-(\theta_\lambda)^2)^2|(2-\theta\At)+\theta\At\shift\omega|^2$ is cancelled. Hence, by reordering the remainder terms, the inequality $|\eta|^2<|\zeta|^2$ is equivalent to
\begin{equation}\label{KL:8}
\begin{split}
(1-\theta^2)((2-\theta_{\lambda}\At)^2-(\theta_{\lambda}\At)^2)|\shift\bar{\omega}-\shift\omega|^2\\
<2(1-(\theta_\lambda)^2)(((2-\theta\At)+\theta\At\shift\omega)((2-\theta_{\lambda}\At)\shift\bar{\omega}+\theta_{\lambda}\At))\cdot(\shift\bar{\omega}-\shift\omega),
\end{split}
\end{equation}
where we have eliminated a factor $(1-\theta^2)>0$. Notice that \eqref{KL:8} can be written as $p(\lambda)<0$ for some (3-degree) polynomial $p$ in $\lambda$. In particular, \eqref{KL:8} can be written as
\begin{equation}\label{KL:9}
\lambda\left(\int_0^1\partial_{\lambda}q(s\lambda)\dif s\right)\cdot(\shift{\bar{\omega}}-\shift{\omega})<-p(0),
\end{equation}
where $p(\lambda)=q(\lambda)\cdot(\shift{\bar{\omega}}-\shift{\omega})$, that is,
\begin{align*}
q(\lambda)
\equiv&\,(1-\theta^2)((2-\theta_{\lambda}\At)^2-(\theta_{\lambda}\At)^2)(\shift{\bar{\omega}}-\shift{\omega})\\
&-2(1-(\theta_\lambda)^2)((2-\theta\At)+\theta\At\shift\omega)((2-\theta_{\lambda}\At)\shift\bar{\omega}+\theta_{\lambda}\At).
\end{align*}
On the one hand, since $|\lambda|,|\theta|,|\At|,|\shift\omega|,|\shift\bar{\omega}|\leq 1$ we can bound
\begin{equation}\label{KL:lhs}
\left|\int_0^1\partial_{\lambda}q(s\lambda)\dif s\right|
\leq C,
\end{equation}
for some constant $C>0$.
On the other hand, $-p(0)=(1-\theta^2)\beta$ where 
we have abbreviated
$$\beta\equiv 2(((2-\theta\At)+\theta\At\shift\omega)((2-\theta\At)\shift\bar{\omega}+\theta\At))\cdot(\shift{\bar{\omega}}-\shift{\omega})
-((2-\theta\At)^2-(\theta\At)^2)|\shift{\bar{\omega}}-\shift{\omega}|^2.$$
Remarkably, using $|T\bar{\omega}|=1$ and abbreviating $a\equiv\frac{\theta\At}{2-\theta\At}$, this term can be greatly simplified
\begin{equation}\label{KL:rhs}
\begin{split}
\beta
&=(2-\theta\At)^2(2((1+a\shift{\omega})(\shift{\bar{\omega}}+a))\cdot(\shift{\bar{\omega}}-\shift{\omega})-(1-a^2)|\shift{\bar{\omega}}-\shift{\omega}|^2)\\
&=(2-\theta\At)^2((1+a^2)(\shift{\bar{\omega}}+\shift{\omega})+2a(1+\shift{\bar{\omega}}\shift{\omega}))\cdot(\shift{\bar{\omega}}-\shift{\omega})\\
&=(2-\theta\At)^2|1+a\shift{\bar{\omega}}|^2(1-|\shift{\omega}|^2)\\
&=4|1+\bar{\omega}\theta\At|^2(1-|\shift{\omega}|^2).
\end{split}
\end{equation}
By applying \eqref{KL:lhs}\eqref{KL:rhs} on \eqref{KL:9}, we deduce \eqref{rel:4} with 
$$\alpha\equiv\frac{1}{4|1+\bar{\omega}\theta\At|^2}\left(\int_0^1\partial_{\lambda}q(s\lambda)\dif s\right)\cdot
\frac{(\shift{\bar{\omega}}-\shift{\omega})}{|\shift{\bar{\omega}}-\shift{\omega}|},$$
 which satisfies $|\alpha|\leq \frac{C}{4(1-|\At|)^2}.$
\end{proof}

\subsection{Proof of Lemma~\ref{H2bounded}}\label{sec:bounded}
As in \cite{RIPM}, the relaxed set $\rel$ is unbounded, thereby preventing from constructing $L^\infty$-solutions from the h-principle applied to $\rel$-valued subsolutions. In order to
find bounded subsets of $\rel$ satisfying (H2) we have to restrict $\Kconstrain$ somehow. In \cite{RIPM} ($\At=0$) 
%after symmetrizing (IPM2-$3_0$) via the change of variables $\mathbf{v}=2\velocity+\theta\im$,
Sz\'ekelyhidi computed explicitly the $\Lambda_0$-convex hull of
$$\Kconstrain_{M}:=\{z\in\Kconstrain\,:\,|\underbrace{2\velocity+\theta\im}_{\equiv\mathbf{v}}|\leq M\}=\{z\in\Kconstrain\,:\,\underbrace{4\velocity\cdot(\velocity+\theta\im)}_{\Lambda_0\textrm{-linear}}\leq M^2-1\},$$
for any $M>1$ (notice $\Kconstrain_{M}\Subset\Kconstrain$) which is given by the following 4 inequalities: 
\begin{subequations}
\label{idL1}
\begin{align}
|2(\aux-\theta\velocity)+(1-\theta^2)\im|&<(1-\theta^2),\label{idL1:1}\\
4\velocity\cdot(\velocity+\theta\im)&<M^2-1,\label{idL1:2}\\
|2(\aux-\velocity)+(1-\theta)\im|&<M(1-\theta),\label{idL1:3}\\
|2(\aux+\velocity)+(1+\theta)\im|&<M(1+\theta).\label{idL1:4}
\end{align}
\end{subequations}
As observed in \cite{RIPM}, these inequalities are linked
by the following identity: 
\begin{subequations}
\label{idL}
\begin{align}
&(1-\theta^2)^2-|2(\aux-\theta\velocity)+(1-\theta^2)\im|^2\label{idL:1}\\
&+(1-\theta^2)(M^2-1-4\velocity\cdot(\velocity+\theta\im))\label{idL:2}\\
&=\frac{1+\theta}{2}(M^2(1-\theta)^2-|2(\aux-\velocity)+(1-\theta)\im|^2)\label{idL:3}\\
&\hspace{0.04cm}+\frac{1-\theta}{2}(M^2(1+\theta)^2-|2(\aux+\velocity)+(1+\theta)\im|^2),\label{idL:4}
\end{align}
\end{subequations}
which is indeed crucial to prove (H2). 

\begin{Rem}\label{Rem:asymmetry} 
In \cite{RIPM} Sz\'ekelyhidi introduces the smart (linear) change of variables $(\theta,\mathbf{v},\mathbf{n})=(\theta,2\velocity+\theta\im,2\aux+\im)$, which simplifies significantly the computations and inequalities in \eqref{idL1}. Under this transformation: 1) the wave cone reads as $\Lambda_0=\{\bar{z}\in\R^5\,:\,|\bar{\theta}|=|\bar{\mathbf{v}}|\}$ because (IPM2-$3_0$) become symmetric, 2) the geometry of $\Kconstrain$ is preserved (given $|\theta|=1$: $\aux=\theta\velocity\Leftrightarrow\mathbf{n}=\theta\mathbf{v}$). After this, Sz\'ekelyhidi computed the $\Lambda_0$-convex hull of $\Kconstrain_M=\{z\in\Kconstrain\,:\,|\mathbf{v}|\leq M\}$.\\ 
\indent For a general $|\At|<1$, the corresponding change of variables that keeps 1) and 2) is $(\theta,\mathbf{v},\mathbf{n})=(\theta,2\velocity+\At\aux+\theta\im,(2+\theta\At)\aux+\im)$, which is not linear in $\mathbf{n}$ for $\At\neq 0$, thereby hampering the plane wave analysis. Thus, for $\At\neq 0$, although  $\mathbf{v}=2\velocity+\At\aux+\theta\im$ symmetrizes (IPM2-$3_\At$), any linear change of variables in $\mathbf{n}$ messes the simplicity of $\Kconstrain$ up. This is why we have chosen not to make a change variables in this work.\\
\indent In this regard, for $\At\neq 0$ it is not evident what restriction of $\Kconstrain$ may return a simple $\Lambda_{\At}$-convex hull as in \eqref{idL1}. To overcome this drawback, inspired by \eqref{idL},
instead of restricting $\Kconstrain$ first,
we start trying to extend properly the identity \eqref{idL} to $|\At|<1$, with the hope that this will reveal the analogous inequalities to \eqref{idL1} that describe the $\Lambda_{\At}$-convex hull of some restriction of $\Kconstrain$. Fortunately, this is the case.
\end{Rem}

\begin{lemma} For every $M\in\R$ and $z=(\theta,\velocity,\aux)\in\set$,
\begin{subequations}
\label{id1}
\begin{align}
&\frac{1}{1-\theta\At}((1-\theta^2)^2|\At\velocity+\im|^2-|2(1-\theta\At)(\aux-\theta\velocity)+(1-\theta^2)(\At\velocity+\im)|^2)\label{id1:1}\\
&+(1-\theta^2)(M^2-1-4\velocity\cdot(\velocity+\At\aux+\theta\im+\At\im))\label{id1:2}\\
&=\frac{1+\theta}{2}((M^2-\At)(1-\theta)^2-(1-\At)|2(\aux-\velocity)+(1-\theta)\im|^2)\label{id1:3}\\
&\hspace{0.04cm}+\frac{1-\theta}{2}((M^2+\At)(1+\theta)^2-(1+\At)|2(\aux+\velocity)+(1+\theta)\im|^2).\label{id1:4}
\end{align} 
\end{subequations}
\end{lemma}
\begin{proof} First notice that, by \eqref{fg}, we have $\eqref{id1:1}=-4g(z)$. On the one hand,
\begin{align*}
\eqref{id1:1}+\eqref{id1:2}
&=4(\theta(\velocity+\At\aux+\theta\im)-(\aux+\At\velocity+\im))\cdot\aux+4\theta(\aux+\At\velocity+\im)\cdot\velocity\\
&+(1-\theta^2)(M^2-1-4\At\velocity_2)-4\velocity\cdot(\velocity+\At\aux+\theta\im)\\
&=-4(1-\theta\At)(|\aux|^2+|\velocity|^2)-8(\At-\theta)\aux\cdot\velocity\\
&+(1-\theta^2)(M^2-1-4(\aux_2+\At\velocity_2)).
\end{align*}
On the other hand,
\begin{align*}
\eqref{id1:3}+\eqref{id1:4}
&=\frac{1+\theta}{2}((M^2-1)(1-\theta)^2-4(1-\At)(|\aux-\velocity|^2+(1-\theta)(\aux_2-\velocity_2)))\\
&+\frac{1-\theta}{2}((M^2-1)(1+\theta)^2-4(1+\At)(|\aux+\velocity|^2+(1+\theta)(\aux_2+\velocity_2)))\\
&=-2((1+\theta)(1-\At)|\aux-\velocity|^2+(1-\theta)(1+\At)|\aux-\velocity|^2)\\
&+(1-\theta^2)(M^2-1-2(1-\At)(\aux_2-\velocity_2)-2(1+\At)(\aux_2+\velocity_2)).
\end{align*}
This concludes the proof.
\end{proof}
Observe that \eqref{id1} generalizes \eqref{idL}.
For any $M>1$, we consider the open set $\rel_{\At,M}$ of states $z\in\set$ given by the following 4 inequalities:
\begin{subequations}
\label{id2}
\begin{align}
|2(1-\theta\At)(\aux-\theta\velocity)+(1-\theta^2)(\At\velocity+\im)|&<(1-\theta^2)|\At\velocity+\im|,\label{id2:1}\\
4\velocity\cdot(\velocity+\At\aux+\theta\im+\At\im)&<M^2-1,\label{id2:2}\\
|2(\aux-\velocity)+(1-\theta)\im|&<M_{-\At}(1-\theta),\label{id2:3}\\
|2(\aux+\velocity)+(1+\theta)\im|&<M_{+\At}(1+\theta),\label{id2:4}
\end{align}
\end{subequations}
where
$$M_{\pm\At}\equiv\sqrt{\frac{M^2\pm\At}{1\pm\At}}.$$
By analogy with \cite{RIPM}, \eqref{id1} suggests that $\rel_{\At,M}$ is the interior of the $\Lambda_\At$-convex hull of
$$\Kconstrain_{\At,M}
:=\{z\in\Kconstrain\,:\,|2\velocity+\theta\im|\leq M_{\theta\At}\}
=\{z\in\Kconstrain\,:\,B_{\At}(z)\leq M^2-1\},$$
where  we have abbreviated
\begin{equation}\label{B}
B_{\At}(z):=|\mathbf{b}_{\At}(z)+\velocity|^2-|\mathbf{b}_{\At}(z)-\velocity|^2=4\velocity\cdot\mathbf{b}_{\At}(z),
\end{equation}
and 
$$\mathbf{b}_{\At}(z):=\velocity+\At\aux+\theta\im+\At\im.$$
Observe that $\Kconstrain_{0,M}=\Kconstrain_{M}$. 
In Section \ref{sec:Lambdahull} we shall prove that both $\Kconstrain^{lc,\Lambda_{\At}}=\bar{\rel}_{\At}$ and $(\Kconstrain_{\At,M})^{lc,\Lambda_{\At}}=\bar{\rel}_{\At,M}$.  Now, let us continue with the proof of Lemma~\ref{H2bounded}. Thus, from now on we shall omit the tag ``$\At$'' wherever we do not need to distinguish between the cases $\At=0$ and $\At\neq 0$.\\
\indent Firstly, let us check that $\rel_M$ is indeed bounded.
\begin{lemma}\label{Linf} Let $M>1$. The set $\rel_M$ is bounded.
\end{lemma}
\begin{proof} Given $z\in\rel$ there is $\omega\in\Sc$ so that \eqref{L2:1} holds. In particular,
$$\At\aux+\theta\im=\frac{\theta+\omega\At}{1+\omega\theta\At}(\At\velocity+\im).$$
Then, by applying \eqref{L2:2}, we have $|\At\aux+\theta\im|\leq|\At\velocity+\im|\leq|\At||\velocity|+1$. Hence, \eqref{id2:2}\eqref{B} imply
$$4|\velocity|^2=B(z)-4\velocity\cdot(\At\aux+\theta\im+\At\im)
<M^2-1+4|\velocity|(|\At||\velocity|+1+|\At|),$$
and so
$$4((1-|\At|)|\velocity|-(1+|\At|))|\velocity|< M^2-1.$$
Thus, necessarily
$$|\velocity|<\frac{(1+|\At|)+\sqrt{(1+|\At|)^2+(1-|\At|)(M^2-1)}}{2(1-|\At|)}.$$
Finally, recall that $\aux$ is controlled by \eqref{L2:4}.
\end{proof}

Secondly, let us show that these $\rel_M$'s contain simpler sets as stated in Lemma~\ref{H2bounded}.

\begin{lemma} For any $R>0$ there is $M>1$ so that
$$\{z\in\rel\,:\, |\velocity|< R\}\subset\rel_M.$$
\end{lemma}
\begin{proof} Let $z=(\theta,\velocity,\aux)\in\rel$ with $|\velocity|< R$. By Lemma~\ref{lemma:U}\ref{Ldirection'}, there is $\omega\in\Sc$ so that
$$\aux=\theta\velocity+(1-\theta^2)\frac{(\At\velocity+\im)\omega}{1+\omega\theta\At}.$$
Thus, for \eqref{id2:3}\eqref{id2:4} we have
$$|2(\aux\pm\velocity)+(1\pm\theta)\im|
=(1\pm\theta)\left|\pm 2\velocity+\im+(1\mp\theta)\frac{(\At\velocity+\im)\omega}{1+\omega\theta\At}\right|
\leq(1\pm\theta)C_{\pm},$$
for some constant $C_{\pm}(\At,R)>0$.
Concerning \eqref{id2:2} we have
$$1+B(z)\leq C,$$
for some constant $C(\At,R)>0$.
Hence, since there is $M(\At,R)>1$ satisfying $C_{\pm}\leq M_{\pm}$ and $C\leq M^2$, we have $z\in\rel_M$.
\end{proof}

Finally, the following lemma completes the proof of Lemma~\ref{H2bounded}. 

\begin{Rem}\label{Rem:M*}
The pinch singularity $\At\velocity+\im=0$ becomes further complicated for $\rel_{\At,M}$ because the new inequalities \eqref{id2:2}-\eqref{id2:4} can interfere with it for the particular value (cf.~\ref{M*singularity})
\begin{equation}\label{M*}
M_*(A):=\sqrt{1+4\left(\frac{1}{\At^2}-1\right)}.
\end{equation}
Notice that $M_*$ is symmetric and strictly decreasing on $(0,1]$ with $M_*(0)=+\infty$ and $M_*(1)=1$.
For simplicity we shall omit this case.
\end{Rem}

\begin{lemma}\label{lemma:ZMseg} Let $1<M\neq M_*(\At)$. 
The set $\rel_M$ satisfies $\mathrm{(H2)}$.
\end{lemma}
\begin{proof} 
Given $(\theta,\velocity)\in(-1,1)\times\R^2$ we consider the subsets of $\R^2$
\begin{align*}
\B(\theta,\velocity)&:=\{\aux\in\R^2\,:\,|2(1-\theta\At)(\aux-\theta\velocity)+(1-\theta^2)(\At\velocity+\im)|<(1-\theta^2)|\At\velocity+\im|\},\\
\HP(\theta,\velocity)&:=\{\aux\in\R^2\,:\,4\velocity\cdot(\velocity+\At\aux+\theta\im+\At\im)<M^2-1\},\\
\B_-(\theta,\velocity)&:=\{\aux\in\R^2\,:\,|2(\aux-\velocity)+(1-\theta)\im|<M_{-}(1-\theta)\},\\
\B_+(\theta,\velocity)&:=\{\aux\in\R^2\,:\,|2(\aux+\velocity)+(1+\theta)\im|<M_{+}(1+\theta)\}.
\end{align*}
By definition, a state $z=(\theta,\velocity,\aux)\in(-1,1)\times\R^2\times\R^2$ belongs to $\rel$ if and only if $\aux$ belongs to the open ball $\rel(\theta,\velocity):=\B(\theta,\velocity)$. Similarly, $z$ belongs to the bounded subset $\rel_M$ if and only if $\aux$ belongs to  $\rel_M(\theta,\velocity):=(\B\cap\HP\cap\B_-\cap\B_+)(\theta,\velocity)$.
Notice that $\B_-(\theta,\velocity)$ and $\B_+(\theta,\velocity)$ are (open) balls. 
The geometry of $\HP_{\At}(\theta,\velocity)$ depends on $\At$ (cf.~Fig.~\ref{fig:geo}). On the one hand, for $\At=0$
the condition defining $\HP_0$ only depends on $(\theta,\velocity)$, namely $\velocity$ must belong to the (open) ball
$$\mathscr{B}(\theta):=\{\velocity\in\R^2\,:\,|2\velocity+\theta\im|^2<M^2-(1-\theta^2)\},$$
i.e.~$\HP_0(\theta,\velocity)=\R^2$ (or $\emptyset$) if $\velocity$ belongs (or not) to $\mathscr{B}(\theta)$.
On the other hand, for $\At\neq 0$, $\HP_{\At}(\theta,\velocity)$ is an (open) half-plane (except $\HP_{\At}(\theta,0)=\R^2$).\\

\indent In order to help better understand the set $\rel_{\At,M}$ we provide several pictures (Fig.~\ref{geo:0}-\ref{geo:M}) of the slices $\rel_{\At,M}(\theta,\velocity)$, for some fixed $\At$, $M$, $\theta$, and different $\velocity$'s moving parallel to the real and imaginary axis. By symmetry ($\rel_{\At,M}(\theta,-\velocity^*)=-\rel_{\At,M}(\theta,\velocity)^*$) it is enough to consider $\Re\velocity\geq 0$.\\
We differentiate three cases: 1) $A=0$, 2) $0<|\At|<1$ coupled with either 2.1) $M>M_*(A)$ or 2.2) $M<M_*(A)$ (cf.~\eqref{M*}).\\

\begin{figure}[h!]
	\centering
	\begin{subfigure}{.4\textwidth}
		\centering
		\includegraphics[height=5.6cm]{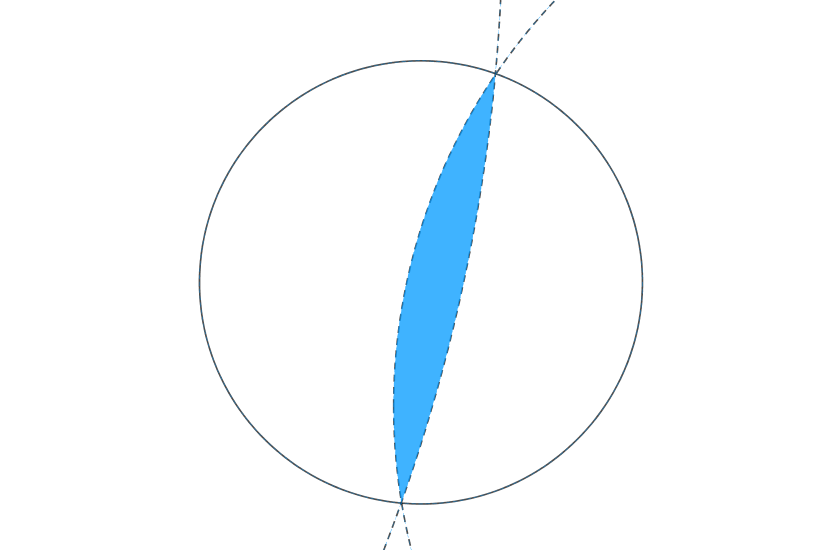}
	\end{subfigure}
	\begin{subfigure}{.5\textwidth}
		\centering
		\includegraphics[height=5.6cm]{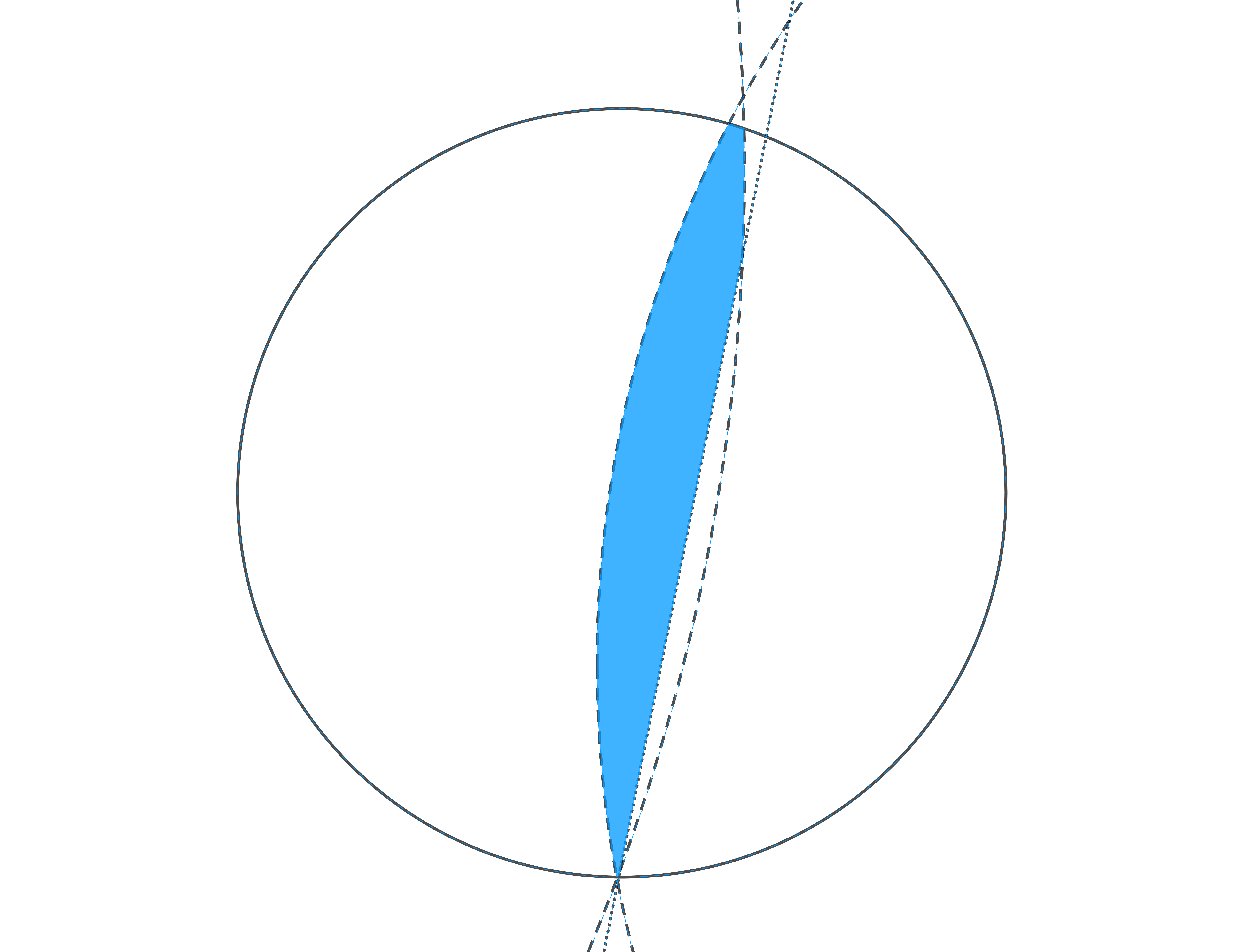}
	\end{subfigure}
	\caption{GeoGebra plot of the region $\rel_{\At,M}(\theta,\velocity)$ (blue) for some $(\theta,\velocity)\in(-1,1)\times\R^2$, $M>1$,
	$\At=0$ (left) and $0<|\At|<1$ (right), where we have added the circles $\partial\B(\theta,\velocity)$ (solid),
	$\partial\B_{-}(\theta,\velocity)$, $\partial\B_{+}(\theta,\velocity)$ (dashed) and, for $\At\neq 0$ (right),  the line $\partial\HP(\theta,\velocity)$ (dotted).}
	\label{fig:geo}
\end{figure}

\indent 1) Let $\At=0$. In this case, the region $\rel_{0,M}(\theta,\velocity)$ does not collapse as $\velocity$ tends to $\partial\mathscr{B}(\theta)$ (cf.~Fig.~\ref{geo:0}). In fact, $\rel_{0,M}(\theta,\velocity)$ collapses if and only if $|\theta|\uparrow 1$ (i.e.~$z$ tends to $\Kconstrain$). In particular, as noted in \cite{RIPM}, $\partial\rel_{0,M}\setminus\Kconstrain$ is locally the graph of a Lipschitz function.

\begin{figure}[h!]
	\centering
	\begin{tabular}{cccc}
		\frame{\includegraphics[width=3.39cm]{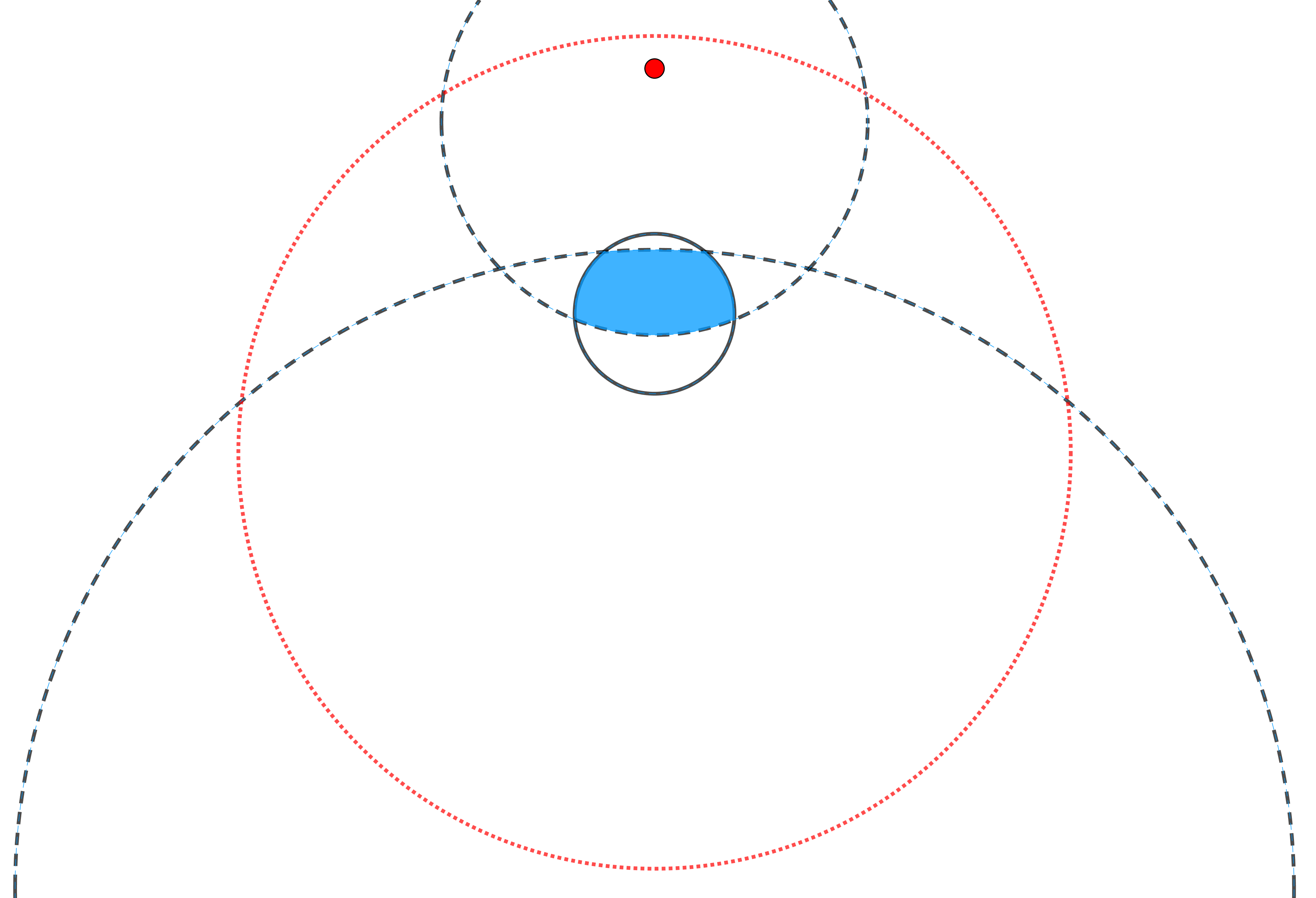}}& & & \\
		$\velocity=\tfrac{1}{2}\theta\im+1.8\im$ & & & \\[0.2cm]
		\frame{\includegraphics[width=3.39cm]{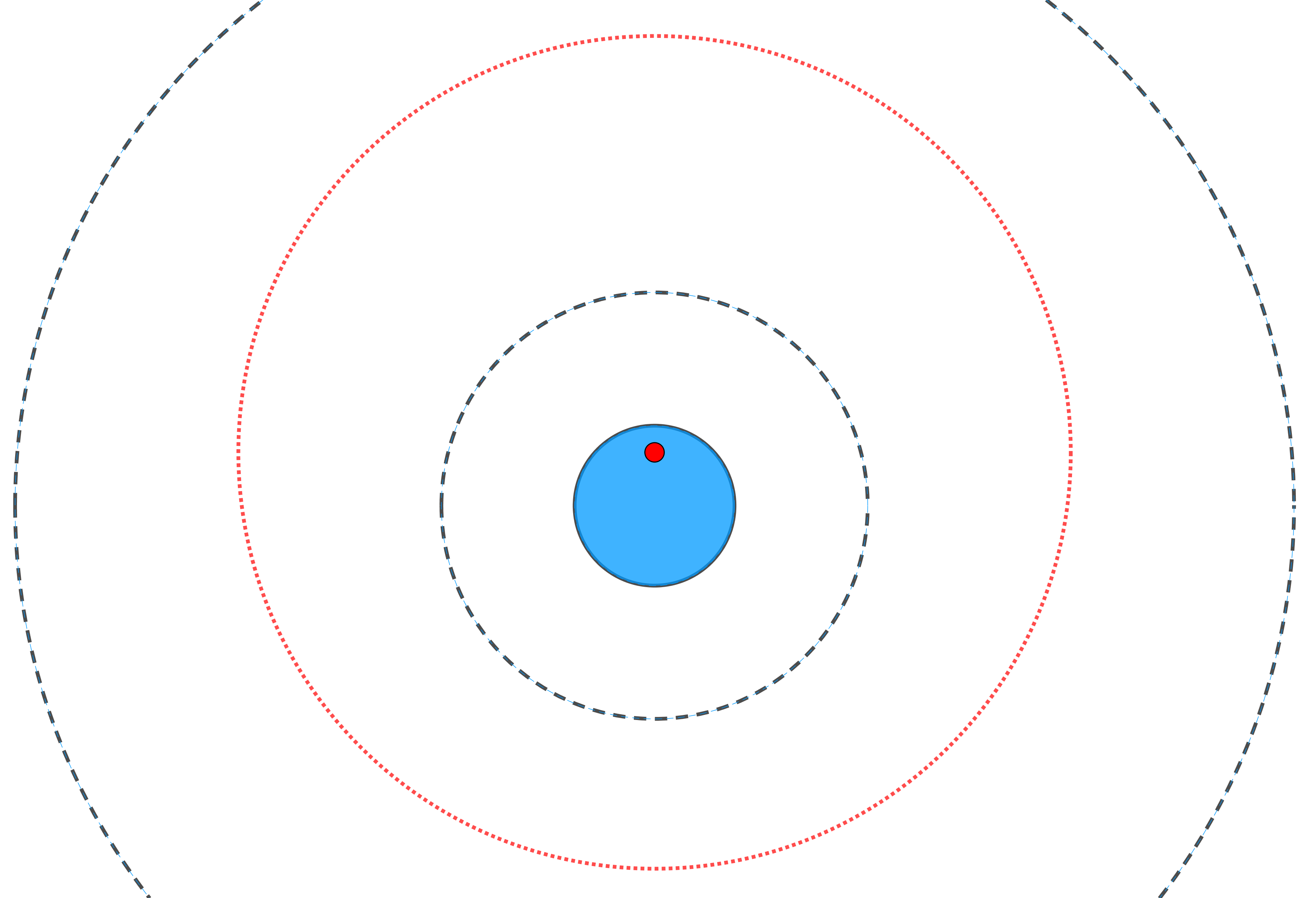}}&
		\frame{\includegraphics[width=3.39cm]{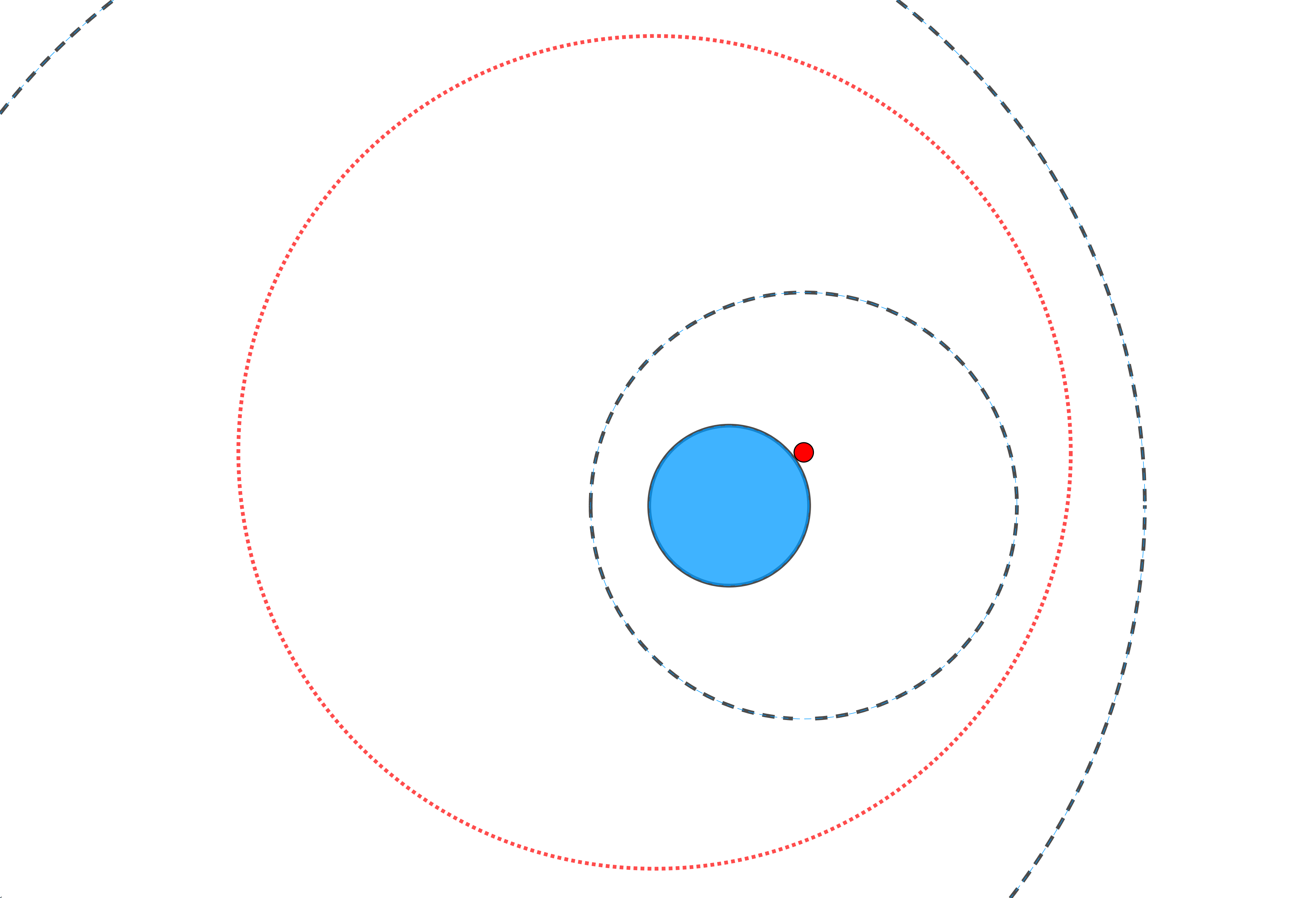}}&
		\frame{\includegraphics[width=3.39cm]{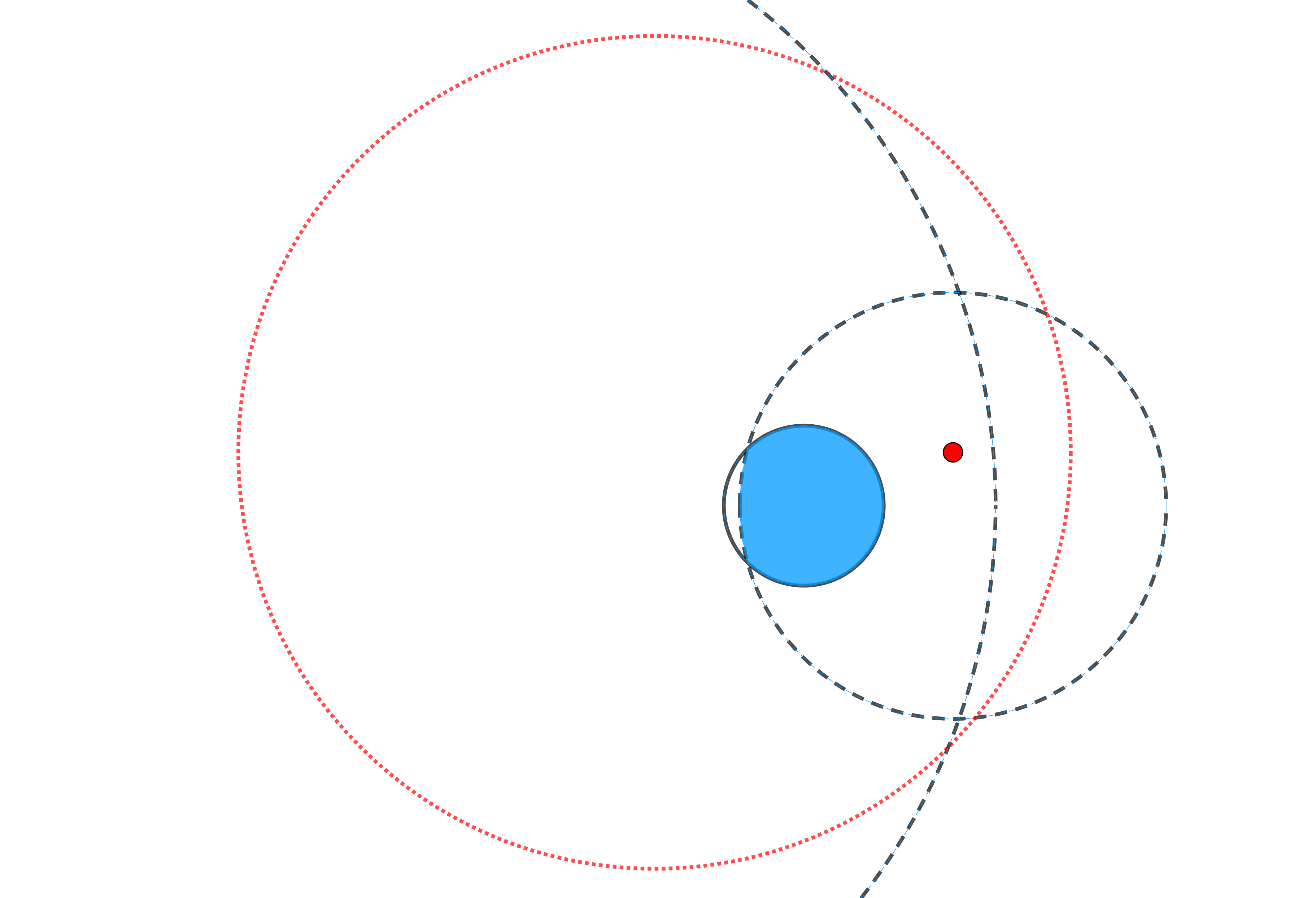}}&
		\frame{\includegraphics[width=3.39cm]{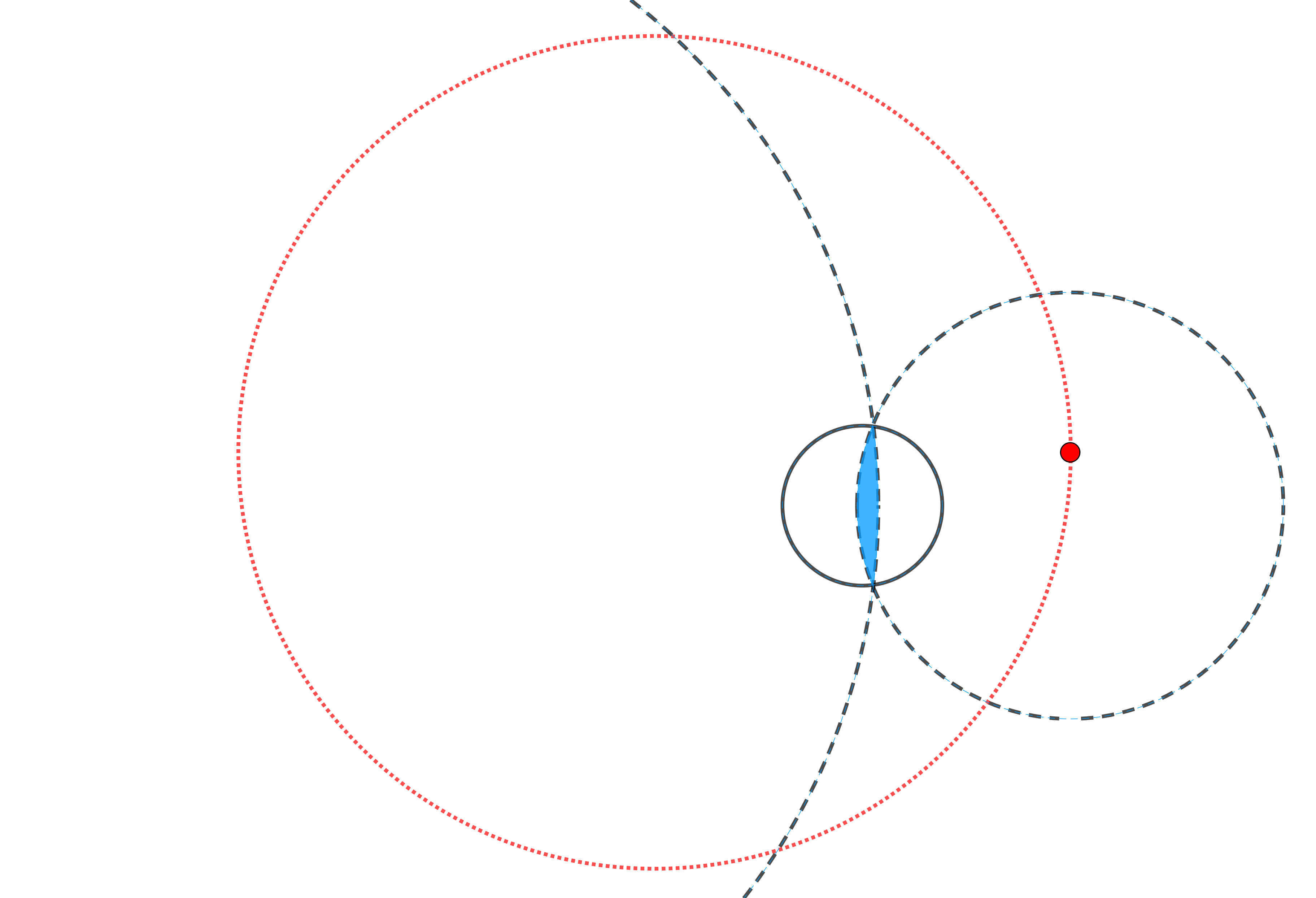}}\\
		$\velocity=\tfrac{1}{2}\theta\im$ & $\velocity=\tfrac{1}{2}\theta\im+0.7$ & $\velocity=\tfrac{1}{2}\theta\im+1.4$ & $\velocity=\tfrac{1}{2}\theta\im+1.95$ \\[0.2cm]
		\frame{\includegraphics[width=3.39cm]{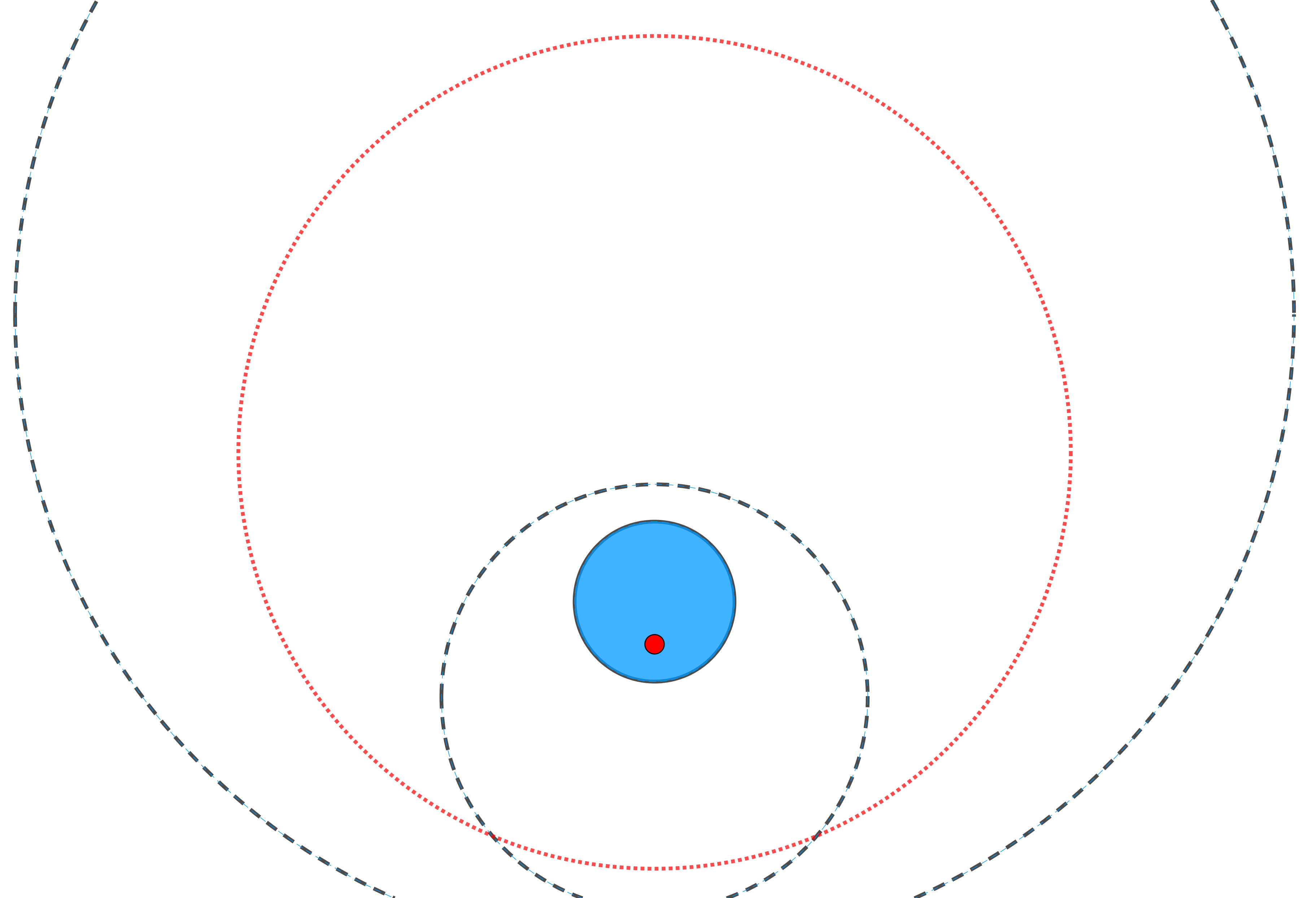}} \\
		$\velocity=\tfrac{1}{2}\theta\im-0.9\im$ & & & 
	\end{tabular}
	\caption{Plots of $\rel_{0,M}(\theta,\velocity)$ (cf.~Fig.~\ref{fig:geo}-left) for $\At = 0$, $M=4$, $\theta=\tfrac{1}{2}$ and different $\velocity$'s (red point) inside the circle $\partial\mathscr{B}(\theta)$ (red dotted).}
	\label{geo:0}
\end{figure}

Since the case $\At=0$ is proved in \cite{RIPM}, from now on we focus on $0<|\At|<1$.\\
\indent 2) Let $0<|\At|<1$. 
On the one hand, the half-plane $\HP(\theta,\velocity)$ causes that $\rel_{\At,M}(\theta,\velocity)$ collapses as $|\velocity|$ grows, in contrast to the case $\At=0$ (cf.~the last column of Fig.~\ref{geo:0} and \ref{geo:MM}).
On the other hand, we have to deal with the pinch singularity $\At\velocity+\im=0$. Given $\gamma>0$ let us denote $S_{\gamma}:=\{z\in\bar{\rel}_{\At}\,:\,|\At\velocity+\im|\leq\gamma\}$. The set $S_0$ ($\gamma=0$) satisfies the following property. Let $(\theta,\velocity,\aux)\in S_0$ with $|\theta|<1$, i.e.~$\At\velocity+\im=0$ and so $\aux=\theta\velocity$. Then, it is straightforward to check that, for any $\square=\HP,\B_{-},\B_{+}$:
\begin{equation}\label{M*singularity}
\aux\in\partial\square(\theta,\velocity)
\quad\Leftrightarrow\quad
M=M_*(A).
\end{equation}
Thus, for the particular value $M=M_*(A)$, the pinch singularity $S_0$ of $\rel_{\At}$ lies in the boundary of all the other new inequalities \eqref{id2:2}-\eqref{id2:4} defining $\rel_{\At,M}$.  For simplicity we omit this case.\\

\indent 2.1) Let $M>M_*(A)$. Then $\rel_{\At,M}(\theta,\velocity)=\B(\theta,\velocity)$ in a neighbourhood of $\velocity=-\tfrac{1}{\At}\im$ (cf.~Fig.~\ref{geo:MM}). Therefore, there is $\gamma(\At,M)>0$ so that $S_{\gamma}\cap\rel_{\At,M}=S_{\gamma}\cap\rel_{\At}$ and thus the $\Lambda$-directions from Lemma~\ref{lemma:Zseg} work in this region.

\begin{figure}[h!]
	\centering
	\begin{tabular}{cccc}
		\frame{\includegraphics[width=3.39cm]{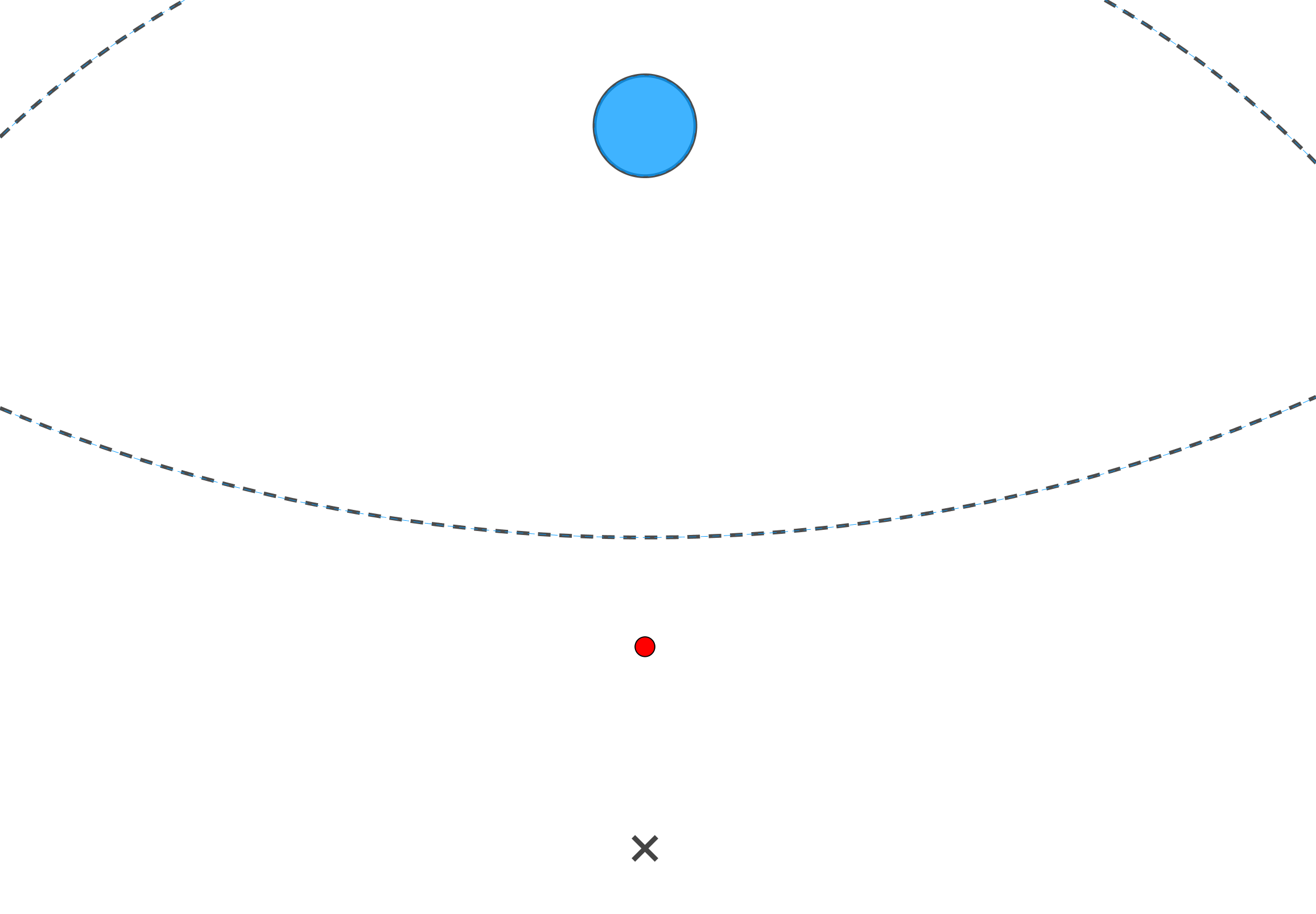}}& & & \\
		$\velocity=-\tfrac{1}{\At}\im+0.3\im$ & & & \\[0.2cm]
		\frame{\includegraphics[width=3.39cm]{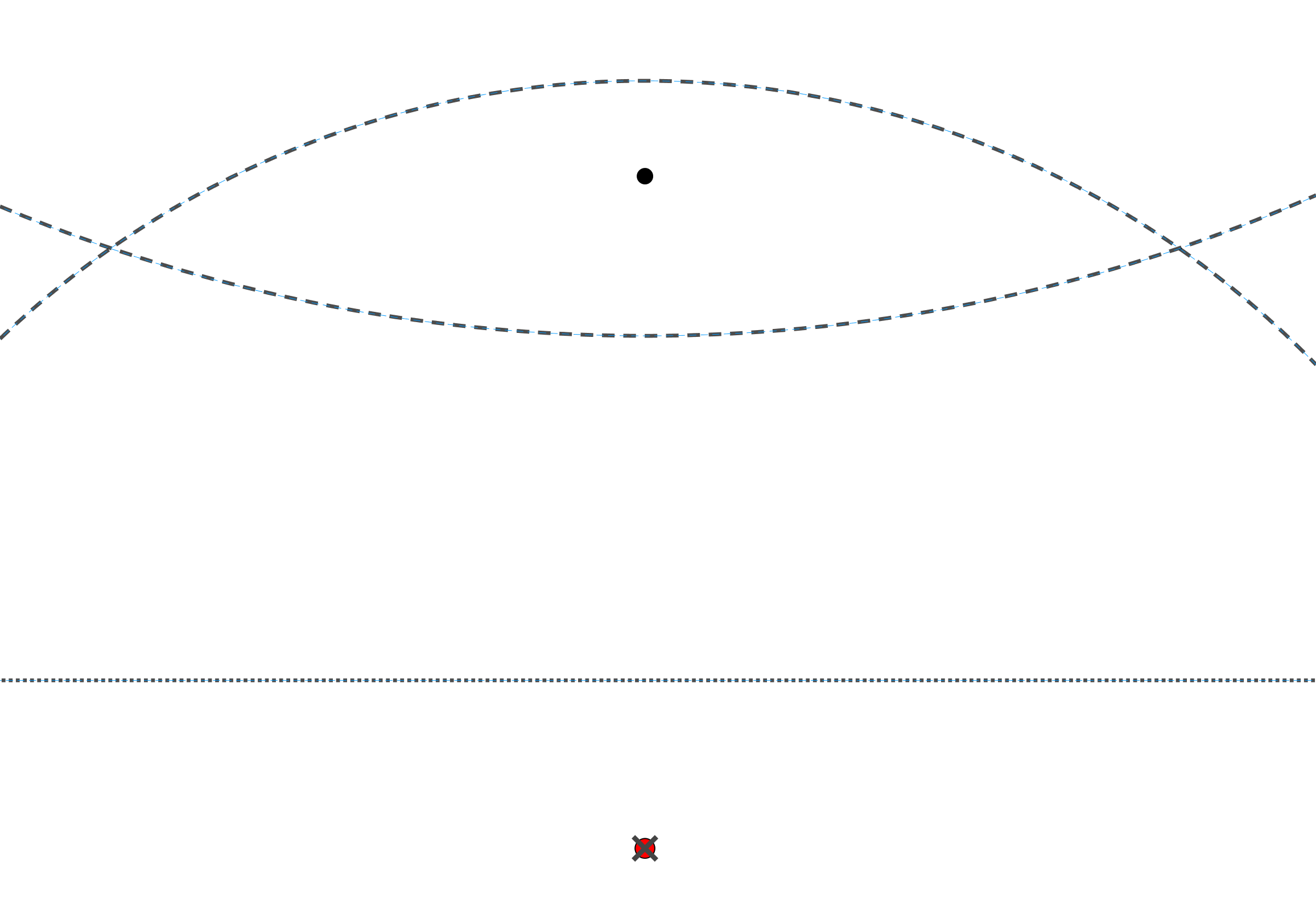}}&
		\frame{\includegraphics[width=3.39cm]{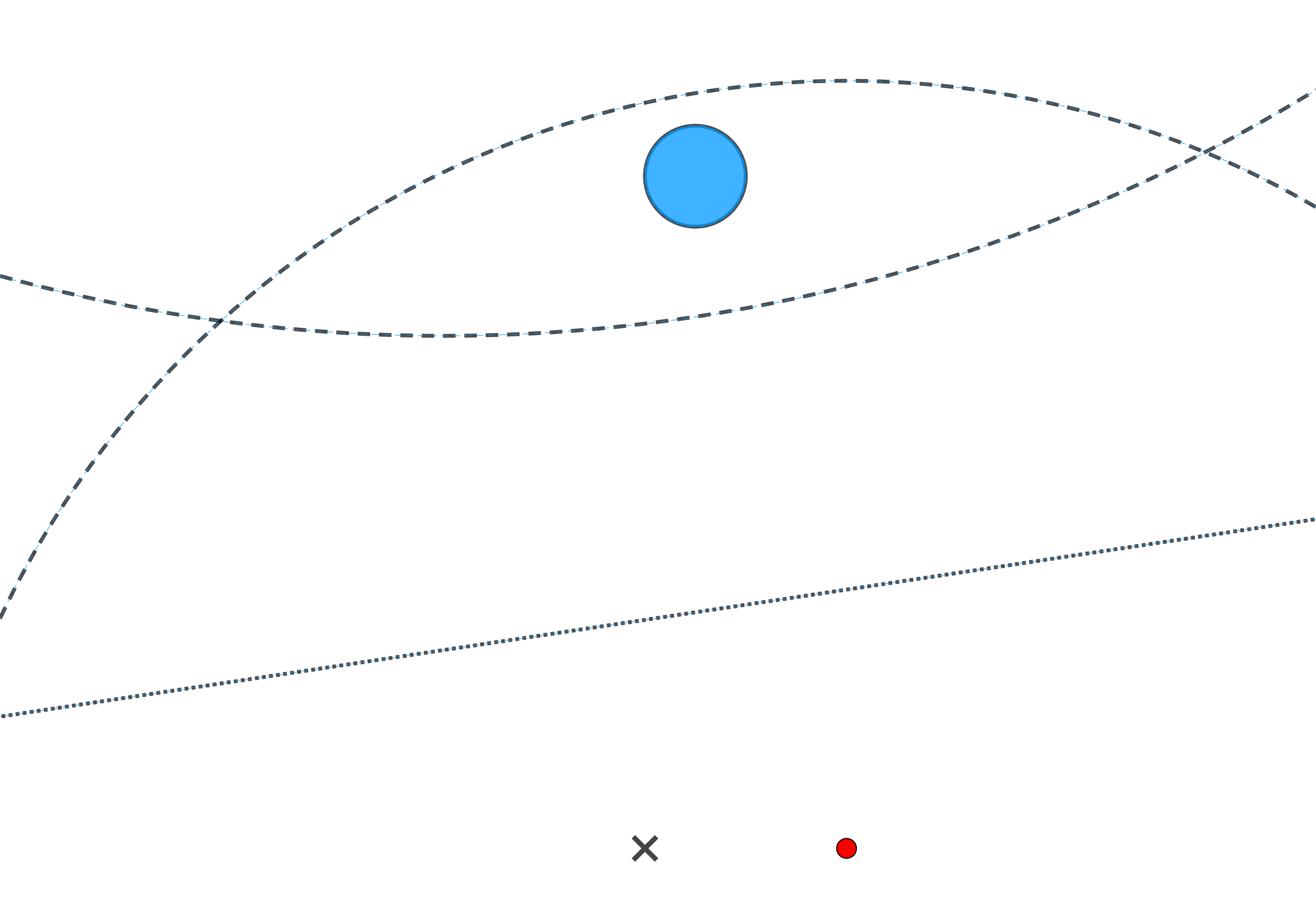}}&
		\frame{\includegraphics[width=3.39cm]{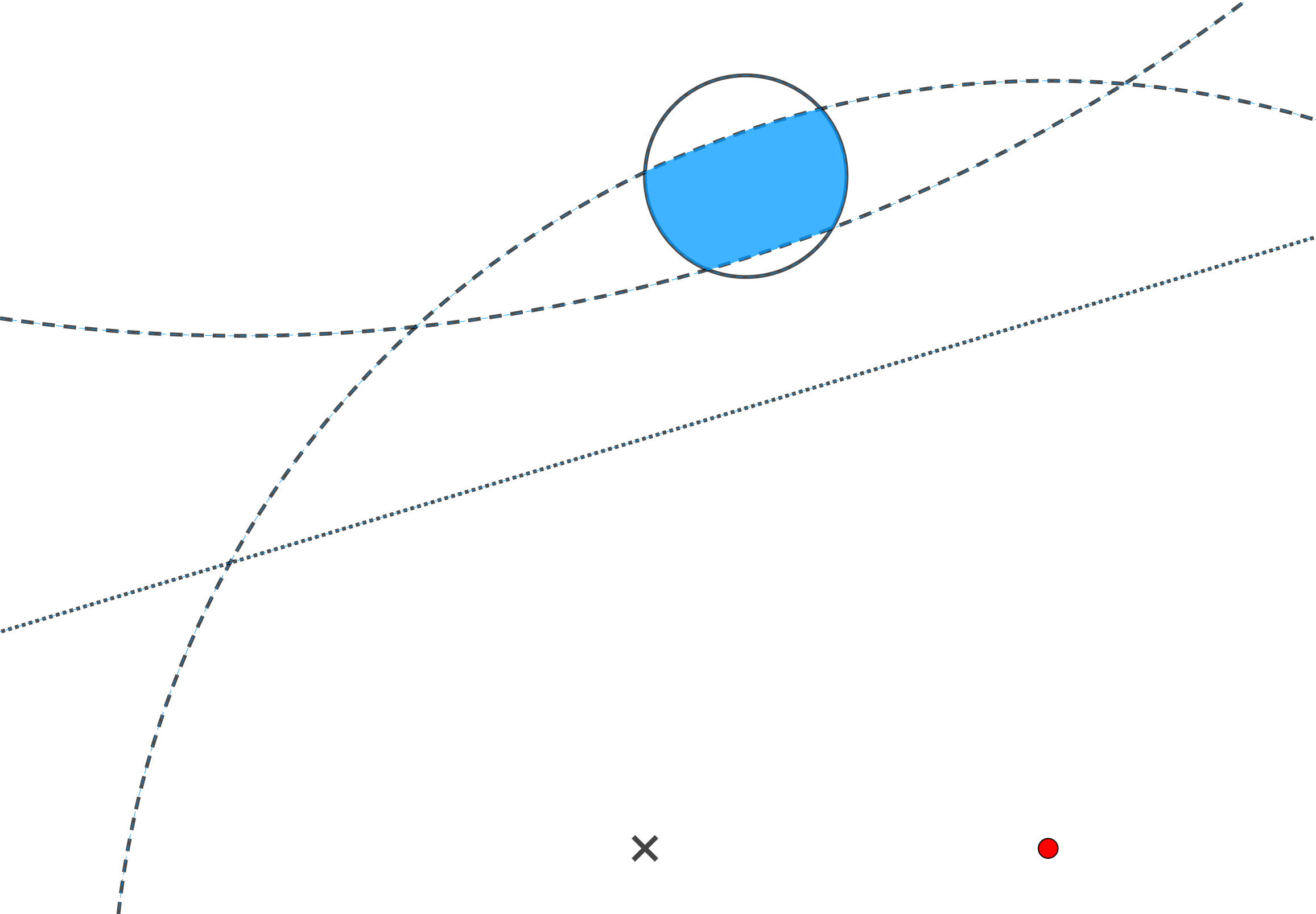}}&
		\frame{\includegraphics[width=3.39cm]{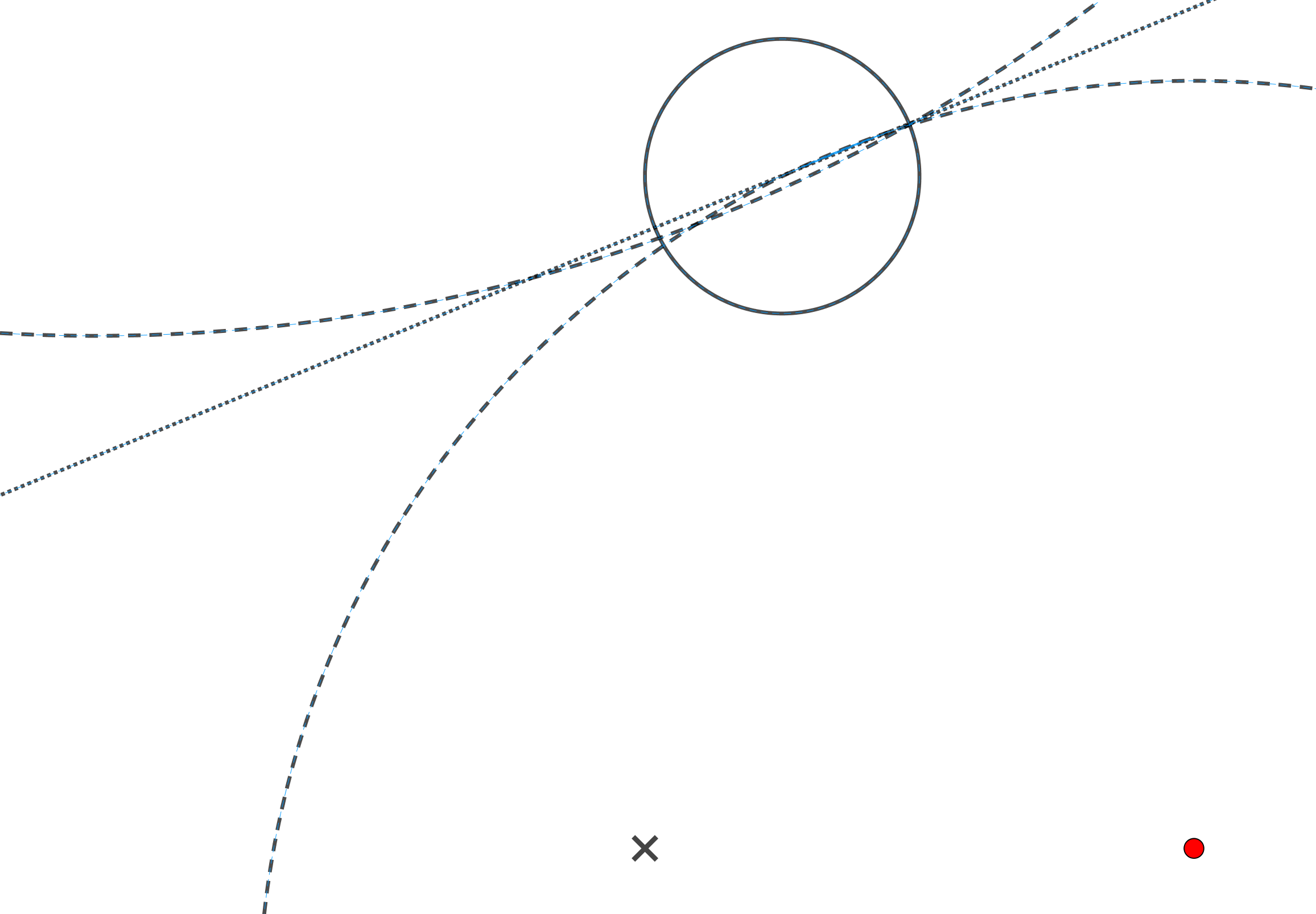}}\\
		$\velocity=-\tfrac{1}{\At}\im$ & $\velocity=-\tfrac{1}{\At}\im+0.3$ & $\velocity=-\tfrac{1}{\At}\im+0.6$ & $\velocity=-\tfrac{1}{\At}\im+0.815$ \\[0.2cm]
		\frame{\includegraphics[width=3.39cm]{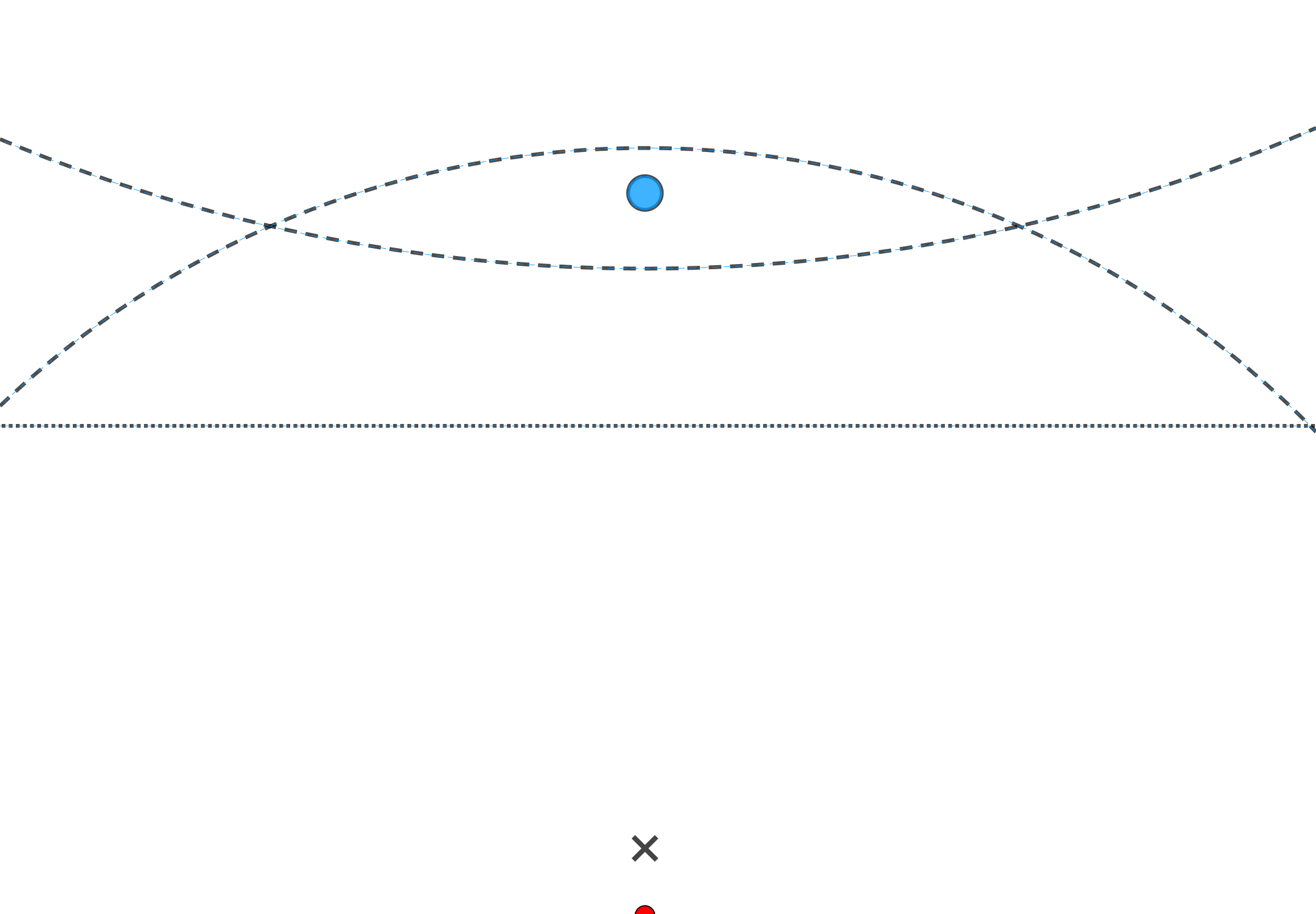}}\\
		$\velocity=-\tfrac{1}{\At}\im-0.1\im$ & & & 
	\end{tabular}
	\caption{Plots of $\rel_{\At,M}(\theta,\velocity)$ (cf.~Fig.~\ref{fig:geo}-right) for $\At=\tfrac{1}{2}$, $M=4>M_*(A)$, $\theta=\tfrac{1}{2}$ and different $\velocity$'s (red point) near the pinch singularity $\At\velocity+\im=0$ (cross) and far from it where $\rel_{\At,M}(\theta,\velocity)$ collapses.}
	\label{geo:MM}
\end{figure}

\indent 2.2) Let $M<M_*(A)$. Then $\rel_{\At,M}(\theta,\velocity)=\emptyset$ in a neighbourhood of $\velocity=-\tfrac{1}{\At}\im$ (cf.~Fig.~\ref{geo:M}). Therefore, there is $\gamma(\At,M)>0$ so that $S_{\gamma}\cap\rel_{\At,M}=\emptyset$.

\begin{figure}[h!]
	\centering
	\begin{tabular}{cccc}
		\frame{\includegraphics[width=3.39cm]{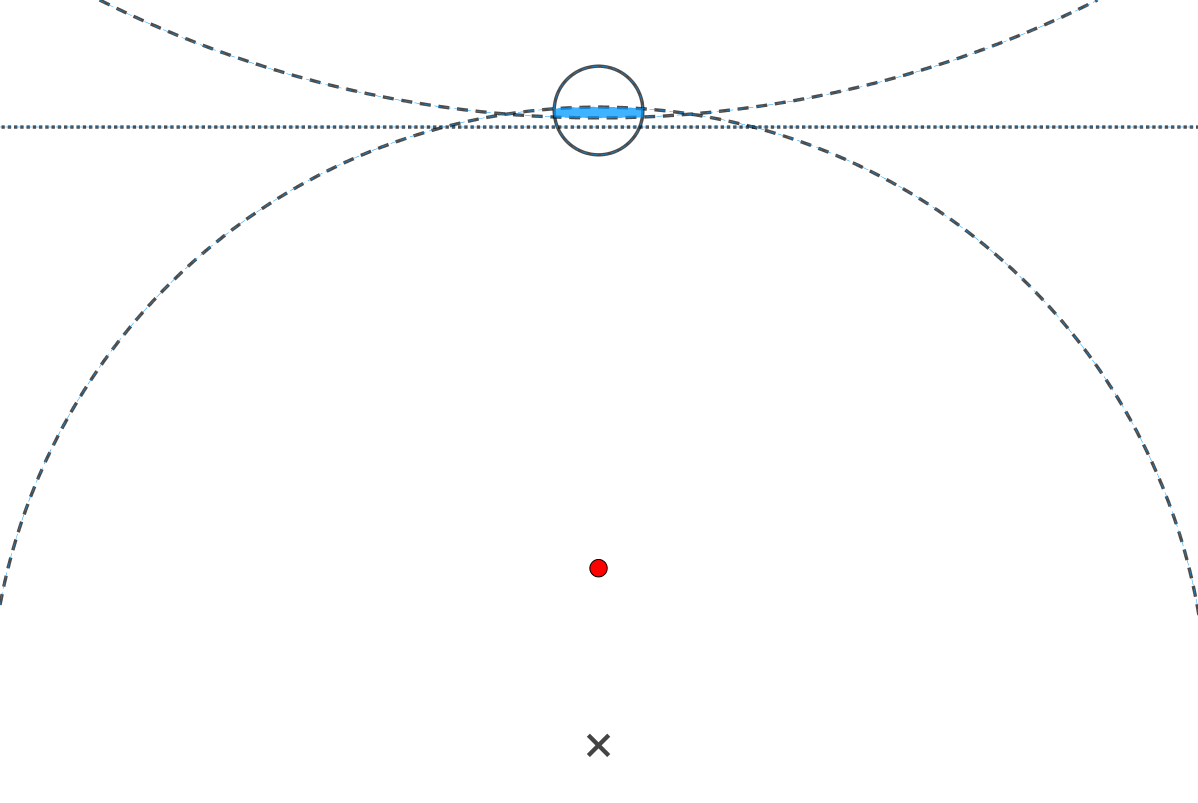}}& & & \\
		$\velocity=-\tfrac{1}{\At}\im+0.3\im$ & & & \\[0.2cm]
		\frame{\includegraphics[width=3.39cm]{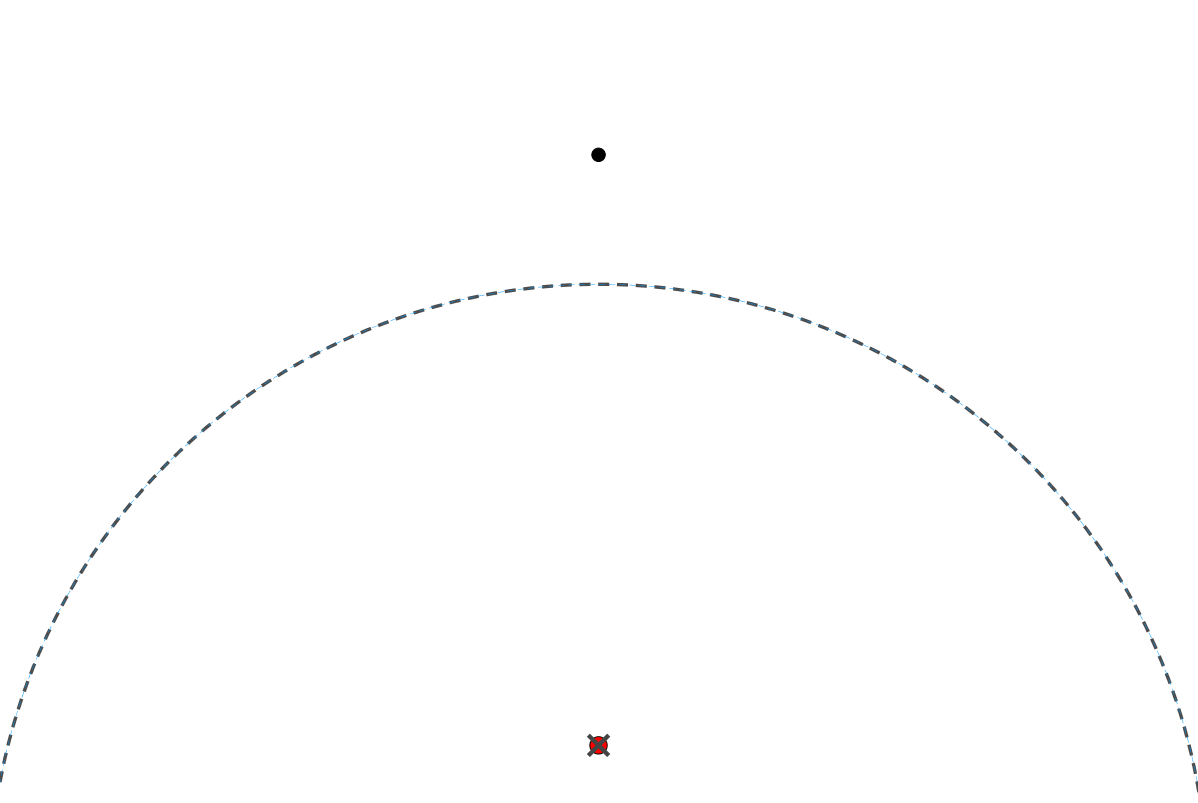}}&
		\frame{\includegraphics[width=3.39cm]{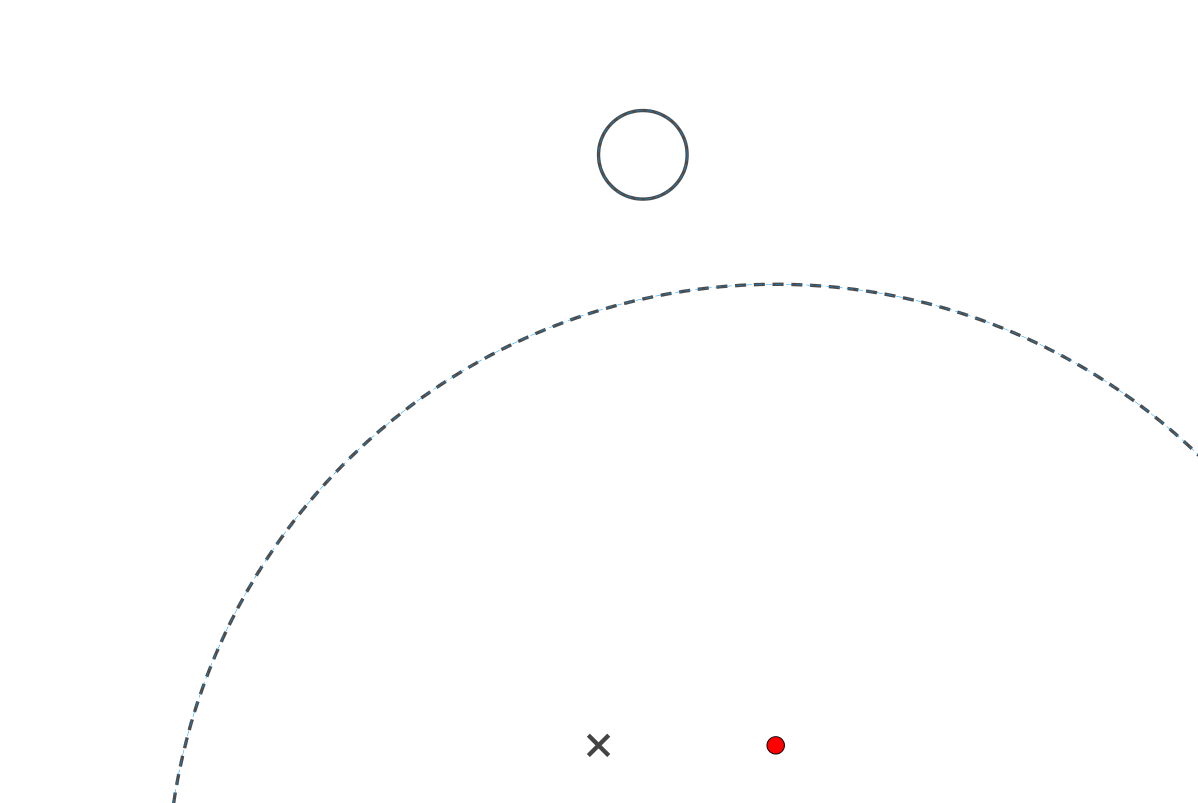}}&
		\frame{\includegraphics[width=3.39cm]{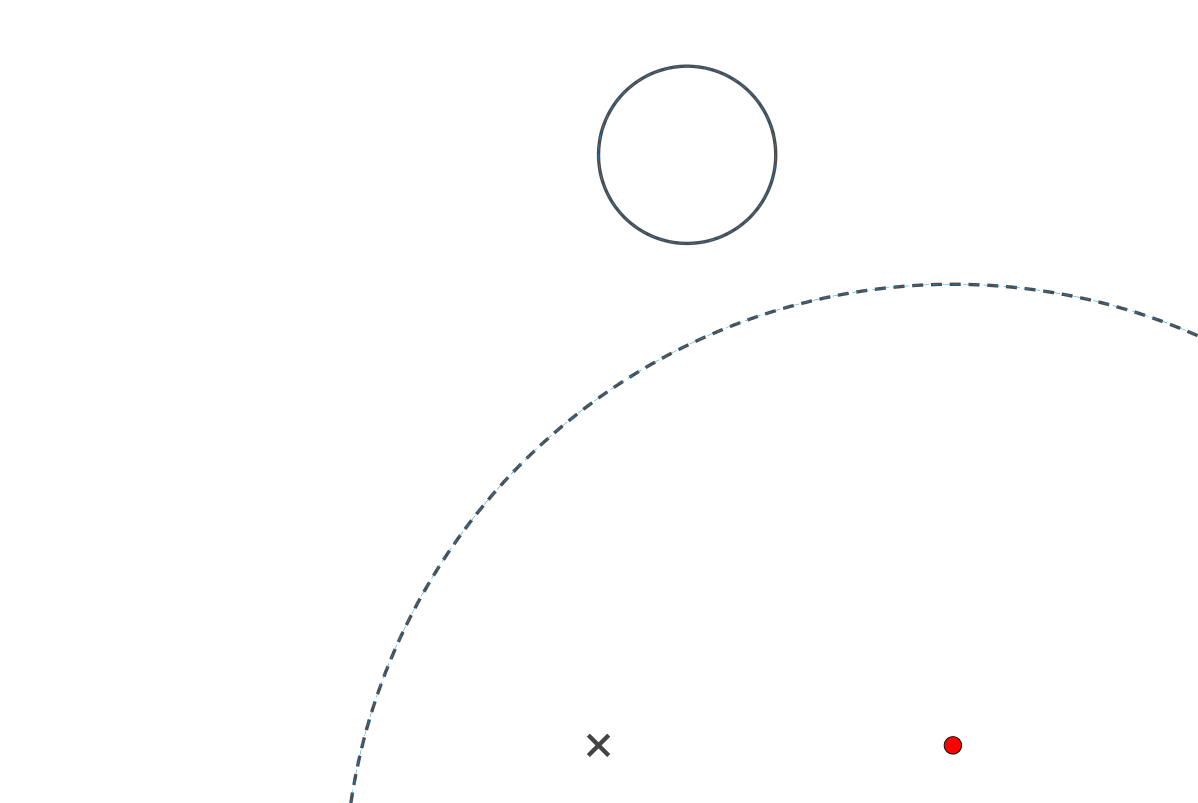}}&
		\frame{\includegraphics[width=3.39cm]{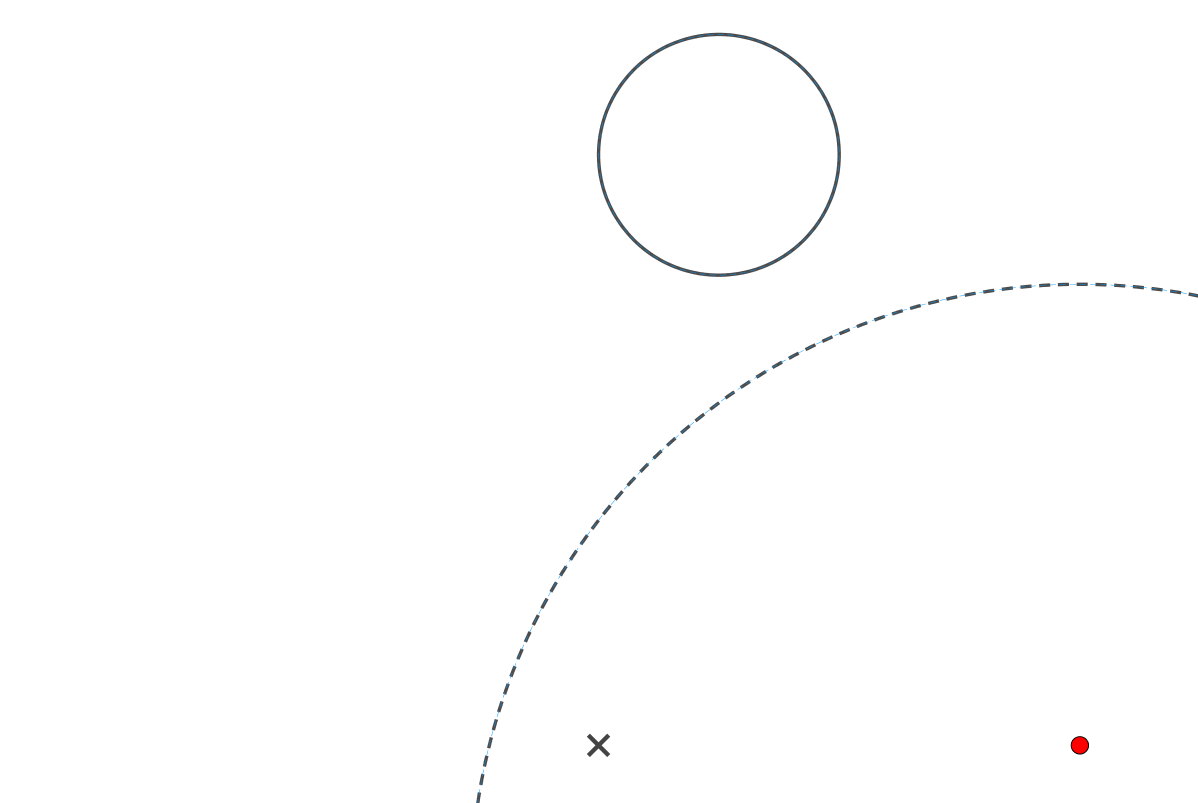}}\\
		$\velocity=-\tfrac{1}{\At}\im$ & $\velocity=-\tfrac{1}{\At}\im+0.3$ & $\velocity=-\tfrac{1}{\At}\im+0.6$ & $\velocity=-\tfrac{1}{\At}\im+0.815$ \\[0.2cm]
		\frame{\includegraphics[width=3.39cm]{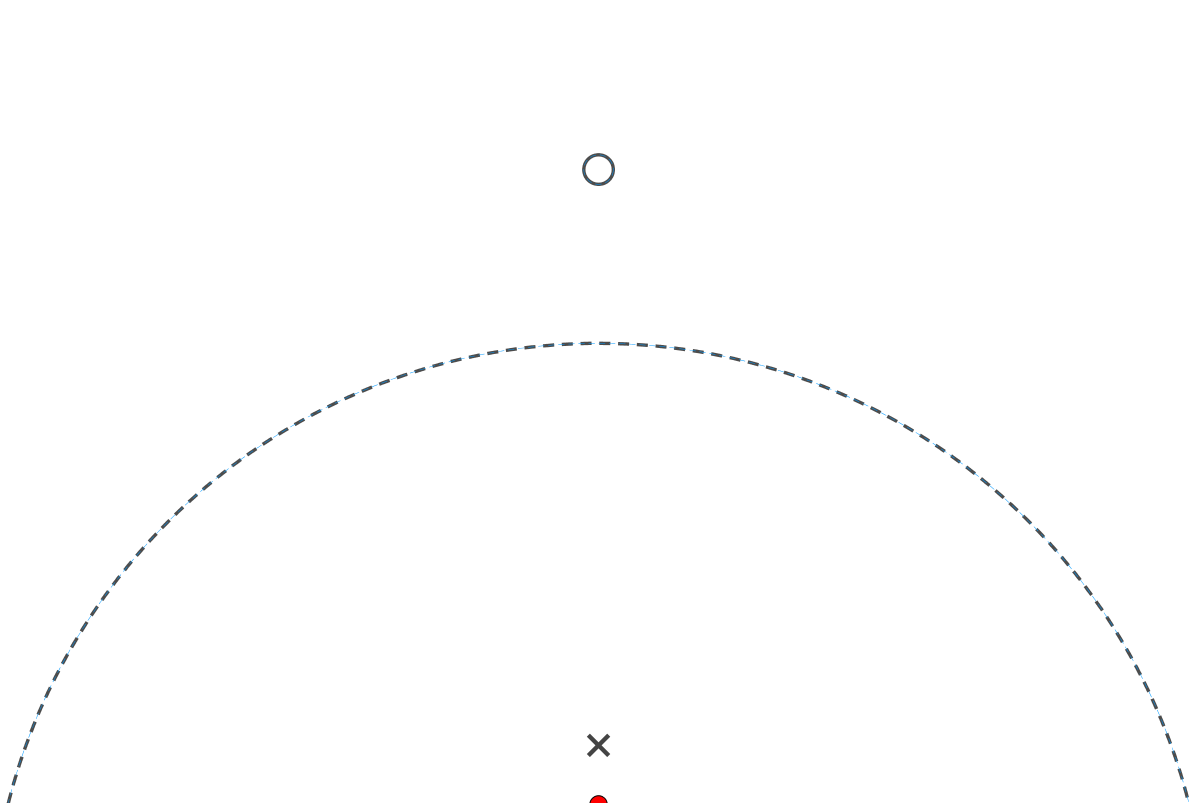}}\\
		$\velocity=-\tfrac{1}{\At}\im-0.1\im$ & & & 
	\end{tabular}
	\caption{Plots of $\rel_{\At,M}(\theta,\velocity)$ (cf.~Fig.~\ref{geo:MM}) for $M=3<M_*(\At)$.}
	\label{geo:M}
\end{figure}

\indent By 2.1) and 2.2), from now on we may assume that $|\At\velocity+\im|>\gamma$ for some fixed $\gamma(\At,M)>0$. We remark in passing that, although we have removed the pinch singularity, it is not clear if  $\partial\rel_{\At,M}\setminus(\Kconstrain\cup S_{\gamma})$ is locally the graph of a Lipschitz function (due to the collapse when $|\velocity|$ grows) thus preventing from following the argument in \cite{RIPM}.\\

\indent\textit{Case $|\At\velocity+\im|>\gamma$}: From now on we focus on states $z=(\theta,\velocity,\aux)\in\rel_M$ with $|\At\velocity+\im|>\gamma$.
In such case, there are $\omega\in\Sc$ and $\sigma_{-},\sigma_{+}\in\D$ so that $\aux$ can be written as
\begin{equation}\label{ZM:1}
\begin{split}
\aux&=\theta\velocity+(1-\theta^2)\frac{\At\velocity+\im}{1+\omega\theta\At}\omega\\
&=\mp\velocity+\frac{1}{2}(1\pm\theta)(M_{\pm}\sigma_{\pm}-\im).
\end{split}
\end{equation}
Thus, $\omega$, $\sigma_{-}$, $\sigma_{+}$ are related via
\begin{equation}\label{ZM:1:3}
\pm\velocity+(1\mp\theta)\frac{\At\velocity+\im}{1+\omega\theta\At}\omega
=\frac{1}{2}(M_{\pm}\sigma_{\pm}-\im).
\end{equation}
By \eqref{ZM:1}, we deduce that the identity \eqref{id1} is equivalent to
\begin{subequations}
\label{id4}
\begin{align}
&(1-\theta^2)\left|\frac{\At\velocity+\im}{1+\omega\theta\At}\right|^2(1-|\shift{\omega}|^2)
+(M^2-1-B(z))\label{id4:1}\\
&=\frac{1-\theta}{2}(M^2-\At)(1-|\sigma_{-}|^2)
+\frac{1+\theta}{2}(M^2+\At)(1-|\sigma_{+}|^2).\label{id4:2}
\end{align}
\end{subequations}
In fact, \eqref{id4} holds for all $z=(\theta,\velocity,\aux)\in\bar{\rel}\setminus\Kconstrain$, with $\omega\in\bar{\Sc}$, $\sigma_-,\sigma_+\in\R^2$ defined via \eqref{ZM:1}.\\

\indent Since $\rel_M$ is open, for every $z\in\rel_M$ and $\bar{z}\in\Lambda$ there is $\epsilon(z,\bar{z},\rel_M)>0$ so that $z_{\lambda}\equiv z+\lambda\bar{z}\in\rel_M$ for all $|\lambda|\leq\epsilon$. However, as in Lemma~\ref{lemma:U}, we must choose $\bar{z}$ carefully in such a way that $\epsilon(1-\theta^2,\At,M)$.
Let us denote $\omega_{\lambda}\in\Sc$ and $\sigma_{\pm,\lambda}\in\D$ by the corresponding points that determine $\aux_{\lambda}$ in the balls $\B(\theta_{\lambda},\velocity_{\lambda})$ and $\B_{\pm}(\theta_{\lambda},\velocity_{\lambda})$ respectively via \eqref{ZM:1}.\\

\indent\textit{Step 1. A change of variables}:
Let $\bar{z}(z)=(1,\bar{\velocity},\bar{\aux})$ be the $\Lambda$-direction we want to construct. Thus, $\bar{\velocity}=\bar{\omega}(\At\bar{\aux}+\im)$ with $(\bar{\aux},\bar{\omega})\in\R^2\times\Ss$ the degrees of freedom. Without loss of generality we take $\bar{\aux}=\LL_{\theta\bar{\omega}}^{-1}(\vect)$ in terms of $\vect\in\R^2$. Inspired by Lemma~\ref{lemma:Zseg}, it is convenient to express w.l.o.g.~this $\vect$ as 
\begin{equation}\label{vect}
\vect(z,\bar{\mathbf{n}},\bar{\omega}):=\velocity+\bar{\mathbf{n}}\frac{\At\velocity+\im}{1+\omega\theta\At}(\bar{\omega}-\omega),
\end{equation}
in terms of some $\bar{\mathbf{n}}\in\R^2$ to be determined. 
Thus, if we denote (recall \eqref{L:1})
\begin{equation}\label{ZM:2}
\mathbf{p}(z,\bar{\mathbf{n}},\bar{\omega}):=\At\bar{\aux}+\im
=\frac{\At\vect+\im}{1+\bar{\omega}\theta\At}
=\frac{\At\velocity+\im}{1+\bar{\omega}\theta\At}
\left(1+\frac{\At\bar{\mathbf{n}}(\bar{\omega}-\omega)}{1+\omega\theta\At}\right),
\end{equation}
the $\Lambda$-direction $\bar{z}$ is written as
\begin{equation}\label{ZM:3}
\bar{\velocity}=\bar{\omega}\mathbf{p},
\quad\quad
\bar{\aux}
=\vect-\theta\bar{\velocity},
\end{equation}
in terms of $(\bar{\mathbf{n}},\bar{\omega})\in\R^2\times\Ss$, which shall be determined in the \textit{step 2} and \textit{3} respectively.\\

\indent\textit{Step 2. Choice of $\bar{\mathbf{n}}$}:
Let us expand the condition $\aux_{\lambda}\in\B_{\pm}(\theta_{\lambda},\velocity_{\lambda})$ in terms of $\lambda$:
\begin{equation}\label{Bpm:1}
\begin{split}
2(\aux_\lambda\pm\velocity_\lambda)+(1\pm\theta_\lambda)\im
&=2(\aux\pm\velocity)+(1\pm\theta)\im
+\lambda(2(\bar{\aux}\pm\bar{\velocity})\pm\im)\\
&=M_\pm(1\pm\theta_\lambda)\sigma_\pm
+\lambda v_{\pm}(z,\bar{z}),
\end{split}
\end{equation}
where we have abbreviated (recall \eqref{ZM:1}-\eqref{ZM:3})
\begin{equation}\label{ZM:5}
\begin{split}
\frac{1}{2}v_{\pm}(z,\bar{z})&:=(\bar{\aux}\pm\bar{\velocity})\mp\frac{1}{2}(M_\pm\sigma_\pm-\im)\\
&\,\,=(\vect-\velocity)\pm(1\mp\theta)\left(\frac{\At\vect+\im}{1+\bar{\omega}\theta\At}\bar{\omega}-\frac{\At\velocity+\im}{1+\omega\theta\At}\omega\right)\\
&\,\,=\frac{1\pm\bar{\omega}\At}{1+\bar{\omega}\theta\At}(\vect-\velocity)
\pm\frac{(1\mp\theta)(\At\velocity+\im)}{(1+\bar{\omega}\theta\At)(1+\omega\theta\At)}(\bar{\omega}-\omega)\\
&\,\,=\frac{(1\pm\bar{\omega}\At)(\At\velocity+\im)}{(1+\bar{\omega}\theta\At)(1+\omega\theta\At)}\left(\bar{\mathbf{n}}\pm\frac{1\mp\theta}{1\pm\bar{\omega}\At}\right)(\bar{\omega}-\omega).
\end{split}
\end{equation}
From \eqref{Bpm:1} we deduce that
\begin{equation}\label{ZM:4}
(1-|\sigma_{\pm,\lambda}|^2)
=(1-|\sigma_{\pm}|^2)
-\tilde{\lambda}_{\pm}v_{\pm}\cdot(2\sigma_{\pm}+\tilde{\lambda}_{\pm}v_{\pm}),
\end{equation}
with
$$\tilde{\lambda}_{\pm}\equiv\frac{\lambda}{M_{\pm}(1\pm\theta_{\lambda})}.$$
Notice that $(1\pm\theta_{\lambda})\geq\tfrac{1}{2}(1\pm\theta)\geq\tfrac{1}{2}(1-|\theta|)$ provided $|\lambda|\leq\tfrac{1}{2}(1-|\theta|)$.\\
The identities \eqref{ZM:5}\eqref{ZM:4} determines a good choice of $\bar{\mathbf{n}}$. More precisely, let us assume w.l.o.g.~that
$|\sigma_{-}|\leq|\sigma_{+}|$ (the case $|\sigma_{+}|<|\sigma_{-}|$ is totally analogous). Then, it is convenient to take (in fact necessary on $(\partial\B_{+}\setminus\partial\B)(\theta,\velocity)$) 
\begin{equation}\label{ZM:8}
\bar{\mathbf{n}}(z,\bar{\omega})=-\frac{1-\theta}{1+\bar{\omega}\At},
\end{equation}
with $\bar{\omega}$ to be determined yet. With this choice of $\bar{\mathbf{n}}$, \eqref{ZM:5} reads as
\begin{equation}\label{ZM:7}
v_{+}(z,\bar{\omega})=0,\quad\quad
v_{-}(z,\bar{\omega})=-\frac{ 4}{1+\bar{\omega}\At}\frac{\At\velocity+\im}{1+\omega\theta\At}(\bar{\omega}-\omega),
\end{equation}
and \eqref{ZM:2} reads as
\begin{equation}\label{q}
\mathbf{p}(z,\bar{\omega})=\frac{1+\omega\At}{1+\bar{\omega}\At}\frac{\At\velocity+\im}{1+\omega\theta\At}
=:\frac{\mathbf{q}(z)}{1+\bar{\omega}\At},
\end{equation}
where we have introduced $\mathbf{q}(z)$ as the part of $\mathbf{p}(z,\bar{\omega})$ independent of $\bar{\omega}$.
Hence, by \eqref{ZM:7}, \eqref{ZM:4} reads as $|\sigma_{+,\lambda}|=|\sigma_{+}|$,
and so $\aux_{\lambda}\in\B_{+}(\theta_{\lambda},\velocity_{\lambda})$ trivially for all $|\lambda|<(1-|\theta|)$.\\
\indent In summary, we have seen that we can take $\bar{\mathbf{n}}$ (depending on whether $|\sigma_-|\leq|\sigma_+|$ or $|\sigma_+|<|\sigma_-|$\footnote{If $|\sigma_+|<|\sigma_-|$ we take $\bar{\mathbf{n}}(z,\bar{\omega})=\frac{1+\theta}{1-\bar{\omega}\At}$ and so \eqref{ZM:5} reads as $v_+(z,\bar{\omega})=\frac{ 4}{1-\bar{\omega}\At}\frac{\At\velocity+\im}{1+\omega\theta\At}(\bar{\omega}-\omega)$, $v_-(z,\bar{\omega})=0$ and \eqref{ZM:2} reads as
$\mathbf{p}(z,\bar{\omega})=\frac{1-\omega\At}{1-\bar{\omega}\At}\frac{\At\velocity+\im}{1+\omega\theta\At}
=:\frac{\mathbf{q}(z)}{1-\bar{\omega}\At}$ for a slightly different $\mathbf{q}$.}) in such a way that the condition $\aux_{\lambda}\in\B_+(\theta_{\lambda},\velocity_{\lambda})$ (or $\B_-(\theta_{\lambda},\velocity_{\lambda})$) holds for all $|\lambda|<(1-|\theta|)$.
Thus, it remains to control the other three inequalities in \eqref{id2}, i.e.~$\B_{-}$, $\B$ and $\HP$.\\

\indent\textit{Step 3. Choice of $\bar{\omega}$}:
By \eqref{ZM:4}\eqref{ZM:7}, the condition $\aux_{\lambda}\in\B_{-}(\theta_{\lambda},\velocity_{\lambda})$ can be written as
\begin{equation}\label{ZM:9}
\tilde{\lambda}_{-}\mathcal{O}(|\shift\bar{\omega}-\shift\omega|)
<(1-|\sigma_{-}|^2).
\end{equation}
Notice that, since  $|\At\velocity+\im|>\gamma$ and $|\theta_{-}|\leq|\theta_{+}|$, the identity \eqref{id4} yields
\begin{equation}\label{ZM:10}
\begin{split}
\frac{1}{4}(1-\theta^2)\gamma^2(1-|\shift\omega|^2)
&\leq(1-\theta^2)\left|\frac{\At\velocity+\im}{1+\omega\theta\At}\right|^2(1-|\shift{\omega}|^2)\\
&\leq\eqref{id4:1}=\eqref{id4:2}\\
&\leq (M^2+|\At|)(1-|\sigma_{-}|^2).
\end{split}
\end{equation}
\indent Since $\mathbf{v}=\velocity+\mathcal{O}(|\shift\bar{\omega}-\shift\omega|)$ \eqref{vect}, by elementary computations as in the \textit{step 2} of the proof of Lemma~\ref{lemma:Zseg}, we deduce that the condition $\aux_{\lambda}\in\B(\theta_{\lambda},\velocity_{\lambda})$ can be written as
\begin{equation}\label{ZM:11}
\lambda\mathcal{O}(|\shift\bar{\omega}-\shift\omega|)
<(1-\theta^2)(1-|\shift\omega|^2).
\end{equation}
\indent In summary, by \eqref{ZM:10}, to guarantee that \eqref{ZM:9}\eqref{ZM:11} hold (for all $|\lambda|$ depending on $(1-\theta^2)$) it is enough to show that we can take $\bar{\omega}\in\Ss$ satisfying $|\shift\bar{\omega}-\shift\omega|\lesssim(1-|\shift\omega|)$ as $|\shift\omega|\uparrow 1$. This suggests to take $\shift\bar{\omega}$ by the projection $\frac{\shift\omega}{|\shift\omega|}$ as in Lemma~\ref{lemma:Zseg}. However, the last inequality \eqref{id2:2} restricts the set of admissible $\bar{\omega}$'s. Let us see it.\\

\indent Let us expand the condition $\aux_{\lambda}\in\HP(\theta_{\lambda},\velocity_{\lambda})$ in terms of $\lambda$:
\begin{equation}\label{ZM:6}
(M^2-1-B(z_{\lambda}))
=(M^2-1-B(z))-\lambda b(z,\bar{\omega}),
\end{equation}
where $b\equiv b_{\At}$ is
\begin{equation}\label{b}
\begin{split}
b(z,\bar{\omega})
&:=4\bar{\velocity}\cdot(\velocity+\At\aux+\theta\im+\At\im)
+4\velocity\cdot(\bar{\velocity}+\At\bar{\aux}+\im)\\
&\,\,=4(\bar{\omega}\mathbf{p})\cdot\mathbf{b}+4\velocity\cdot((\bar{\omega}+1)\mathbf{p})\\
&\,\,=2\mathbf{p}\cdot(\shift{\bar{\omega}}^*(\mathbf{b}+\velocity)-(\mathbf{b}-\velocity)).
\end{split}
\end{equation}

\indent Before continuing with the choice of $\bar{\omega}$, let us remark a difference to the case of equal viscosities. For $\At=0$, the functions $B_0$, $\mathbf{b}_0$ and $b_0$ do not depend on $\aux$ (equiv. $\omega$). As a result,
given $(\theta,\velocity)\in(-1,1)\times\R^2$, the set of  $\bar{\omega}$'s that can be used as $B_0(\theta,\velocity)\uparrow M^2-1$ (i.e.~$\velocity$ tends to $\partial\mathscr{B}(\theta)$) is more explicit, namely this is $\Omega_0(\theta,\velocity)
=\{\bar{\omega}\in\Ss\,:\,\aux_{\bar{\omega}}\equiv\theta\velocity+(1-\theta^2)\bar{\omega}\im\in(\bar{\B}_-\cap\bar{\B}_+)(\theta,\velocity)\}$
(i.e.~$\aux_{\bar{\omega}}\in(\partial\B\cap\bar{\B}_-\cap\bar{\B}_+)(\theta,\velocity)$), independently of $\aux$. Thus, for each $\aux\in\rel_{0,M}(\theta,\velocity)$, the choice of $\bar{\omega}$ in \cite{RIPM} is the minimizer of $|\bar{\omega}-\omega|$ in $\Omega_0(\theta,\velocity)$. To conclude, Sz\'ekelyhidi checked that the circles $\partial\B_{\pm}(\theta,\velocity)$ intersect $\partial\B(\theta,\velocity)$ transversally.
For $\At\neq 0$, the analogous set of $\bar{\omega}$'s depends on $(\theta,\velocity,\aux)$, in terms of the proximity to the boundary of the half-plane $\HP(\theta,\velocity)$, and it is less explicit. In this regard, for $\At\neq 0$, instead of figuring out how is $\Omega_{\At}(\theta,\velocity,\aux)$, we design a suitable $\bar{\omega}$ for each $z$ separately.\\

\indent As in \cite{RIPM}, in order to choose $\bar{\omega}$ we distinguish three cases (see Fig.~\ref{fig:delta}) depending on some parameter $0<\delta(1-\theta^2,\At,M,\gamma)<M^2-1$ which shall be determined in the \textit{step 4}.

\begin{figure}[h!]
	\centering
	\includegraphics[height=5.6cm]{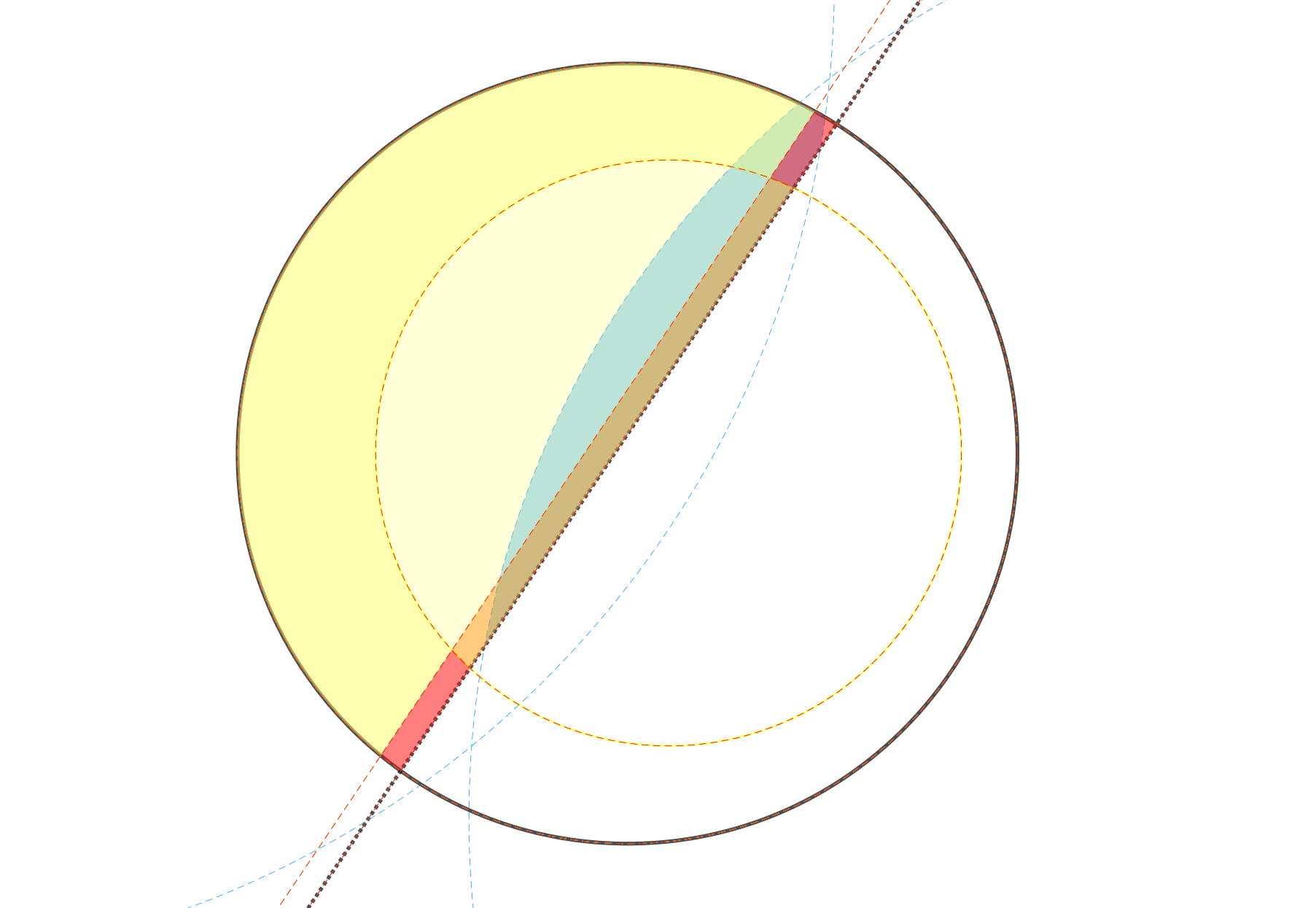}
	\caption{Plot of the various regions dividing $\rel_{\At,M}(\theta,\velocity)$ in terms of some $\delta>0$ small, for some $0<|\At|<1$, $M>1$, $(\theta,\velocity)\in(-1,1)\times\R^2$. Over $\rel_{\At,M}(\theta,\velocity)$ (cf.~Fig.~\ref{fig:geo}-right) we have overlapped: the circle\protect\footnotemark $\,\,(1-|\shift\omega|)=\delta$, the line $M^2-1-B(z)=\delta$, and the regions: 1) $M^2-1-B(z)>\delta$ (yellow: lighter if $(1-|\shift\omega|)>\delta$, darker if $(1-|\shift\omega|)\leq\delta$), 2) $M^2-1-B(z)\leq\delta$ coupled with either 2.1) $(1-|\shift\omega|)>\delta$ (orange) or 2.2) $(1-|\shift\omega|)\leq\delta$ (red).}
	\label{fig:delta}
\end{figure}
\footnotetext{Recall that $\varphi_{\theta\At}\in\mathrm{Aut}(\Sc)$ ($\subset$ M\"obius transformations) and so it preserves circles.}

\indent 1) If $M^2-1-B(z)>\delta$ (cf.~Fig.~\ref{fig:delta}-yellow) we can take directly $\bar{\omega}\in\Ss$ as in Lemma~\ref{lemma:Zseg}, that is $\bar{\omega}=0$ if $|\shift\omega|\leq\tfrac{1}{2}$ and $\shift\bar{\omega}=\frac{\shift\omega}{|\shift\omega|}$ if $\tfrac{1}{2}<|\shift\omega|<1$ (clearly $|\shift\bar{\omega}-\shift\omega|\lesssim(1-|\shift\omega|)$). Notice that there is $B(\At,M)>0$ so that $|b(z,\bar{\omega})|\leq B$. Hence, by \eqref{ZM:6}, $\aux_{\lambda}\in\HP(\theta_{\lambda},\velocity_{\lambda})$ for all $|\lambda|<\delta/B$.\\ 
\indent 2) Now let us suppose that $M^2-1-B(z)\leq\delta$.\\ 
\indent 2.1) In this case, if $(1-|\shift\omega|)>\delta$ (cf.~Fig.~\ref{fig:delta}-orange), then \eqref{ZM:9}\eqref{ZM:11} hold for all $|\lambda|\lesssim(1-\theta^2)^2\delta$. Thus, as we shall see in \textit{step 4}, there exists $\bar{\omega}$ satisfying $b(z,\bar{\omega})=0$.  With such choice, \eqref{ZM:6} reads as $B(z_{\lambda})=B(z)$, and so $\aux_{\lambda}\in\HP(\theta_{\lambda},\velocity_{\lambda})$ trivially for all $|\lambda|<(1-|\theta|)$. \\ 
\indent 2.2) Finally let us suppose that $(1-|\shift\omega|)\leq\delta$ (cf.~Fig.~\ref{fig:delta}-red). 
As we have seen, on the one hand, if $\aux\in\partial\HP(\theta,\velocity)$ we have to take $\bar{\omega}$ satisfying $b(z,\bar{\omega})=0=:\alpha_{\HP}(z)$. On the other hand, if $\aux\in\partial\B(\theta,\velocity)$ we have to take $\bar{\omega}=\omega$. Furthermore, for any $\aux\in\partial\B(\theta,\velocity)$ (not necessarily on $\bar{\rel}_M(\theta,\velocity)$) by applying $\vect(z,\omega)=\velocity$, $v_{\pm}(z,\omega)=0$,  Lemma~\ref{lemma:K1}\ref{K1:1}, \eqref{ZM:4} and \eqref{ZM:6}, the coefficient of order 1 in $\lambda$ of the identity \eqref{id4} reads as
$$b(z,\omega)
=\frac{1}{2}((M^2-\At)(1-|\sigma_{-}|^2)
-(M^2+\At)(1-|\sigma_{+}|^2))=:\alpha_{\B}(z).$$
Hence, both cases are compatible because, if $\aux\in(\partial\B\cap\partial\HP)(\theta,\velocity)$, the identity \eqref{id4} implies that  $\aux\in(\partial\B_-\cap\partial\B_+)(\theta,\velocity)$ too (cf.~Fig.~\ref{fig:geo}) and so $\alpha_{\B}(z)=0=\alpha_{\HP}(z)$.\\
For states near the boundary, what we would like is to find $\bar{\omega}\in\Ss$ satisfying
\begin{equation}\label{ZM:beq}
b(z,\bar{\omega})=\alpha(z),
\end{equation}
for some suitable interpolation $\alpha(z)$ from the values that $b$ must take on the walls $\partial\HP(\theta,\velocity)$ and $\partial\B(\theta,\velocity)$. In this regard, here we consider a convex combination of $\alpha_{\B}$ and $\alpha_{\HP}$
\begin{equation}\label{alpha}
\begin{split}
\alpha(z)
:=&\,\,\frac{(M^2-1-B(z))+d(z)}{\eqref{id4:1}+d(z)}\alpha_{\B}(z)
+\frac{\displaystyle(1-\theta^2)\left|\frac{\At\velocity+\im}{1+\omega\theta\At}\right|^2(1-|\shift{\omega}|^2)+d(z)}{\eqref{id4:1}+d(z)}\alpha_{\HP}(z)\\
=&\,\,\frac{(M^2-1-B(z))+d(z)}{\eqref{id4:1}+d(z)}
\frac{1}{2}((M^2-\At)(1-|\sigma_{-}|^2)
-(M^2+\At)(1-|\sigma_{+}|^2)),
\end{split}
\end{equation}
where we have introduced $d(z):=8(1\vee|\At\velocity|)\mathrm{dist}(\aux;\rel_M(\theta,\velocity))$ to extend $\alpha$ on $\bar{\B}(\theta,\velocity)\setminus\rel_M(\theta,\velocity)$  (notice that $d(z)\geq 2|M^2-1-B(z)|$ on $\bar{\B}(\theta,\velocity)\setminus\rel_M(\theta,\velocity)$).
For instance, if $\aux\in\partial\B_{\pm}(\theta,\velocity)$ we have
$$\pm\alpha(z)
=\frac{M^2-1-B(z)}{1\mp\theta}
=\frac{1}{2}(M^2\mp\At)(1-|\sigma_{\mp}|^2)-(1\pm\theta)\left|\frac{\At\velocity+\im}{1+\omega\theta\At}\right|^2(1-|\shift{\omega}|^2).$$
Hence, if there is such $\bar{\omega}\in\Ss$ satisfying \eqref{ZM:beq} for \eqref{alpha}, then \eqref{ZM:6} reads as
$$M^2-1-B(z_{\lambda})
=\frac{M^2-1-B(z)}{\eqref{id4:2}}
\left(\frac{1-\theta_{\lambda}}{2}(M^2-\At)(1-|\sigma_{-}|^2)
+\frac{1+\theta_{\lambda}}{2}(M^2+\At)(1-|\sigma_{+}|^2)\right),$$
and so $\aux_{\lambda}\in\HP(\theta_{\lambda},\velocity_{\lambda})$ for all $|\lambda|<(1-|\theta|)$. 
Thus, it remains to show that there is $\bar{\omega}\in\Ss$ satisfying \eqref{ZM:beq} and that the corresponding map $\omega\mapsto\bar{\omega}$ is Lipschitz (see \eqref{ZM:Lip1}\eqref{ZM:Lip2}). \\

\indent\textit{Step 4. Lipschitz solution to $b(z,\bar{\omega})=\alpha$}: 
Firstly, let us determine the solvability of $b(z,\bar{\omega})=\alpha$ for states  $\aux\in\bar{\B}(\theta,\velocity)$ and $\alpha\in\R$. By \eqref{q}\eqref{b}, there is such $\bar{\omega}\in\R^2$ if and only if 
$$\frac{\shift\bar{\omega}^*(\mathbf{b}+\velocity)-(\mathbf{b}-\velocity)}{1+\At\bar{\omega}^*}
=\frac{1}{2}\frac{\alpha+\beta\im}{\mathbf{q}^*},$$
or equivalently
$$(4\mathbf{q}^*(\mathbf{b}+\velocity)-\At(\alpha+\beta\im))\shift\bar{\omega}^*
=4\mathbf{q}^*(\mathbf{b}-\velocity)
+(2-\At)(\alpha+\beta\im),$$
for some real $\beta$.
Since we require $\bar{\omega}\in\Ss$, necessarily
$$|4\mathbf{q}^*(\mathbf{b}+\velocity)-\At(\alpha+\beta\im)|
=|4\mathbf{q}^*(\mathbf{b}-\velocity)
+(2-\At)(\alpha+\beta\im)|,$$
which turns out to be a quadratic equation for $\beta$, $a_2\beta^2+a_1\beta+a_0=0$, where
\begin{align*}
a_2&=(1-\At)>0,\\
a_1&=4((1-\At)\mathbf{b}-\velocity)\cdot\mathbf{q}^\perp,\\
a_0&=(1-\At)\alpha^2
+4((1-\At)\mathbf{b}-\velocity)\cdot\mathbf{q}^\perp\alpha
-4B(z)|\mathbf{q}|^2.
\end{align*}
The discriminant of this quadratic equation verifies
$$
\Delta(z,\alpha)
=a_1^2-4a_2a_0
\geq 16(1-\At)B(z)|\mathbf{q}(z)|^2+\mathcal{O}(\alpha).
$$
In particular, if $B(z)\geq M^2-1-\delta>0$, for $\alpha=0$ we have $\Delta(z,0)>0$ and so there exists $\bar{\omega}\in\Ss$ satisfying $b(z,\bar{\omega})=0$.
Now let $\alpha(z)$ given in \eqref{alpha}. Notice that this can be bounded by
$$|\alpha(z)|
\leq\frac{1}{2}(M^2+|\At|)(|1-|\sigma_{-}|^2|+|1-|\sigma_{+}|^2|).$$
Hence, since $|\mathbf{q}(z)|\geq\frac{1-|\At|}{1+|\At|}\gamma$, there is a constant $C(A,M,\gamma)>0$ so that 
$$\Delta(z,\alpha(z))\geq 4(1-\At)(M^2-1)\left(\frac{1-|\At|}{1+|\At|}\gamma\right)^2>0,$$
for all $\aux\in\bar{\B}(\theta,\velocity)$ in the intersection of the half-plane $B(z)\geq\tfrac{1}{2}(M^2-1)$ and the annuli $|1-|\sigma_{-}|^2|,|1-|\sigma_{+}|^2|\leq C$. Therefore, in this region $L\equiv L_{\At,M,\gamma}$ 
$$L(\theta,\velocity):=\{\aux\in\bar{\B}(\theta,\velocity)\,:\,B(z)\geq\tfrac{1}{2}(M^2-1),\,|1-|\sigma_{-}|^2|,|1-|\sigma_{+}|^2|\leq C\}$$ 
there are two ($s\in\{-1,1\}$) solutions $\shift\bar{\omega}_s=q_s(z)$ to $b(z,\bar{\omega})=\alpha(z)$ given by
\begin{equation}\label{balphasol}
q_s(z)
:=\frac{4\mathbf{q}(z)(\mathbf{b}(z)-\velocity)^*
	+(2-\At)(\alpha-\beta_{s}\im)(z)}
{4\mathbf{q}(z)(\mathbf{b}(z)+\velocity)^*-\At(\alpha-\beta_{s}\im)(z)},
\end{equation}
where
$$\beta_{s}(z):=\frac{-a_1(z)+s\sqrt{\Delta(z,\alpha(z))}}{2a_2}.$$
Furthermore, since $\Delta(z,\alpha(z))\gg 0$, the square root of $\Delta$ gives no problem and so the map $\shift\omega\mapsto q_s(\theta,\velocity;\shift\omega)$ is Lipschitz in this region.
In particular, we select the sign $s\in\{-1,1\}$ that minimizes $|\shift\bar{\omega}_s-\shift\omega|$.\\ 
\indent Finally, let $\aux\in\rel_M(\theta,\velocity)$ with $|\sigma_{-}|\leq|\sigma_{+}|$ and $M^2-1-B(z),1-|\shift\omega|\leq\delta$. Notice that the identity \eqref{id4} yields
$$\frac{1\pm\theta}{2}(M^2\pm\At)(1-|\sigma_{\pm}|^2)
\leq\eqref{id4:2}=\eqref{id4:1}=\mathcal{O}(\delta).$$
Hence, we can take $\delta=D(1-\theta^2)(M^2-1)$ for some constant $0<D(\At,M,\gamma)<\tfrac{1}{2}$ in such a way that $\aux\in L(\theta,\velocity)$. In addition, we can take $D$ so that the projection $\aux_0$ of $\aux$ into $\partial\B(\theta,\velocity)$ given by $\shift\omega_0=\frac{\shift\omega}{|\shift\omega|}$ also satisfies $\aux_0\in L(\theta,\velocity)$. Recall that, by construction, $b(z_0,\omega_0)=\alpha(z_0)$ since $\aux_0\in\partial\B(\theta,\velocity)$. Thus, for some $s(z)\in\{-1,1\}$,
\begin{equation}\label{ZM:Lip1}
|\shift\bar{\omega}-\shift\omega_0|
=|q_s(z)-q_s(z_0)|\lesssim|\shift\omega-\shift\omega_0|,
\end{equation}
and so
\begin{equation}\label{ZM:Lip2}
|\shift\bar{\omega}-\shift\omega|
\leq|\shift\bar{\omega}-\shift\omega_0|+|\shift\omega-\shift\omega_0|
\lesssim|\shift\omega-\shift\omega_0|=(1-|\shift\omega|).
\end{equation}
If $|\sigma_+|<|\sigma_-|$ the formulas in \textit{step 4} are slightly different but the argument does not change.
This concludes the proof.
\end{proof}

\subsection{The $\Lambda$-lamination hull}\label{sec:Lambdahull}
In this section we prove that $\Kconstrain^{lc,\Lambda}=\bar{\rel}$ and $(\Kconstrain_M)^{lc,\Lambda}=\bar{\rel}_M$.

\begin{lemma}\label{lemma:KK1} Let $z_0\in\Kconstrain$ and $z_1\in\Kconstrain^{1,\Lambda}$ satisfying $z_1-z_0\in\Lambda$. Then, the segment $[z_0,z_1]=\{z_0+\tau(z_1-z_0)\,:\,\tau\in[0,1]\}$ lies in $\bar{\rel}$.
\end{lemma}
\begin{proof}
Recall that, by Lemma~\ref{lemma:K1}\ref{Lposition}: $z_0,z_1\in\Kconstrain$ s.t. $z_1-z_0\in\Lambda$ $\Rightarrow[z_0,z_1]\subset\partial\rel$.\\
Now, let $z_0=(\theta_0,\velocity_0,\aux_0)\in\Kconstrain$ and $z_1=(\theta_1,\velocity_1,\aux_1)\in\Kconstrain^{1,\Lambda}\setminus\Kconstrain$, that is,  $|\theta_0|=1$, $|\theta_1|<1$ and
\begin{equation}\label{KK1:1}
\aux_0=\theta_0\velocity_0,\quad\quad
\aux_1=\theta_1\velocity_1+\frac{(1-(\theta_1)^2)(\At\velocity_1+\im)}{1+\bar{\omega}_1\theta_1\At}\bar{\omega}_1,
\end{equation}
for some $\bar{\omega}_1\in\Ss$. 
Let us suppose that $\bar{z}\equiv z_1-z_0\in\Lambda$, that is, $\bar{\velocity}=
\bar{\omega}(\At\bar{\aux}+\bar{\theta}\im)$ for some $\bar{\omega}\in\Ss$.
We want to show that the intermediate states $z_{\tau}\equiv z_0+\tau\bar{z}$ belong to $\bar{\rel}$ for all $\tau\in(0,1)$.
We split the proof in two steps. Firstly (\textit{step 1}) we prove the statement by assuming a claim. Secondly (\textit{step 2}) this claim is proved by elementary computations.\\
\indent\textit{Step 1. Claim}: Given $\tau\in(0,1)$, there is $\omega_{\tau}\in\R^2$ satisfying
\begin{equation}\label{KK1:2}
(1+\omega_\tau\theta_\tau\At)(\aux_\tau-\theta_\tau\velocity_\tau)
=(1-(\theta_\tau)^2)(\At\velocity_\tau+\im)\omega_\tau,
\end{equation}
if and only if
\begin{equation}\label{KK1:3}
(\At\velocity_1+\im)((\beta_\tau-\beta)\omega_{\tau}-(\beta_\tau\bar{\omega}-\beta\bar{\omega}_1))=0,
\end{equation}
where we have abbreviated
\begin{subequations}
\label{alphabeta}
\begin{align}
\beta_\tau&\equiv(\theta_1-\theta_0)\alpha_1(1-\tau),
&\alpha_1&\equiv 1-\bar{\omega}_1\theta_0\At,\\
\beta&\equiv(\theta_1+\theta_0)\alpha,
&\alpha&\equiv 1-\bar{\omega}\theta_0\At.
\end{align}
\end{subequations}
(Notice that $\alpha,\alpha_1,\beta,\beta_\tau\neq 0$).
We shall prove this equivalence in the \textit{step 2}.\\
Assume that this claim is true.
Then, if $\At\velocity_1+\im=0$, \eqref{KK1:3} holds trivially for every $\omega_{\tau}\in\Ss$ ($\Rightarrow z_{\tau}\in\partial\rel$ by Lemma~\ref{lemma:K1}\ref{Lposition}). Now let us assume that $\At\velocity_1+\im\neq0$. Hence, \eqref{KK1:3} holds if and only if
$$(\beta_\tau-\beta)\omega_{\tau}=\beta_\tau\bar{\omega}-\beta\bar{\omega}_1,$$
or equivalently (by applying the translation operator $\shift$ \eqref{shift})
\begin{equation}\label{KK1:6}
(\beta_\tau-\beta)\shift\omega_{\tau}=\beta_\tau\shift\bar{\omega}-\beta\shift\bar{\omega}_1.
\end{equation}
A priori there could be some (unique) $\tau\in (0,1)$ satisfying $\beta_\tau=\beta$. However, since $\bar{\rel}$ is closed and $\tau\mapsto z_{\tau}$ is continuous, it is enough to prove the statement for the remainder $\tau$'s satisfying $\beta_\tau\neq\beta$. For those $\tau$'s, \eqref{KK1:6} determines $\omega_{\tau}$:
$$\shift\omega_{\tau}
=\frac{\beta_\tau\shift\bar{\omega}-\beta\shift\bar{\omega}_1}{\beta_\tau-\beta}.$$
Hence, since $|\shift\bar{\omega}|=|\shift\bar{\omega}_1|=1$, we have (recall \eqref{alphabeta})
\begin{equation}
\label{KK1:4}
\begin{split}
|\shift\omega_\tau|^2
&=1+2\frac{\beta_\tau\cdot\beta-(\beta_\tau\shift\bar{\omega})\cdot(\beta\shift\bar{\omega}_1)}{|\beta_\tau-\beta|^2}\\
&=1-2(1-\tau)(1-(\theta_1)^2)\frac{\alpha_1\cdot\alpha-(\alpha_1\shift\bar{\omega})\cdot(\alpha\shift\bar{\omega}_1)}{|\beta_\tau-\beta|^2}.
\end{split}
\end{equation}
Finally, by applying
\begin{align*}
4\alpha_1\alpha^*
&=(2+(1-\shift\bar{\omega}_1)\theta_0A)(2+(1-\shift\bar{\omega}^*)\theta_0A)\\
&=(2+\theta_0\At)^2
+(\theta_0\At)^2\shift\bar{\omega}_1\shift\bar{\omega}^*
-\theta_0\At(2+\theta_0\At)(\shift\bar{\omega}_1+\shift\bar{\omega}^*),
\end{align*}
we get
\begin{align*}
\alpha_1\cdot\alpha-(\alpha_1\shift\bar{\omega})\cdot(\alpha\shift\bar{\omega}_1)
&=\Re((\alpha_1\alpha^*)(1-\shift\bar{\omega}\shift\bar{\omega}_1^*))\\
&=\underbrace{\frac{1}{4}((2+\theta_0\At)^2-(\theta_0\At)^2)}_{=1+\theta_0\At}
(1-\shift\bar{\omega}\cdot\shift\bar{\omega}_1)
\geq 0.
\end{align*}
Therefore, \eqref{KK1:4} yields $|\shift\omega_\tau|\leq 1$ ($\Rightarrow z_\tau\in\bar{\rel}$ by Lemmas~\ref{lemma:K1}\ref{Ldirection} and \ref{lemma:U}\ref{Ldirection'}).\\
\indent\textit{Step 2. Proof of the claim}:
On the one hand, $\bar{\theta}=\theta_1-\theta_0$, $\bar{\velocity}=\velocity_1-\velocity_0$ and, by \eqref{KK1:1},
\begin{equation}\label{KK1:7}
\bar{\aux}
=\aux_1-\aux_0
=\theta_0\bar{\velocity}
+\bar{\theta}\left(\velocity_1-\frac{(\theta_1+\theta_0)(\At\velocity_1+\im)}{1+\bar{\omega}_1\theta_1\At}\bar{\omega}_1\right).
\end{equation}
On the other hand, by applying \eqref{KK1:7} into the condition $\bar{\velocity}=\bar{\omega}(\At\bar{\aux}+\bar{\theta}\im)$ we get
\begin{equation}\label{KK1:8}
(\underbrace{1-\bar{\omega}\theta_0\At}_{=\alpha})\bar{\velocity}
=\bar{\omega}\bar{\theta}\frac{(\overbrace{1-\bar{\omega}_1\theta_0\At}^{=\alpha_1})(\At\velocity_1+\im)}{1+\bar{\omega}_1\theta_1\At}.
\end{equation}
Let us abbreviate $\expected{z}\equiv z_1+z_0$ and
$$\mathbf{f}\equiv\frac{\At\velocity_1+\im}{\alpha(1+\bar{\omega}_1\theta_1\At)}.$$
(Notice that: $\mathbf{f}=0\Leftrightarrow\At\velocity_1+\im=0$).
Thus, \eqref{KK1:7}\eqref{KK1:8} read as
$$
\bar{\velocity}
=\bar{\theta}\alpha_1\bar{\omega}\mathbf{f},\quad\quad
\bar{\aux}
=\theta_0\bar{\velocity}+\bar{\theta}(\velocity_1-\expected{\theta}\alpha\bar{\omega}_1\mathbf{f}).
$$
Let us expand the factors of \eqref{KK1:2} in terms of $\tau$. They are
$$\aux_\tau-\theta_\tau\velocity_\tau
=\underbrace{\aux_0-\theta_0\velocity_0}_{=0}
+\tau(\underbrace{\bar{\aux}-\theta_0\bar{\velocity}-\bar{\theta}\velocity_0}_{=\bar{\theta}(\bar{\velocity}-\expected{\theta}\alpha\bar{\omega}_1\mathbf{f})})-\tau^2\bar{\theta}\bar{\velocity}
=\tau\bar{\theta}(\underbrace{(1-\tau)\dev{\theta}\alpha_1}_{=\beta_\tau}\bar{\omega}-\underbrace{\expected{\theta}\alpha}_{=\beta}\bar{\omega}_1)\mathbf{f},$$
$$1-(\theta_\tau)^2
=(\underbrace{\theta_0-\theta_{\tau}}_{=-\tau\bar{\theta}})(\theta_0+\theta_{\tau}),$$
and
\begin{align*}
\At\velocity_\tau+\im
&=(\At\velocity_1+\im)-\At(1-\tau)\bar{\velocity}\\
&=(\alpha(1+\bar{\omega}_1\theta_1\At)
-\At(1-\tau)\bar{\theta}\alpha_1\bar{\omega})\mathbf{f}\\
&=(\alpha\alpha_1-\At(\underbrace{(1-\tau)\dev{\theta}\alpha_1}_{=\beta_\tau}\bar{\omega}-\underbrace{\expected{\theta}\alpha}_{=\beta}\bar{\omega}_1))\mathbf{f}.
\end{align*}
Hence, the equation \eqref{KK1:2} reads as
$$(1+\omega_\tau\theta_\tau\At)
\tau\bar{\theta}(\beta_\tau\bar{\omega}-\beta\bar{\omega}_1)\mathbf{f}
=\tau\bar{\theta}(\theta_0+\theta_{\tau})(\At(\beta_\tau\bar{\omega}-\beta\bar{\omega}_1)-\alpha\alpha_1)\mathbf{f}\omega_\tau,$$
or equivalently ($\tau\bar{\theta}\neq 0$)
\begin{equation}\label{KK1:5}
(\beta_\tau\bar{\omega}-\beta\bar{\omega}_1)\mathbf{f}
=(\theta_0\At(\beta_\tau\bar{\omega}-\beta\bar{\omega}_1)
-(\theta_0+\theta_{\tau})\alpha\alpha_1)\omega_{\tau}\mathbf{f}.
\end{equation}
Finally, by splitting $(\theta_0+\theta_{\tau})=\expected{\theta}-(1-\tau)\dev{\theta}$, we have
$$(\theta_0+\theta_{\tau})\alpha\alpha_1
=(1-\bar{\omega}_1\theta_0\At)\underbrace{\expected{\theta}\alpha}_{=\beta}
-(1-\bar{\omega}\theta_0\At)\underbrace{(1-\tau)\bar{\theta}\alpha_1}_{=\beta_\tau},$$
and so \eqref{KK1:5} is equivalent to \eqref{KK1:3}.
\end{proof}

\begin{prop}\label{KL2} 
$\Kconstrain^{lc,\Lambda}
=\Kconstrain^{2,\Lambda}
=\bar{\rel}.$
\end{prop}
\begin{proof} 
Firstly (\textit{step 1}) we prove that $\Kconstrain^{2,\Lambda}
=\bar{\rel}$. Secondly (\textit{step 2}) we deduce that $\Kconstrain^{lc,\Lambda}
=\Kconstrain^{2,\Lambda}$.\\
\indent\textit{Step 1}. $\Kconstrain^{2,\Lambda}
=\bar{\rel}$: 
Since $\rel$ is open and $\partial\rel=\Kconstrain^{1,\Lambda}$ (Lemma~\ref{lemma:K1}), this is equivalent to prove that $\Kconstrain^{2,\Lambda}\setminus\Kconstrain^{1,\Lambda}=\rel$.\\
By definition \eqref{def:K1L} and Lemma~\ref{lemma:K1}\ref{gA} a state $z=(\theta,\velocity,\aux)\in\set$ belongs to $\Kconstrain^{2,\Lambda}\setminus\Kconstrain^{1,\Lambda}$ if and only if $g(z)\neq 0$ and there are $0\neq\bar{z}\in\Lambda$ and $\lambda_-<0<\lambda_+$ satisfying
$$
|\theta_{\lambda_{\pm}}|\leq 1,\quad\quad
g(z_{\lambda_{\pm}})=0,
$$
where $z_\lambda\equiv z+\lambda\bar{z}$. Since $\bar{z}\in\Lambda$, notice that
$$\det\mathbf{T}(z_{\lambda})
=\textrm{quadratic}+\lambda^3\underbrace{\det\mathbf{T}(\bar{z})}_{=0}.$$
Then, by Lemma~\ref{lemma:K1}\ref{gA}, the polynomial $p:\lambda\mapsto g(z_{\lambda})$ is cubic
\begin{equation}\label{KL2:1}
p(z,\bar{z};\lambda):=g(z_\lambda)=\sum_{j=0}^3a_j(z,\bar{z})\lambda^j.
\end{equation}
\indent\textit{Step 1.1}. $\bar{\rel}\subset\Kconstrain^{2,\Lambda}$: The analysis of \eqref{KL2:1} is easier for $\bar{z}\in\varLambda_0$ because $p$ is quadratic ($a_3=0$) in such case.
Moreover, since $\bar{\theta}=0$ and $\bar{\velocity}=-\At\bar{\aux}$, the second coefficient is strictly positive
$$a_2=(\bar{\aux}+\At\bar{\velocity}-\theta(\bar{\velocity}+\At\bar{\aux}))\cdot(\bar{\aux}-\theta\bar{\velocity})=(1-\At^2)(1+\theta\At)|\bar{\aux}|^2>0.$$
Hence, $p$ has two real roots of different sign if and only if $g(z)=a_0<0$ ($z\in\rel$). Therefore, $\bar{\rel}=(\Kconstrain^{1,\Lambda})^{1,\varLambda_0}\subset\Kconstrain^{2,\Lambda}$ 
(As a curiosity observe that, since $g$ is $\varLambda_0$-convex, $\bar{\rel}=(\Kconstrain^{1,\Lambda})^{\varLambda_0}$).\\
\indent\textit{Step 2}: $\Kconstrain^{2,\Lambda}\subset\bar{\rel}$.
Since $\Kconstrain^{1,\Lambda}=\Kconstrain^{1,\varLambda_1}$ (Lemma~\ref{lemma:K1}), by the \textit{step 1} we only need to check that $\Kconstrain^{2,\varLambda_1}\setminus\Kconstrain^{1,\varLambda_1}\subset\rel$.\\
Let $z=(\theta,\velocity,\aux)\in\Kconstrain^{2,\varLambda_1}\setminus\Kconstrain^{1,\varLambda_1}$. 
By hypothesis, $g(z)\neq 0$ and there are $0\neq\bar{z}\in\varLambda_1$ with $\bar{\theta}=1$ and $\lambda_-<0<\lambda_+$ satisfying $|\theta+\lambda_{\pm}|\leq 1$ and $p(z,\bar{z};\lambda_{\pm})=g(z_{\lambda_{\pm}})=0$. 
Notice that necessarily $|\theta|<1$.
If we abbreviate $z^{\pm}\equiv z_{(\pm 1-\theta)}=z+(\pm 1-\theta)\bar{z}$, then $\theta^{\pm}=\pm 1$ and Lemma~\ref{lemma:K1}\ref{gA} yields
$$p(z,\bar{z};\pm 1-\theta)
=g(z^{\pm})
=(1\mp\At)|\aux^\pm\mp\velocity^\pm|^2\geq 0.$$
If both $p(z,\bar{z};\pm 1-\theta)>0$ necessarily $g(z)=p(z,\bar{z};0)<0$, otherwise we would deduce that $p'(z,\bar{z};\cdot)$ has at least $3$ roots in $[-1-\theta,1-\theta]$. If both $p(z,\bar{z};\pm 1-\theta)=0$ we would have $z^{\pm}=(\pm 1,\velocity^{\pm},\pm\velocity^{\pm})\in\Kconstrain$, and so $z\in\Kconstrain^{1,\varLambda_1}$. If only one of $p(z,\bar{z};\pm 1-\theta)$ is zero, then $z$ is a $\Lambda$-convex combination of a state in $\Kconstrain$ and other in $\Kconstrain^{1,\varLambda_1}\setminus\Kconstrain$. Thus, by Lemma~\ref{lemma:KK1}, $z\in\bar{\rel}$.\\
\indent\textit{Step 2}. $\Kconstrain^{lc,\Lambda}=\Kconstrain^{2,\Lambda}$: It is a general fact in Lamination Theory that, for any closed $K$, the following holds: $K^{1,\Lambda}\setminus K
=(\partial K)^{1,\Lambda}\setminus K$.
Hence, since $\partial(\Kconstrain^{2,\Lambda})
=\partial\rel
=\Kconstrain^{1,\Lambda}$, we deduce that
$\Kconstrain^{3,\Lambda}\setminus\Kconstrain^{2,\Lambda}=\emptyset$. Therefore, inductively $\Kconstrain^{n,\Lambda}=\Kconstrain^{2,\Lambda}$ for all $n\geq 3$.
\end{proof}

\begin{prop}\label{KML} Let $M>1$. Then $(\Kconstrain_M)^{lc,\Lambda}=\bar{\rel}_M$.
\end{prop}
\begin{proof} \textit{Step 1. $(\Kconstrain_M)^{lc,\Lambda}\subset\bar{\rel}_M$}: It follows from: $\bar{\rel}$ is $\Lambda$-lamination convex, \eqref{id2:2} defines the sublevel set of a $\Lambda$-convex (indeed $\Lambda$-affine) function, and \eqref{id2:3}-\eqref{id2:4} define sublevel sets of convex functions.\\
\indent\textit{Step 2. $\bar{\rel}_M\subset(\Kconstrain_M)^{lc,\Lambda}$}: As in \cite{RIPM}, it follows from the Krein-Milman
type theorem in the context of $\Lambda$-convexity \cite[Lemma 4.16]{Rigidity}, because, as we saw in Lemma~\ref{lemma:ZMseg}, for all $z\in\partial\rel_M\setminus\Kconstrain_M$ there is $0\neq\bar{z}\in\Lambda$ such that $z\pm\bar{z}\in\bar{\rel}_M$ (i.e.~$z$ is not an extreme point of $\bar{\rel}_M$). More precisely, let $z=(\theta,\velocity,\aux)\in\partial\rel_M\setminus\Kconstrain_M$. As in \textit{step 1} in the proof of Lemma~\ref{lemma:ZMseg}, we take $\bar{z}$ in terms of $(\bar{\mathbf{n}},\bar{\omega})\in\R^2\times\Ss$ to be determined. 
If $\aux\in\partial\B$ ($\omega\in\Ss$) it is enough to take $\bar{\omega}=\omega$. Otherwise ($\aux\notin\partial\B$) we may assume w.l.o.g.~that $\aux\notin\partial\B_{-}$. If $\aux\in\partial\HP$ we take $\bar{\omega}$ satisfying $b(z,\bar{\omega})=0$ \eqref{ZM:6}. If $\aux\in\partial\B_{+}$ we take $\bar{\mathbf{n}}$ as in \eqref{ZM:8}.
\end{proof}

\begin{Rem}
Notice that we are not excluding the case $M=M_*(\At)$ in Proposition~\ref{KML}. Although we believe that in this case Lemma~\ref{lemma:ZMseg} holds too, we have chosen to exclude it in Lemma~\ref{lemma:ZMseg} for simplicity.
\end{Rem}

\begin{Rem}\label{Rem:KLambda?}
In \cite{RIPM} the identity $\bar{\rel}_0=\Kconstrain^{\Lambda_0}$
(and also $\bar{\rel}_{0,M}=(\Kconstrain_{0,M})^{\Lambda_0}$) follows from the fact that $f_0$ is $\Lambda_0$-convex. However, $f_\At$ (Lemma~\ref{lemma:K1}\ref{fA}) is not $\Lambda_\At$-convex for $\At\neq 0$: Let $z_0=(0,-\im/\At,0)\in\bar{\rel}$ and $\bar{z}_0=(1,0,0)\in\Lambda$ ($(\bar{\aux},\bar{\omega})=(0,0)\in\R^2\times\Ss$). Then, the function
\begin{equation}\label{counterexample}
h_{\At}(\lambda):=f_{\At}(z_0+\lambda\bar{z}_0)=2|1-\lambda\At|\frac{|\lambda|}{|\At|},
\end{equation}
is not convex since
$$\partial_{\lambda}^2h_{\At}(\lambda)
=-4\mathrm{sgn}(\lambda\At),
\quad 0<|\lambda|<1/\At.
$$
Notice that this does not imply that $\bar{\rel}_\At\subsetneq\Kconstrain^{\Lambda_\At}$. In general, $\bar{\rel}_{\At}$ can be expressed as $\{z\in\set\,:\,c_{\At}(z)f_{\At}(z)\leq 0\}$ for all $c_{\At}>0$ on $\set$. Thus, to prove that $\bar{\rel}_\At=\Kconstrain^{\Lambda_\At}$ it is enough to find a correcting factor $c_\At>0$ making $c_\At f_{\At}$ $\Lambda_{\At}$-convex on $\set$. For instance, $c_{\At}(z)=1/(1-\theta\At)$ repairs the counterexample \eqref{counterexample} since $(c_\At f_{\At})(z_0+\lambda\bar{z})=2|\lambda|/|\At|$.  However, it seems hard to check if $c_\At f_{\At}$ is $\Lambda_{\At}$-convex.
Still we conjecture that $\bar{\rel}_{\At}$ is indeed $\Kconstrain^{\Lambda_{\At}}$ and also closed under weak*-convergence, thus representing the full relaxation of $(\textrm{IPM}_{\At})$ in analogy with the case $A=0$.
\end{Rem}

\appendix
\section{Toy random walk}\label{sec:RW}

\begin{figure}[h!]
	\centering
	\includegraphics[height=5.67cm]{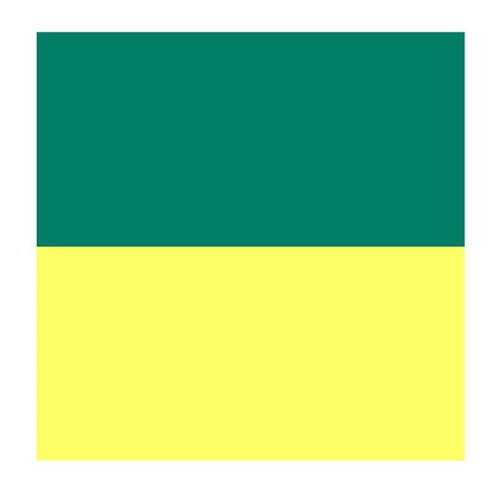}
	\caption{The unstable planar phase distribution.}
	\label{fig:theta0}
\end{figure}

\begin{figure}[h!]\centering
	\includegraphics[height=6.87cm]{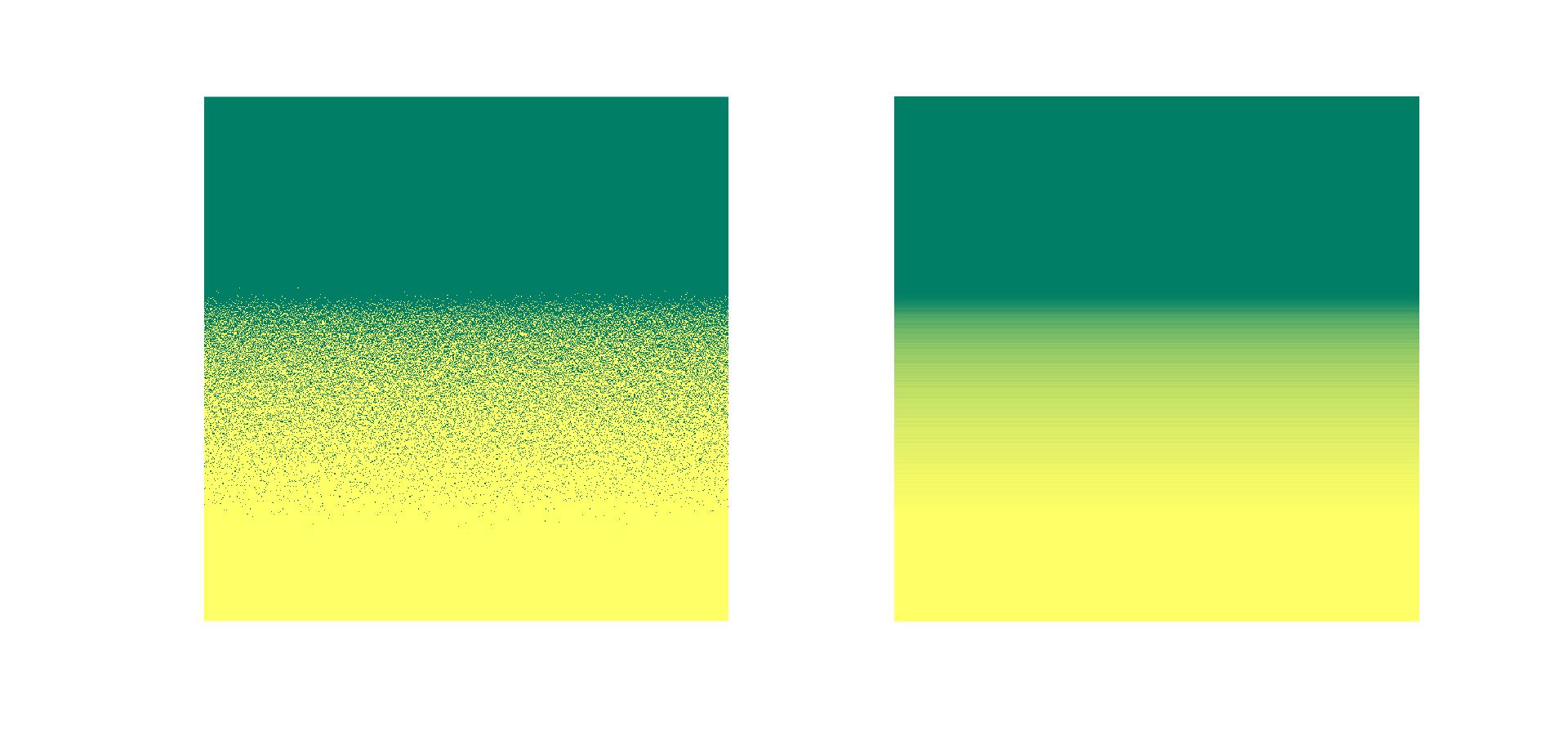}
	\includegraphics[height=6.87cm]{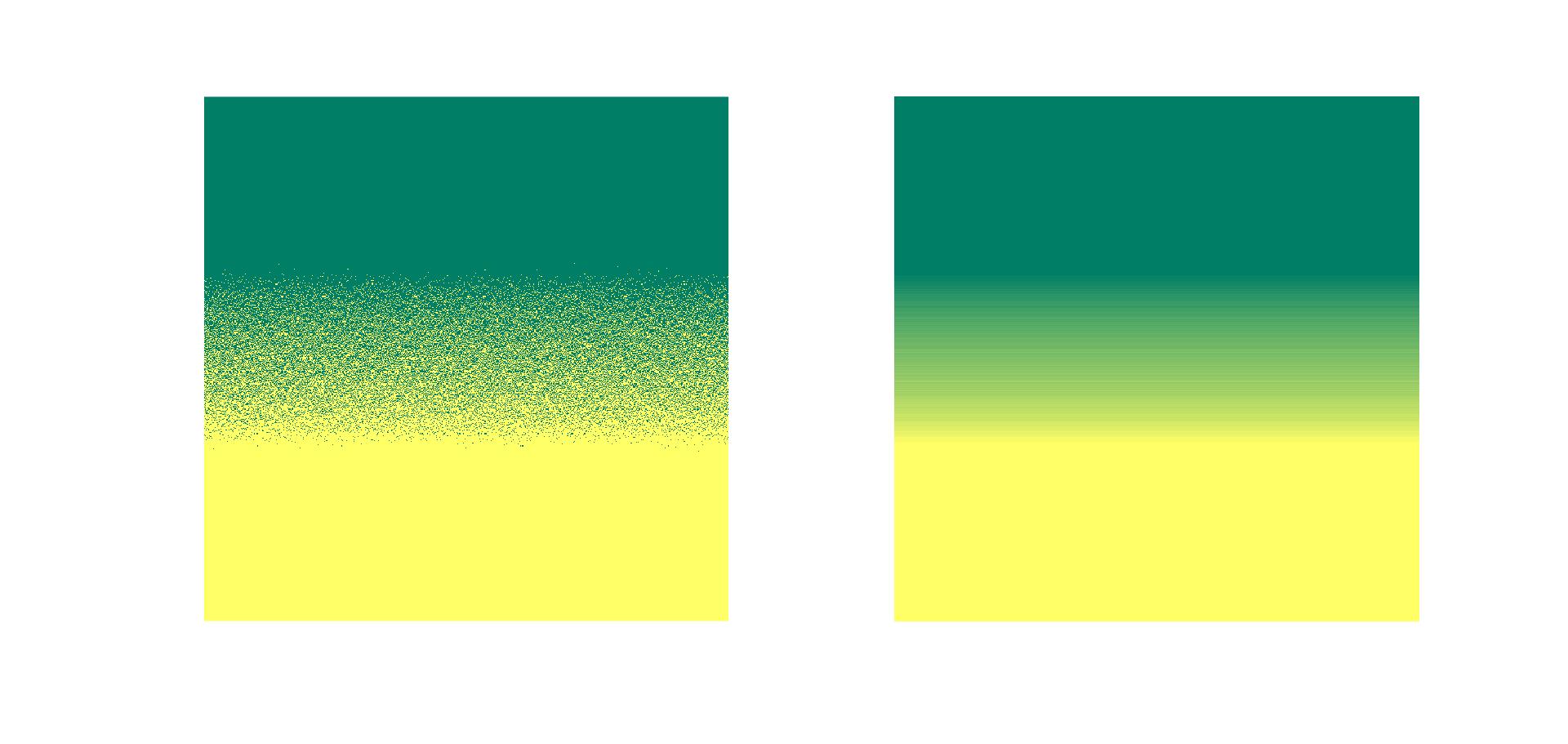}
	\includegraphics[height=6.87cm]{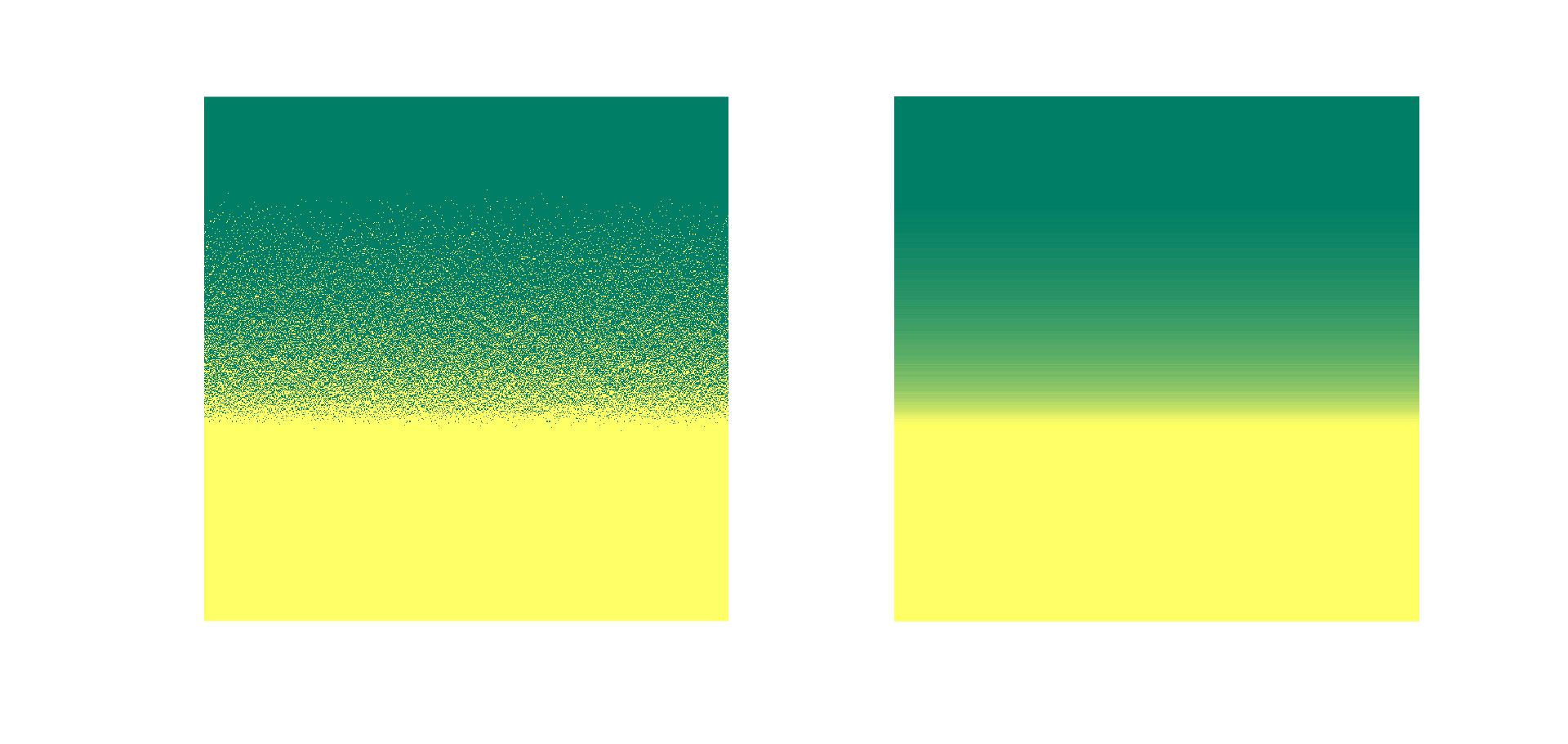}
	\caption{On the left hand column we see a Matlab simulation (``solution'') of this random walk stopped at some time starting from Fig.~\ref{fig:theta0}, while the right hand column shows  the average over lines (``subsolution'') of the previous picture. From top to bottom, the corresponding Atwood number $\At$ is $-\tfrac{1}{2}$, $0$ and $\tfrac{1}{2}$ respectively (cf.~Fig.~\ref{fig:ThetaA}).}
	\label{fig:mixprofile}
\end{figure}

\begin{figure}[h!]\centering
	\includegraphics[width=9.0cm]{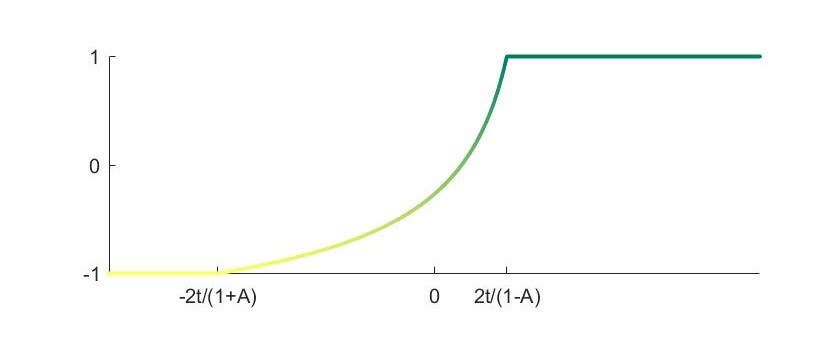}\\[0.3cm]
	\includegraphics[width=9.0cm]{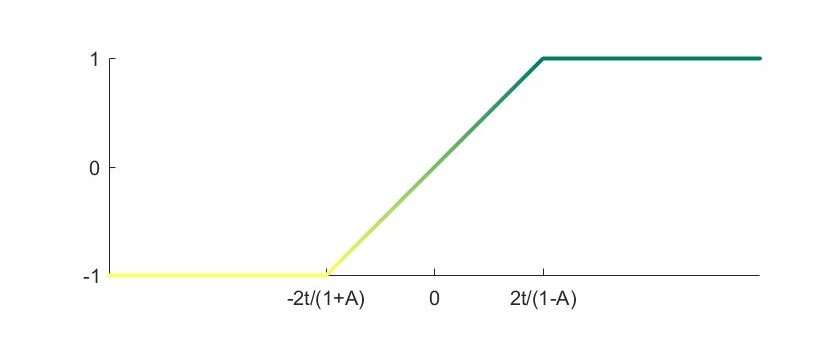}\\[0.3cm]
	\includegraphics[width=9.0cm]{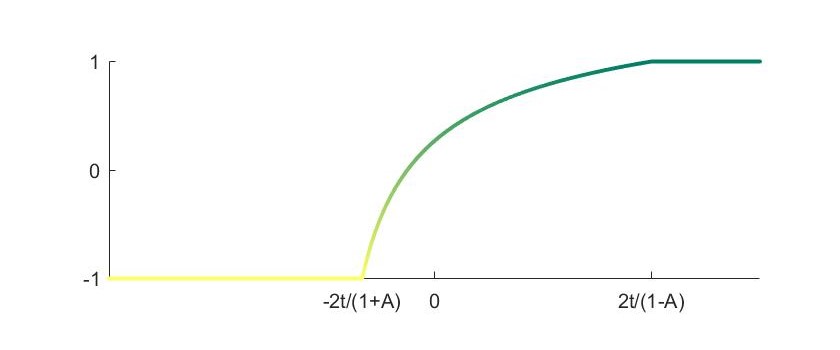}
	\caption{From top to bottom, we see the mixing profile $\Theta(t,x_2)$ at time $t=\tfrac{1}{2}$ for the Atwood number $\At$ equal to $-\tfrac{1}{2}$, $0$ and $\tfrac{1}{2}$ respectively.}
	\label{fig:ThetaA}
\end{figure}

\begin{figure}[h!]\centering
	\includegraphics[height=6.7cm]{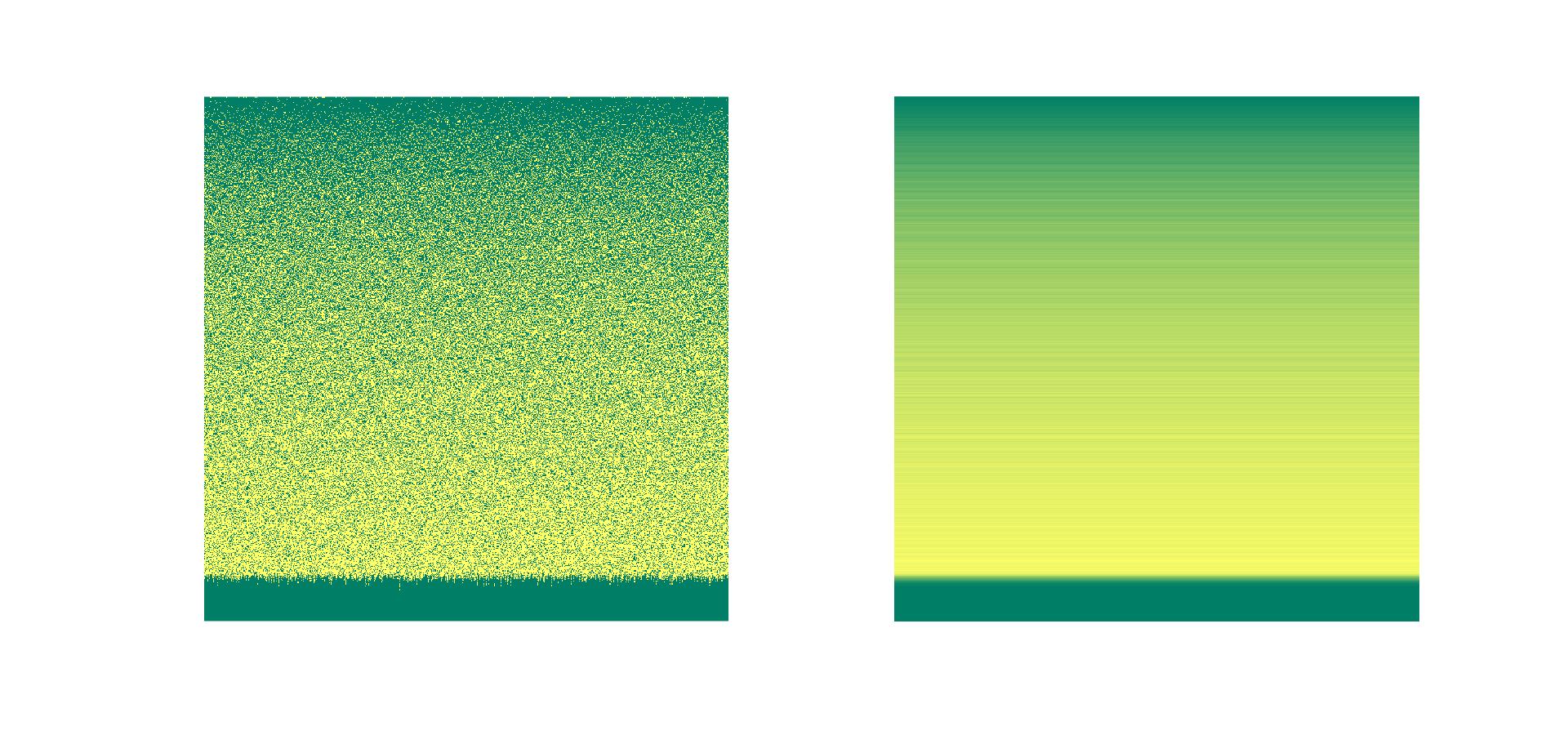}
	\caption{Evolution of Figure~\ref{fig:mixprofile} for $\At=-\frac{1}{2}$ at some $(c_{\At}^{-})^{-1}<t<(c_{\At}^{+})^{-1}$.}
	\label{fig:lower}
\end{figure}

\begin{figure}[h!]\centering
	\includegraphics[height=6.7cm]{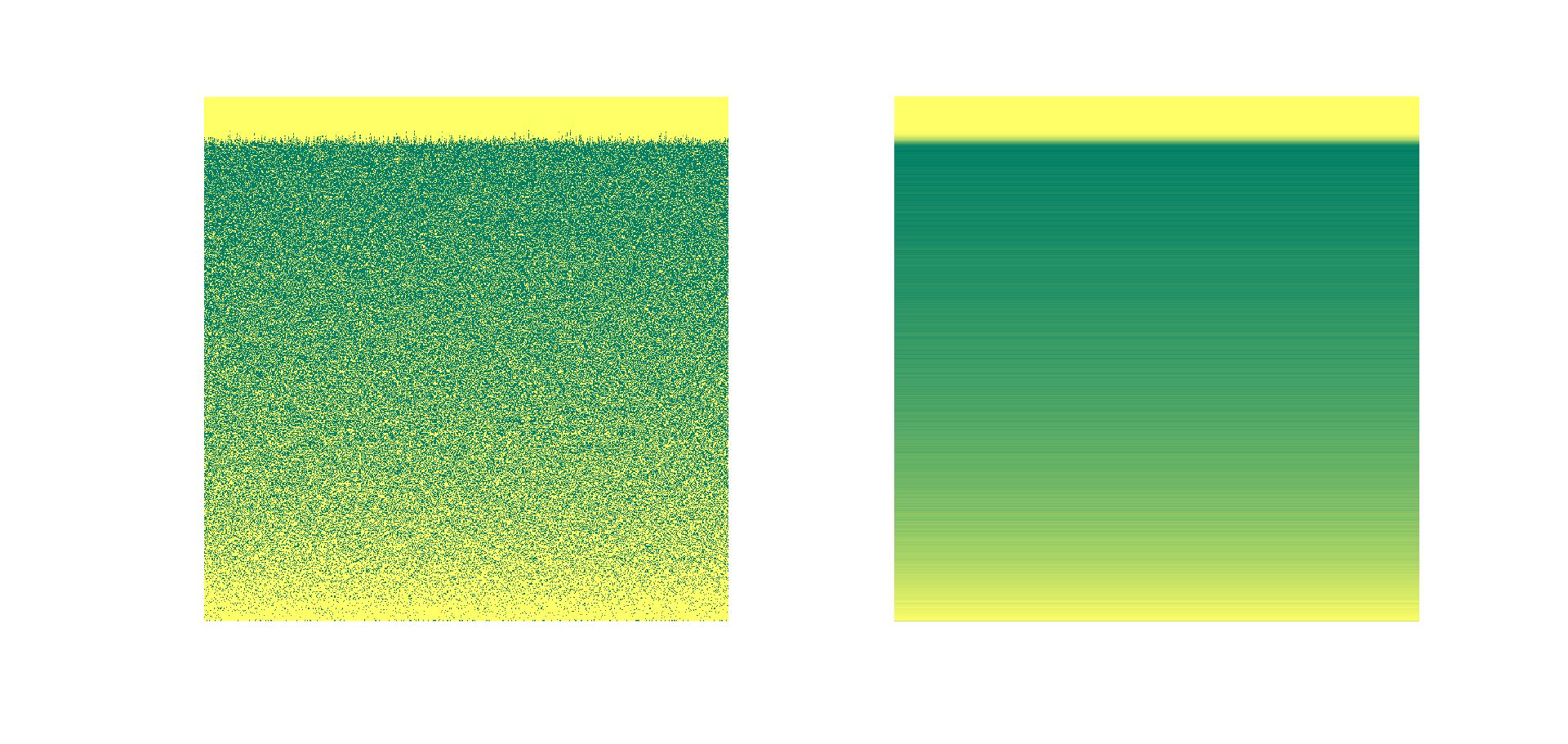}
	\caption{Evolution of Figure~\ref{fig:mixprofile} for $\At=\frac{1}{2}$ at some $(c_{\At}^{+})^{-1}<t<(c_{\At}^{-})^{-1}$.}
	\label{fig:upper}
\end{figure}

\begin{figure}[h!]\centering
	\includegraphics[height=6.7cm]{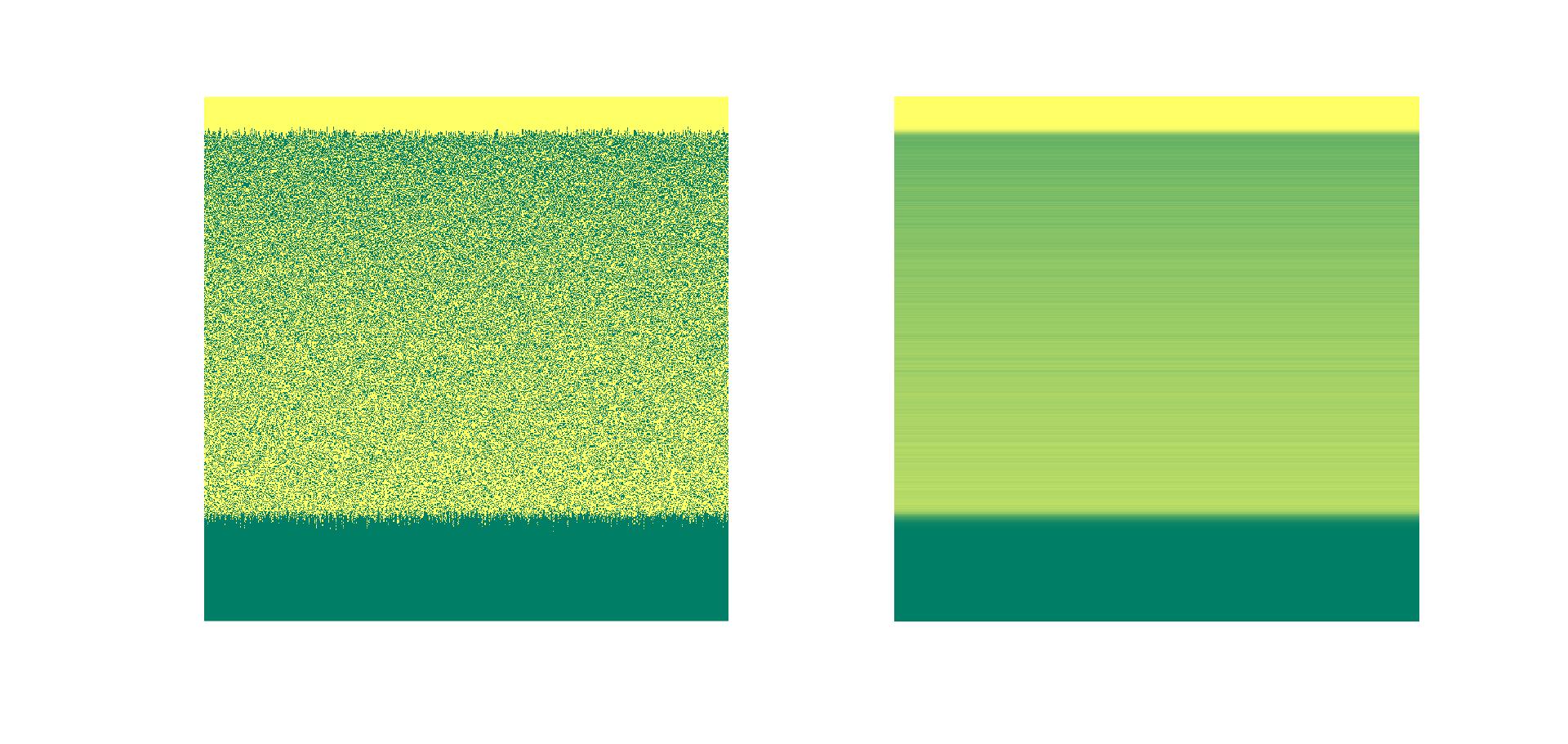}
	\includegraphics[height=6.7cm]{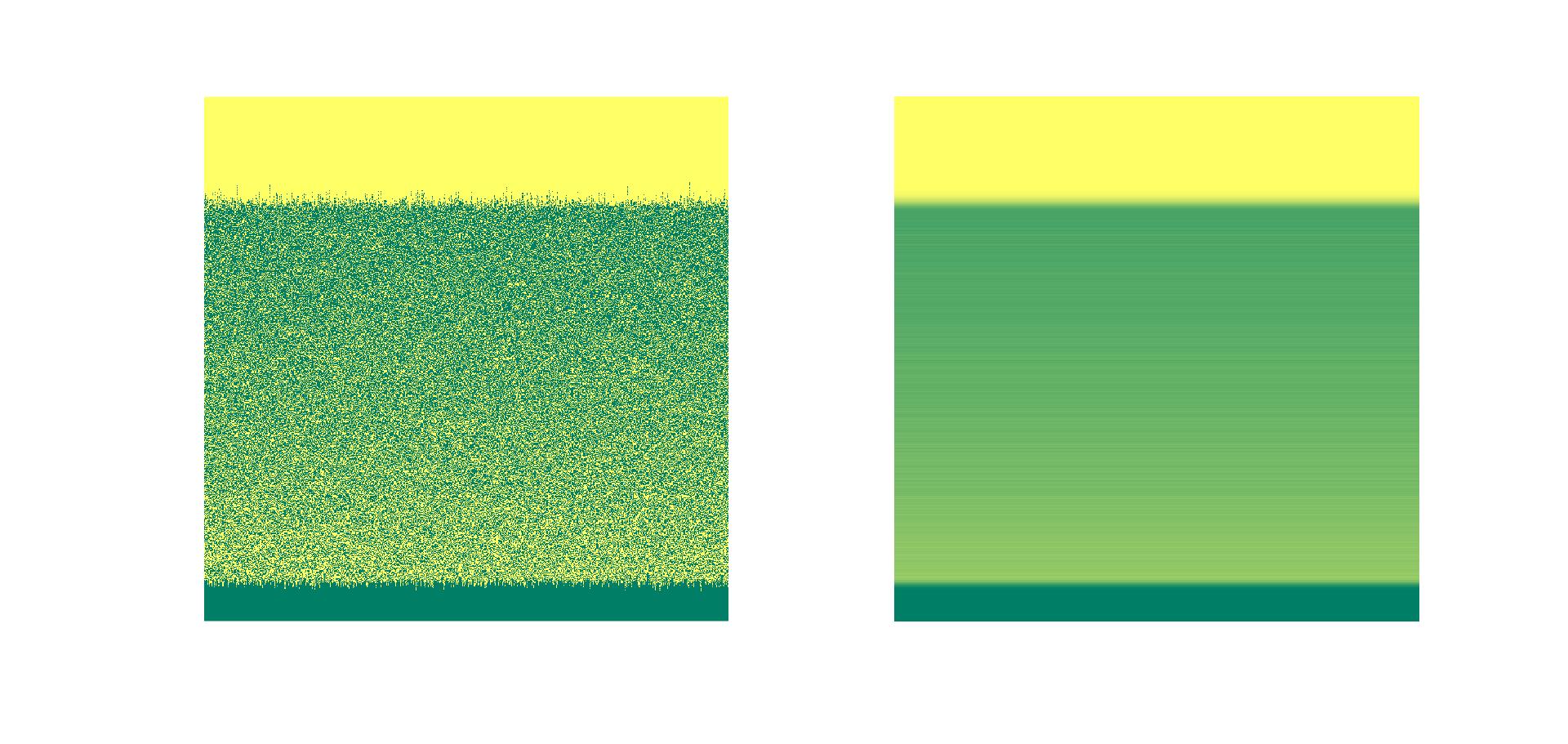}
	\caption{Evolution of Figure~\ref{fig:mixprofile} for $\At=-\frac{1}{2},\frac{1}{2}$ at some $t>(c_{\At}^{-})^{-1}\vee(c_{\At}^{+})^{-1}$.}
	\label{fig:after}
\end{figure}

In this section we introduce a toy random walk to illustrate how these $\Theta_{\At}$-mixing solutions may look like (see Fig.~\ref{fig:theta0}-\ref{fig:after}) and, at the same time, to give somehow an intuitive idea of the interplay between the unpredictable nature at the microscopic level of the mixing phenomenon and the deterministic point of view at the mesoscopic scale. This is also motivated by the relaxation approach of Otto \cite[\S2]{O}: \\

\textbf{Otto's approach.} Roughly speaking, by passing from the Eulerian (phase $\theta(t,\x)$) to the Lagrangian (flow map $\Phi(t,\x)$) point of view, Otto rewrote the Muskat problem as a gradient flux for $\Phi$ w.r.t.~the gravitational potential energy $E$ with the following physical interpretation: ``Given $\theta_0$ \eqref{flat}, the phase distribution $\theta$ advected by the flow ($\theta(t,\Phi(t))=\theta_0$) aims at minimizing $E$ by transforming it into kinetic energy, which then is dissipated by friction when forcing the fluid through the porous medium''. A natural discretization in time intervals of size $h$ yields a recurrence $\Phi_h^k\rightsquigarrow\Phi_h^{k+1}$ starting from $\Phi_h^0=\mathrm{id}$ that leads an approximate time-discrete solution $\Phi_h=(\Phi_h^k)_k$, where $\Phi_h^{k+1}$ is the unique solution of a variational problem defined in terms of $\Phi_h^k$. 
As he noted, $\Phi_h^1$ is not one-to-one, thus preventing (a priori) from
defining the corresponding $\theta_h^1$ by advection. %$\theta_h^1=\theta_0\circ(\Phi_h^1)^{-1}$. 
Nevertheless, by subdividing the space in a grid of size $r$, each $\Phi_h^k$ can be approximated by a (minimizing) sequence of permutations $\Phi_{h,r}^k$ of this partition. Then, each $\Phi_{h,r}^k$ defines a $\{-1,1\}$-valued discrete phase distribution $\theta_{h,r}^k=(\theta_0)_r\circ(\Phi_{h,r}^k)^{-1}$ where $(\theta_0)_r$ is a sampling of $\theta_0$. It is interesting to point that $\Phi_{h,r}^1$ (and so $\theta_{h,r}^1$) breaks the planar symmetry of \eqref{flat} and consequently is not unique. Despite this lack of uniqueness, Otto showed that $\theta_{h,r}^k\overset{*}{\rightharpoonup}\theta_h^k=(\Phi_h^k)^\sharp(\theta_0)$ $\equiv$ push-forward of $\theta_0$ under $\Phi_h^k$, which allows to interpret $\theta_h^k$ as the average in space of the actual phase distribution.
At the same time, $\theta_h^k$ is the unique solution of a convex variational problem, linked with the one for $\Phi_h^k$ through Optimal Transport Theory. To conclude Otto proved that $\theta_h$ converges in $L^\infty_tL^1$ to the (unique) entropy solution $\Theta_{\At}$ \eqref{entropy} of the conservation law \eqref{conservationlaw}.\\

\textbf{Toy random walk}. As in \cite{O}, we discretize in time intervals of size $h=\triangle t$ and we subdivide the domain in a grid of size $r=\triangle x_i$ whose center points form the lattice $r(\Z^2+\tfrac{1}{2}\im)=\{\x_{s,j}\equiv r(s,j+\tfrac{1}{2})\,:\,s,j\in\Z\}$. Take a sample of $\theta_0$ \eqref{flat}
\begin{equation}\label{flatdis}
\theta^{(0)}(\x_{s,j})=
\left\lbrace\begin{array}{cl}
+1, & j>0,\\[0.1cm]
-1, & j<0.
\end{array}\right.
\end{equation} 

\noindent Then, we interpret the conservation of mass and volume by setting that two close
different ``molecules'' may interchange their positions if the heavier is above the lighter, i.e.~if their state is unstable due to gravity. Darcy's law is interpreted by setting that such
interchange happens with some probability 
\begin{equation}\label{prob:1}
p_{j}^{(k)}\equiv\textrm{probability of interchange between lines $j\leftrightarrow j-1$ at time $k+1$},
\end{equation}
depending on the Atwood number $\At$ and in terms of the
proximity to the rest molecules of the same fluid respectively. Note that, by simplicity, we are considering $p_j^{(k)}$ independent of $s$ due to the planar symmetry of \eqref{flatdis}. This induces a time-discrete stochastic process $\{\theta_{s,j}^{(k)}\}_{k\geq 0}$ where $\theta_{s,j}^{(k)}\equiv\theta^{(k)}(\x_{s,j})$ is the $\{-1,1\}$-valued random variable. In this way, \eqref{prob:1} reads as
$$p_{j}^{(k)}=P(\textrm{interchange }\cmark\,|\,\theta_{s,j}^{(k)}=1,\,\theta_{s,j-1}^{(k)}=-1),$$
while the probability of interchange in the remaining situations is zero. We are interested in the deterministic value
\begin{equation}\label{d-}
\breve{\theta}^{(k)}_j:=E(\theta_{s,j}^{(k)})
=d_{j,+}^{(k)}-d_{j,-}^{(k)},
\end{equation}
where
$$d_{j,\pm}^{(k)}
:=P(\theta_{s,j}^{(k)}=\pm 1).$$
This $d_{j}^{(k)}$ can be computed recursively
\begin{align*}
d_{j,+}^{(k+1)}&=\underbrace{d_{j,+}^{(k)}d_{j-1,+}^{(k)}}_{\textrm{stable}}+
\underbrace{(1-p_j^{(k)})d_{j,+}^{(k)}d_{j-1,-}^{(k)}}_{\substack{\textrm{unstable} \\ \textrm{interchange }\xmark}}
+\underbrace{p_{j+1}^{(k)}d_{j+1,+}^{(k)}d_{j,-}^{(k)}}_{\substack{\textrm{unstable} \\ \textrm{interchange }\cmark}}\\
&=\underbrace{d_{j,+}^{(k)}}_{\xmark}+\underbrace{p_{j+1}^{(k)}d_{j+1,+}^{(k)}d_{j,-}^{(k)}}_{\substack{\textrm{interchange }\cmark \\ +\textrm{ increases}}}-\underbrace{p_{j}^{(k)}d_{j,+}^{(k)}d_{j-1,-}^{(k)}}_{\substack{\textrm{interchange }\cmark \\ +\textrm{ reduces}}},
\end{align*}
and analogously
\begin{align*}
d_{j,-}^{(k+1)}&=\underbrace{d_{j+1,-}^{(k)}d_{j,-}^{(k)}}_{\textrm{stable}}+
\underbrace{(1-p_{j+1}^{(k)})d_{j+1,+}^{(k)}d_{j,-}^{(k)}}_{\substack{\textrm{unstable} \\ \textrm{interchange }\xmark}}
+\underbrace{p_{j}^{(k)}d_{j,+}^{(k)}d_{j-1,-}^{(k)}}_{\substack{\textrm{unstable} \\ \textrm{interchange }\cmark}}\\
&=\underbrace{d_{j,-}^{(k)}}_{\xmark}-\underbrace{p_{j+1}^{(k)}d_{j+1,+}^{(k)}d_{j,-}^{(k)}}_{\substack{\textrm{interchange }\cmark \\ -\textrm{ reduces}}}+\underbrace{p_{j}^{(k)}d_{j,+}^{(k)}d_{j-1,-}^{(k)}}_{\substack{\textrm{interchange }\cmark \\ -\textrm{ increases}}}.
\end{align*}
In summary, the dynamic is given by
\begin{equation}\label{dynd}
d_{j,\pm}^{(k+1)}=d_{j,\pm}^{(k)}\pm(p_{j+1}^{(k)}d_{j+1,+}^{(k)}d_{j,-}^{(k)}-p_{j}^{(k)}d_{j,+}^{(k)}d_{j-1,-}^{(k)}).
\end{equation}
Then, by \eqref{d-} and $d_{j,+}^{(k)}+d_{j,-}^{(k)}=1$ we get
\begin{equation}\label{A:1}
d_{j,\pm}^{(k)}=\tfrac{1}{2}(1\pm\breve{\theta}_{j}^{(k)}),
\end{equation} 
and consequently the recurrence \eqref{dynd} can be written in terms of $\breve{\theta}_{j}^{(k)}$ as
\begin{equation}\label{dyntheta} 
\breve{\theta}_{j}^{(k+1)}=\breve{\theta}_{j}^{(k)}+\tfrac{1}{2}(p_{j+1}^{(k)}(1+\breve{\theta}_{j+1}^{(k)})(1-\breve{\theta}_{j}^{(k)})-p_{j}^{(k)}(1+\breve{\theta}_{j}^{(k)})(1-\breve{\theta}_{j-1}^{(k)})).
\end{equation}

\noindent With \cite{O} in mind, we declare
\begin{equation}
p_{j}^{(k)}=\frac{1}{2}\frac{\mu^+\wedge\mu^-}{d_{j,-}^{(k)}\mu^++d_{j,+}^{(k)}\mu^-}\in[0,\tfrac{1}{2}].
\end{equation}
In the balanced case $\mu^+=\mu^-$ ($\At=0$), we have $p_{j}^{(k)}=\tfrac{1}{2}$ independently of $j,k$. In the case of viscosity jump $\mu^+\neq\mu^-$ ($0<|\At|<1$), the probability of interchange at time $k+1$ depends on the relative position in terms of the mobility quotient $B=\mu^+/\mu^-$ (cf.~\S\ref{sec:Otto}). For instance, when $\mu^+>\mu^-$ the lighter molecules rise through the heavier ones without many difficulties ($p_{j}^{(k)}\uparrow\tfrac{1}{2}$ as $\theta_{j}^{(k)}\uparrow 1$), whereas the molecules of the heavier fluid sink with lower speed because the fluid with phase $+$ has smaller mobility ($p_{j}^{(k)}\downarrow\tfrac{1}{2}B^{-1}$ as $\theta_{j}^{(k)}\downarrow -1$). The case $\mu^+<\mu^-$ follows analogously ($p_{j}^{(k)}\downarrow\tfrac{1}{2}B$ as $\theta_{j}^{(k)}\uparrow 1$ and $p_{j}^{(k)}\uparrow\tfrac{1}{2}$ as $\theta_{j}^{(k)}\downarrow -1$).
A simple calculation yields
\begin{equation}
p_{j}^{(k)}=\frac{a}{1-\breve{\theta}_{j}^{(k)}\At}
\quad\textrm{where}\quad
a=\frac{\mu^+\wedge\mu^-}{\mu^++\mu^-}
=\frac{1-|\At|}{2}
=\frac{1}{c_\At^+\vee c_\At^-}.
\end{equation}
Thus, if we scale the discretization as $r=ch$ for some $c>0$, the recurrence \eqref{dyntheta} can be written as a finite difference equation
\begin{equation}\label{fde}
\frac{\breve{\theta}_{j}^{(k+1)}-\breve{\theta}_{j}^{(k)}}{\triangle t}
=ca\left(\frac{(1+\breve{\theta}_{j+1}^{(k)})(1-\breve{\theta}_{j}^{(k)})}{1-\breve{\theta}_{j+1}^{(k)}\At}-\frac{(1+\breve{\theta}_{j}^{(k)})(1-\breve{\theta}_{j-1}^{(k)})}{1-\breve{\theta}_{j}^{(k)}\At}\right)\Big/\triangle x_2.
\end{equation}
Notice that, by construction, there is not interchange of molecules outside $\{(t,\x)\,:\, |x_2|<ct\}$. When $h\downarrow0$, the scheme \eqref{fde}
converges formally to the Burgers type equation \eqref{conservationlaw:alpha}
where $\alpha=ca$ is the mixing speed. Since $0<\alpha<1$, necessarily
$$0<c< a^{-1}=c_\At^+\vee c_\At^-.$$
\newpage

\newpage

As we have mentioned, the aim of this stochastic process is just to give a simple way to outline the mixing phenomenon for the flat case. Similarly to the approach of Otto, while this random walk provides infinitely many trajectories $\theta_h=\{\theta_{s,j}^{(k)}\}$ starting from  \eqref{flatdis} (for different mixing speeds $0<\alpha<1$), the simulations evidence that $\theta_{h}\overset{*}{\rightharpoonup}\breve{\theta}_{\At,\alpha}$. In other words, when $h\approx 0$, although each simulation yields a different picture, at the macroscopic level we can not distinguish them. Moreover, $\breve{\theta}_{\At,\alpha}$ can be (almost) recovered from each experiment separately by averaging it over lines as in Remark~\ref{diagram} 
$$\tfrac{1}{2N+1}\sum_{|s|\leq N}\theta_{s,j}^{(k)}\underset{N\rightarrow\infty}{\longrightarrow} \breve{\theta}_j^{(k)},$$
due to the Central Limit Theorem.

\section{The function $\Theta_{\At}$}\label{sec:Otto}

Since the derivation of  \eqref{entropy}\eqref{conservationlaw} from \cite{O} involves some parameters and computations, we have considered appropriate to give a brief explanation of it in order to save time to the reader. In \cite{O} the phase ``$s$'' introduced by Otto takes values in $\{0,1\}$, while in this paper the phase $\theta$ takes values in $\{-1,1\}$. Both are related via: $s=0\leftrightarrow\theta=1$ and $s=1\leftrightarrow\theta=-1$. Thus, the density $\rho$ and the \textbf{mobility} $m=\mu^{-1}$ are described in terms of the phase $s$ as
\begin{equation}\label{IPMO:0}\tag{IPM0}
a(t,\x)
%=a^-s(t,\x)+a^+(1-s(t,\x))
=a^++(a^--a^+)s(t,\x),\quad a=\rho,m.
\end{equation}
After rescaling in time, Otto considered the (normalized) IPM system
\begin{align}
\partial_ts+\nabla\cdot(s\mathbf{v}) & = 0, \label{IPMO:1}\tag{IPM1}\\
\Div\mathbf{v} & = 0,\label{IPMO:2}\tag{IPM2}\\
\Curl((B^{-1}s+(1-s))\mathbf{v}-s\im) & =0, \label{IPMO:3}\tag{$\textrm{IPM3}^B$}
\end{align}
in $\R_+\times\Domain$, starting from the unstable planar phase $s_0=\tfrac{1-\theta_0}{2}$ \eqref{flat},
where $B$ is the \textbf{mobility quotient}
$$B=\frac{m^-}{m^+}=\frac{\mu^+}{\mu^-}
=\frac{1+\At}{1-\At}>0\quad
\leftrightarrow\quad \At=\frac{B-1}{B+1}\in(-1,1).$$
Thus, one can easily check that $(s,\mathbf{v})$ is a solution to $(\textrm{IPM}^B)$ if and only if $(\theta,\velocity)$ given by
$$\theta(t,\x)=1-2s(\alpha t,\x),
\quad\quad
\velocity(t,\x)=\alpha\mathbf{v}(\alpha t,\x),$$ 
with $\alpha=1+B^{-1}$,
solves $(\textrm{IPM}_\At)$. After the relaxation explained in Appendix \ref{sec:RW}, Otto obtained the entropy solution
$$S_{B}(t,\x)
=\left\lbrace
\begin{array}{cl}
	0, & \hspace{1.0cm}x_2>B t,\\
	\frac{B t-x_2}{B t+(B-1)x_2+\sqrt{B^2t(B t+(B-1)x_2)}}, & -t< x_2<B t,\\
	1, &  -t>x_2,
\end{array}\right.$$
of the scalar conservation law
$$
\partial_tS+\partial_{x_2}\left(\frac{S(1-S)}{S+B^{-1}(1-S)}\right)=0,\quad\quad
S|_{t=0}=s_0.
$$
Hence, since
$$
\alpha=1+B^{-1}=c_\At^-,\quad\quad
B\alpha=1+B=c_\At^+,
$$
the function $\Theta_{\At}(t,\x)=1-2S_B(\alpha t,\x)$ is the entropy solution of the scalar conservation law \eqref{conservationlaw}. Clearly $\Theta_{\At}(t,\x)=\pm 1$ in $\Omega_{\pm}=\{(x,t)\in\R_+\times\Domain\,:\,\pm x_2> c_\At^{\pm}t\}$. Inside the mixing zone $\Mixzone=\{(t,\x)\in\R_+\times\Domain\,:\,-c_\At^-t<x_2< c_\At^+t\}$, for $\At=0$ we have
$$\Theta_0(t,\x)=\frac{x_2}{2t},$$
and for $0<|\At|<1$ it is not difficult to check the following identities
\begin{align*}
\Theta_{\At}(t,\x)&=\frac{(x_2- t)+\sqrt{B t(t+\At x_2)}}{ t+\At  x_2+\sqrt{B t(t+\At x_2)}}\\
&=\frac{x_2+\At t}{ t+\At  x_2+\sqrt{(1-\At^2) t(t+\At x_2)}}\\
&=\frac{1}{\At}\left(1-\sqrt{\frac{(1-\At^2)t}{t+\At x_2}}\right).
\end{align*}

\begin{prop}\label{Oprop} For $\Domain=\R^2$, $\Theta_{\At}$ satisfies the following properties. At each time slice:
\begin{enumerate}[(i)]
	\item\label{Oprop:1} $\Theta_{\At}(t,\cdot)$ is continuous and smooth in $\Mixzone(t)$.
	\item\label{Oprop:2} $\Theta_{\At}(t,\cdot)$ is strictly $x_2$-increasing and concave (convex) for $\At>0$ ($\At<0$) in $\Mixzone(t)$. 
	\item\label{Oprop:3} $\Theta_{\At}(t,\x)=\Theta_{\At}(\tau,\frac{\tau}{t}\x)$ for all $\tau>0$ and $\x\in\R^2$.
	\item\label{Oprop:4} $\Theta_{-\At}(t,\x)=-\Theta_{\At}(t,-\x)$.
	\item\label{Oprop:5} For every $L=(l_1,l_2)\subset\alpha(-c_\At^-,c_\At^+)$,
	$$\expected{L}_{\At,\alpha}=\dashint_{L}\Theta_{\At}(\alpha,x_2)\dif x_2=
	\left\lbrace\begin{array}{cl}
	\displaystyle\frac{1}{\At}\left(1-\frac{2\sqrt{(1-\At^2)\alpha}}{\sqrt{\alpha+\At l_1}+\sqrt{\alpha+\At l_2}}\right), & \At\neq 0,\\[0.1cm]
	\displaystyle\frac{l_1+l_2}{4\alpha}, & \At=0.
	\end{array}\right.$$
\end{enumerate}
For $\Domain=(-1,1)^2$ see Section \ref{sec:transition}.
\end{prop}
\begin{proof} \ref{Oprop:1} is a straightforward computation. \ref{Oprop:2} is a consequence of
\begin{align*}
\partial_{x_2}\Theta_{\At}(t,\x)&=\tfrac{1}{2}\sqrt{(1-\At^2)t}(t+\At x_2)^{-\tfrac{3}{2}}>0,\\
\partial_{x_2}^2\Theta_{\At}(t,\x)&=-\tfrac{3}{4}\At\sqrt{(1-\At^2)t}(t+\At x_2)^{-\tfrac{5}{2}}.
\end{align*}
\ref{Oprop:3}\ref{Oprop:4} follow from \eqref{entropy}. \ref{Oprop:5} is due to, for $\At=0$
$$\int\Theta_{0}(\alpha,x_2)\dif x_2=\frac{x_2^2}{4\alpha},$$
and for $\At\neq 0$
$$\int\Theta_{\At}(\alpha,x_2)\dif x_2=
\frac{1}{\At^2}\left(\At x_2-2\sqrt{(1-\At^2)\alpha(\alpha+\At x_2)}\right).$$
\end{proof}

\begin{Rem}\label{uncernainity} To conclude we recall briefly the ``uncertainty principle'' presented in \cite{Degraded}. On the one hand, for $a=\rho,\mu$ given in terms of a $\Theta_{\At}$-mixing solution $\theta$ via (IPM0), the Lebesgue Differentiation Theorem implies
$$\lim_{\substack{\Mixzone(t)\supset R\downarrow \{\x_0\}  \\ R\textrm{ regular}}}\dashint_R a(t,\x)\dif x =a(t,\x_0),$$
for a.e.~$\x_0\in\Domain$ at each time slice $t\in\R_+$, where $a$ jumps unpredictably between $a^+$ and $a^-$ due to Thm.~\ref{HP}\ref{Mixingprop}. On the other hand, for every  rectangle $R=S\times tL\subset \Mixzone(t)$ either large or close enough to the (space-time) boundary of the mixing zone, we have
$$\dashint_R a(t,\x)\dif\x
\approx\tfrac{a^++a^-}{2}+\tfrac{a^+-a^-}{2}\expected{L}_{\At,\alpha},$$
at each time slice $t\in\R_+$, due to Thm.~\ref{HP}\ref{degraded}. In other words, either the position is localized $\{\x_0\}$ and so the phase is unpredictable or it is averaged in a suitable region $R$.
\end{Rem}

\subsection{Transition to the stable planar phase}\label{sec:transition}

In this section we describe $\Theta_{\At}$ in the confined domain $\Domain=(-1,1)^2$ once the mixing zone hits the lower or upper boundary. Immediately after the heavier fluid attains $x_2=-1$ ($c_{\At}^{-}t>1$) the bottom of the tank begins to be filled up with it and the phases
begin to separate\\
$$\Theta_{\At}(t,\x)
=\left\lbrace
\begin{array}{cl}
\frac{x_2+\At t}{ t+\At  x_2+\sqrt{(1-\At^2) t(t+\At x_2)}}, & d_{\At}^{-}(t)<x_2<0,\\[0.1cm]
+1, &  d_{\At}^{-}(t)>x_2,
\end{array}\right.$$
 
\noindent and the same happens once the lighter one attains $x_2=1$ ($c_{\At}^{+}t>1$)
$$\Theta_{\At}(t,\x)
=\left\lbrace
\begin{array}{cl}
-1, & \hspace{0.7cm}x_2>d_{\At}^{+}(t),\\[0.1cm]
\frac{x_2+\At t}{ t+\At  x_2+\sqrt{(1-\At^2) t(t+\At x_2)}}, & 0< x_2<d_{\At}^{+}(t),
\end{array}\right.$$

\noindent where $d_{\At}^{\pm}$ are the free boundaries, to be determined.

\noindent By taking $\breve{\velocity}_{\At}=0$ and $\breve{\aux}_{\At}$ as in \eqref{Msubsol} ($\alpha=1$), ($\mathbf{T}2$-$3_{\At}$) is automatically satisfied while ($\mathbf{T}1$) is equivalent to
\begin{equation}\label{eqh:1}
[\Theta_{\At}]_{\pm}\partial_td_{\At}^{\pm}=[\breve{\aux}_{\At}]_{\pm},
\end{equation}
where $[\cdot]_{\pm}$ denotes the jump discontinuity at $x_2=d_{\At}^{\pm}$ respectively. By writing $d_{\At}^{\pm}=\pm(1-f_{\At}^{\pm})$, \eqref{eqh:1} turns out to be a Cauchy problem for $f_{\At}^{\pm}$
\begin{equation}
\label{CP:f}
\begin{split}
\partial_tf_{\At}^{\pm}&=F_{\At}^{\pm}(t,f_{\At}^{\pm}),\\
f_{\At}^{\pm}|_{c_{\At}^{\pm}t=1}&=0,
\end{split}
\end{equation}
where 
$$F_{\At}^{\pm}(t,f)=\frac{1\mp\Theta_{\At}(t,\pm(1-f(t)))}{1-\Theta_{\At}(t,\pm(1-f(t)))\At}.$$
By the Picard-Lindel\"{o}f Theorem, there is a unique solution to \eqref{CP:f}. Furthermore, it is strictly increasing with $f_{\At}^{\pm}(t_{\At}^{\pm})=1$ ($d_{\At}^{\pm}(t_{\At}^{\pm})=0$) at some $1<c_{\At}^{\pm}t_{\At}^{\pm}<\infty$. Since ($\mathbf{T}1$) implies $\int\Theta_{\At}(t,\x)\dif x=0$ for all times, necessarily $t_{\At}^{\pm}=t_{\At}$. That is, the mixing zone collapses at this finite time $t_{\At}$ and the stable planar phase is reached. For $\At=0$ this is explicit
$$f_0(t)=1+2t-2\sqrt{2t},$$
for all $c_0^{-1}=\frac{1}{2}\leq t\leq 2=t_0.$

\section*{Acknowledgements}

The author thanks \'Angel Castro, Daniel Faraco and Sauli Lindberg for their valuable comments during the preparation of this work, and also thanks Lorena Romero for the Matlab simulations and Elena Mengual for her help with the GeoGebra pictures.
This work was partially supported by the Spanish Ministry of Economy through the ICMAT Severo Ochoa project SEV-2015-0554, 
the grant MTM2017-85934-C3-2-P (Spain) and the ERC grant  307179-GFTIPFD, ERC grant 834728-QUAMAP.

%\vfill
\centering

\begin{flushleft}
	\quad\\
	\address{Departamento de Matem\'aticas, Universidad Aut\'onoma de Madrid, E-28049 Madrid, Spain; Instituto de Ciencias Matem\'aticas (CSIC-UAM-UC3M-UCM), E-28049  Madrid, Spain.}
	\email{francisco.mengual@uam.es}
\end{flushleft}

\end{document}